\documentclass[hidelinks,onefignum,onetabnum]{siamart250211}

\usepackage{lipsum}
\usepackage{amsfonts}
\usepackage{graphicx}
\usepackage{epstopdf}
\usepackage{algorithmic}
\usepackage{amsmath}
\usepackage{amssymb}
\usepackage{cite}
\usepackage{xfrac}
\ifpdf
  \DeclareGraphicsExtensions{.eps,.pdf,.png,.jpg}
\else
  \DeclareGraphicsExtensions{.eps}
\fi

\newsiamremark{remark}{Remark}
\newsiamremark{hypothesis}{Hypothesis}
\crefname{hypothesis}{Hypothesis}{Hypotheses}
\newsiamthm{claim}{Claim}
\newsiamremark{fact}{Fact}
\crefname{fact}{Fact}{Facts}

\graphicspath{{./figs/}}

\headers{Theory of ADER-DG method for ODEs}{I. S. Popov}
\title{Theory and internal structure of ADER-DG method for ordinary differential equations\thanks{Submitted to the editors DATE.}}

\author{Ivan S. Popov\thanks{Department of Theoretical Physics, Dostoevsky Omsk State University, Omsk, Russia\\ (\email{diphosgen@mail.ru}, \email{popovis@omsu.ru}).}}

\usepackage{amsopn}

\newcommand{\corrtext}[1]{#1}

\ifpdf\hypersetup{
  pdftitle={Theory and internal structure of ADER-DG method for ordinary differential equations},
  pdfauthor={I. S. Popov}
}\fi

\begin{document}

\maketitle

\begin{abstract}
Investigation of the approximation properties, convergence, and stability of the ADER-DG method for solving an ODE system is carried out. \corrtext{The ADER-DG method is $A$- and $AN$-stable, $L$-stable, $B$- and $BN$-stable, and algebraically stable. Several other relations useful for an application and implementation of the ADER-DG method are proved. Applications of the ADER-DG method demonstrated compliance with the expected theoretical results.}
\end{abstract}

\section{Introduction}

In this paper, a study of the approximation properties, convergence and stability of arbitrary high order ADER (``Arbitrary-higher-order-DERivati\-ves'') discontinuous Galerkin (DG) method with local DG predictor for solving the initial value problem (IVP) for a system of ordinary differential equations (ODE) is carried out. The IVP for ODE system is chosen in the following classical form:
\begin{equation}\label{eq:ivp_ode_diff_src}
\frac{d\mathbf{u}}{dt} = \mathbf{F}\left(\mathbf{u},\, t\right),\quad t\in\Omega = \left\{t\, |\, t \in [t_{0},\ t_{f}]\right\},\quad \mathbf{u}(t_{0}) = \mathbf{u}_{0},
\end{equation}
where $\mathbf{u}: \Omega \rightarrow \mathcal{R}^{D}$ is a desired function, $\mathbf{F}: \mathcal{R}^{D}\times\Omega \rightarrow \mathcal{R}^{D}$ is a right side function, which is given. The initial condition $\mathbf{u}_{0}$ on the desired function $\mathbf{u}$ was chosen at the point $t_{0}$. The classical ODE theory shows that in case $\mathbf{F} \in \mathcal{C}_{1}(\mathcal{R}^{D}\times\Omega)$ the solution of the problem exists and is unique. The ODE systems of the second and higher orders, uniquely solvable with respect to higher derivatives, can be represented in the form of the first-order ODE system (\ref{eq:ivp_ode_diff_src}).

ODE systems are used in a wide range of fields of science and technology. There are now many numerical methods for solving ODE systems~\cite{Butcher_book_2016, Hairer_book_1, Hairer_book_2}. In recent years, there has been increased scientific interest in studying DG methods and the associated superconvergence in solving ODE systems, see books~\cite{Babuska_book_2001, Wahlbin_lectures_1995} and recent articles~\cite{dg_ivp_ode_4, dg_ivp_ode_5, dg_ivp_ode_6, dg_ivp_ode_1, dg_ivp_ode_2, dg_ivp_ode_3, ader_dg_ivp_ode}. DG methods were proposed in~\cite{lasl_rep_dg_1973} for the numerical description of neutron transport. Delfour \textit{et al} \cite{Delfour_1981} proposed a DG method for solving ODE systems, then a rigorous theory on this basis was developed in~\cite{Delfour_1986}. Cockburn, Shu \textit{et al.} in a series of works~\cite{Cockburn_base_1, Cockburn_base_2, Cockburn_base_3, Cockburn_base_4, Cockburn_base_5} created a rigorous mathematical basis for DG methods, which stimulated their further development and use for solving a wide class of problems.

The ADER-DG methods are DG methods based on the ADER paradigm proposed by Titarev and Toro in~\cite{ader_init_1, ader_init_2} in the context of finite volume methods. The modern version of ADER involves the use of a local DG predictor, which was proposed by Dumbser \textit{et al}~\cite{ader_stiff_1, ader_stiff_2}. The use of ADER in DG methods has allowed obtaining unrivaled results in the numerical solution of systems of partial differential equations (PDE)~\cite{PNPM_DG_2009, PNPM_DG_2010, ader_dg_ideal_flows, ader_dg_PNPM, ader_dg_seiemic_underwater, ader_dg_seiemic, ader_dg_wb_shwater_2022, dg_entropy, dg_entropy_add, ader_stiff_3, ader_stiff_4, ader_rev_2024, ader_dg_gr_z4_2024, ader_eff_blas}. Han Veiga \textit{et al}~\cite{dec_vs_ader_2021} and Micalizzi \textit{et al}~\cite{dec_vs_ader_2023} showed that the methods of the ADER family are significantly interconnected with numerical methods on the deferred correction (DeC) paradigm.

\corrtext{The use of the ADER-DG method for solving ODE systems was proposed in~\cite{ader_dg_ivp_ode, ader_improving_2024}, based on the ideas in~\cite{ader_stiff_1, PNPM_DG_2010}. ADER-DG has demonstrated high accuracy and stability in solving ODE systems, as well as systems of differential-algebraic equations~\cite{ader_dg_ivp_dae}, however, these results are only empirical, based on qualitative considerations, numerical calculations, and solving a set of test problems.}

\corrtext{Han Veiga \textit{et al}~\cite{ader_improving_2024} showed that the ADER-DG method is equivalent to a new implicit Runge-Kutta (RK) method, and rigorously proved that it has approximations of order $2N+1$ for the numerical solution at nodes and $N+1$ for the local solution in the space between nodes, if polynomials of degree $N$ are used to represent the solution. \"Offner \textit{et al}~\cite{ader_proofs_2025} rigorously proved the $A$-stability of the ADER-DG method and showed that the stability function has the form of a Pad\'{e} approximation $R_{N, N+1}(z)$ of $\exp(z)$. This concludes the existing rigorous analysis of the ADER-DG method.}

\corrtext{An important requirement for the applicability of the ADER-DG method to solving a wide class of problems (\ref{eq:ivp_ode_diff_src}) is the analog of the $A$-stability property for equations with variable coefficients --- $AN$-stability, and the properties of nonlinear stability --- $B$- and $BN$-stability~\cite{Butcher_book_2016, Hairer_book_1, Hairer_book_2}. The properties of $AN$-, $B$-, $BN$-stability and algebraic stability of the ADER-DG method are rigorously proven and investigated for ADER-DG method. The study is carried out for an arbitrary value of the degree of polynomials $N$. Some already known properties of approximation and stability are also proven by different tools compared to~\cite{ader_improving_2024, ader_proofs_2025}.}

\section{Formulation of ADER-DG method}

A description of the ADER-DG method for solving problem (\ref{eq:ivp_ode_diff_src}) was presented in~\cite{ader_dg_ivp_ode, ader_dg_ivp_dae}, therefore the description will be limited to what is necessary for this paper.

The numerical solution is found on a one-dimensional grid that discretizes the domain of definition $\Omega$ by a finite set $\{\Omega_{n}\}$ of non-overlapping discretization domains $\Omega_{n} = \left\{t\, |\, t \in [t_{n},\ t_{n+1}]\right\}$, where $t_{n}$ are the grid nodes, ${\Delta t}_{n} = t_{n+1} - t_{n}$ is the discretization step, and the discretization domains $\Omega_{n}$ themselves are called the space between the nodes. The numerical method is completely single-step, therefore the formula apparatus of the method itself is formulated only for one discretization domain $\Omega_{n}$, and allows the use of a variable discretization step ${\Delta t}_{n}$. Due to the high accuracy of the ADER-DG method even on coarse grids~\cite{ader_dg_ivp_ode, ader_dg_ivp_dae}, it is possible to allow the use of a grid containing only one discretization domain $\Omega_{1} = \Omega$.

The ADER-DG method is based on the integral form of the ODE system (\ref{eq:ivp_ode_diff_src}):
\begin{equation}\label{eq:ivp_ode_int_src}
\mathbf{u}_{n+1} = \mathbf{u}_{n} + \int\limits_{t_{n}}^{t_{n+1}} \mathbf{F}\left(\mathbf{u}(t),\, t\right) dt,
\end{equation}
where $\mathbf{u}_{n}$ and $\mathbf{u}_{n+1}$ are the numerical solution at grid nodes $t_{n}$ and $t_{n+1}$, respectively. In the integrand $\mathbf{F}$ the function $\mathbf{u} = \mathbf{u}(t)$, explicitly dependent on $t$, is defined in the space $\Omega_{n}$ between the grid nodes $t_{n}$ and $t_{n+1}$. 
To formulate the ADER-DG method and isolating its main structural components, the discretization domain $\Omega_{n}$ is mapped $t\mapsto\tau(t)$ onto the reference domain $\omega = \{\tau\, |\, \tau\in[0,\ 1]\}$:
\begin{equation}
t(\tau) = t_{n} + \tau\cdot{\Delta t}_{n}\in\Omega_{n},\quad \tau(t) = \frac{t - t_{n}}{{\Delta t}_{n}}\in\omega,
\end{equation}
where $\tau\in[0,\ 1]$ is the coordinate of the reference domain $\omega$. The ADER-DG numerical method is based on the use of a local solution $\mathbf{q}_{n}: \omega \rightarrow \mathcal{R}^{D}$, which is defined in the space between the grid nodes $\Omega_{n}$:
\begin{equation}
\mathbf{u} = \mathbf{u}(t(\tau)),\, t\in\Omega_{n}\ \mapsto\ \mathbf{q}_{n} = \mathbf{q}_{n}(\tau),
\end{equation}
and is used to calculate the integral in (\ref{eq:ivp_ode_int_src}). The local solution is also defined for the definition domain $\mathbf{u}_{L}: \Omega \rightarrow \mathcal{R}^{D}$ by a piecewise assembly of local solutions $\mathbf{q}_{n}$ for each $\Omega_{n}$:
\begin{equation}\label{eq:local_sol_assembly}
\mathbf{u}_{L}(t) = \sum\limits_{n} \chi_{n}(t)\cdot\mathbf{q}_{n}\left(\frac{t - t_{n}}{{\Delta t}_{n}}\right),
\end{equation}
where $\chi_{n}: \Omega\rightarrow\{0,\, 1\}$ is the indicator function of the discretization domain $\Omega_{n}\subseteq\Omega$.

The original version of the ADER numerical methods~\cite{ader_init_1, ader_init_2} chose a local solution as expansion according to the Cauchy-Kovalevskaya procedure, which in the case of ODE systems is close to the families of high-order Taylor numerical methods~\cite{ivp_ode_taylor_series_soft_2005, ivp_ode_taylor_series_2017}. The modern formulation of the ADER numerical methods~\cite{ader_stiff_1, ader_stiff_2} is based on the use of a DG predictor, in which the local solution is represented in the form of an expansion over a set $\{\varphi_{p}\}$ of basis functions $\varphi_{p}: \omega\rightarrow\mathcal{R}$ of the following form:
\begin{equation}\label{eq:local_sol}
\mathbf{q}_{n}(\tau) = \sum\limits_{p} \hat{\mathbf{q}}_{n, p}\varphi_{p}(\tau),
\end{equation}
and is found as a solution to the weak form of the ODE system (\ref{eq:ivp_ode_diff_src}) in $\Omega_{n}$:
\begin{equation}\label{eq:ivp_ode_weak}
\int\limits_{0}^{1}\varphi_{p}(\tau)\left[\frac{d\mathbf{q}_{n}(\tau)}{d\tau} - {\Delta t}_{n}\mathbf{F}(\mathbf{q}_{n}(\tau), t(\tau))\right] = 0,\quad
\mathbf{q}_{n}(0) = \mathbf{u}_{n}.
\end{equation}
The basis functions are selected in the form of Lagrange interpolation polynomials $\{\varphi_{p}\}_{p = 0}^{N}$ with nodal points at the roots of the shifted Legendre polynomials $\tilde{P}_{N+1}$:
\begin{equation}\label{eq:basis_funcs_def}
\varphi_{p}(\tau) = \sum\limits_{k = 0}^{N} \varphi_{p, k}\tau^{k} = \prod\limits_{k \neq p} \frac{\tau - \tau_{k}}{\tau_{p} - \tau_{k}},\quad
\varphi_{p}(\tau_{q}) = \delta_{pq},
\end{equation}
where $\{\tau_{k}\}_{k = 0}^{N}$ is an ascending order set of roots of the shifted Legendre polynomials $\tilde{P}_{N+1}$, $\{\varphi_{p, k}\}$ is a set of coefficients of the basis polynomials $\{\varphi_{p}\}$. The choice of the basis allowed to use the Gauss-Legendre (GL) quadrature formula:
\begin{equation}\label{eq:gl_rule}
\int\limits_{0}^{1} f(\tau)d\tau \approx \sum\limits_{p = 0}^{N} w_{p} f(\tau_{p}),\ 
w_{p} = \int\limits_{0}^{1} \varphi_{p}^{2}(\tau) d\tau = \int\limits_{0}^{1} \varphi_{p}(\tau) d\tau > 0,\ 
\sum\limits_{p = 0}^{N} w_{p} = 1,
\end{equation}
where $\{w_{p}\}$ is the weights, to calculate integrals (\ref{eq:ivp_ode_int_src}) and obtain point-wise evaluation
\begin{equation}
\mathbf{F}(\mathbf{q}_{n}(\tau), t(\tau)) \approx \sum\limits_{p = 0}^{N} \hat{\mathbf{F}}_{n, p}\varphi_{p}(\tau)\ \mapsto\ 
\sum\limits_{p = 0}^{N} \mathbf{F}(\hat{\mathbf{q}}_{n, p}, t(\tau_{n}))\varphi_{p}(\tau)\corrtext{,}
\end{equation}
\corrtext{which is based on the representation of the expansion coefficients $\hat{\mathbf{F}}_{n, p}$ in basis $\{\varphi_{p}\}$:
\begin{equation}\nonumber
\begin{split}
\hat{\mathbf{F}}_{n, p} &= \frac{1}{w_{p}} \int\limits_{0}^{1} \mathbf{F}(\mathbf{q}_{n}(\tau), t(\tau)) \varphi_{p}(\tau) d\tau \approx
\frac{1}{w_{p}} \sum\limits_{q = 0}^{N} w_{q} \mathbf{F}\left(\sum\limits_{r = 0}^{N} \hat{\mathbf{q}}_{n, r}\varphi_{r}(\tau_{q}), t(\tau_{q})\right) \varphi_{p}(\tau_{q})\\
&= \frac{1}{w_{p}} \sum\limits_{q = 0}^{N} w_{q} \mathbf{F}\left(\sum\limits_{r = 0}^{N} \hat{\mathbf{q}}_{n, r}\delta_{rq}, t(\tau_{q})\right) \delta_{pq} = 
\frac{w_{p}}{w_{p}} \mathbf{F}(\hat{\mathbf{q}}_{n, p}, t(\tau_{n})) = \mathbf{F}(\hat{\mathbf{q}}_{n, p}, t(\tau_{n})).
\end{split}
\end{equation}}
Evaluating the expression (\ref{eq:ivp_ode_weak}) leads to a system of nonlinear algebraic equations:
\begin{equation}\label{eq:ivp_ode_weak_rewr_dst}
\sum\limits_{q = 0}^{N}\Big[\kappa_{pq}\hat{\mathbf{q}}_{n, q} - {\Delta t}_{n} \mu_{pq}\mathbf{F}(\hat{\mathbf{q}}_{n, q}, t(\tau_{q}))\Big] = 
\smash[b]{\underline{\varphi}}_{p}\mathbf{u}_{n},
\end{equation}
where the coefficients can be collected into matrices $\kappa = \lbrace\kappa_{pq}\rbrace_{p,q}$, $\mu = \lbrace\mu_{pq}\rbrace_{p,q}$ and calculated by the following expressions:
\begin{equation}\label{eq:kappa_mu_defs}
\begin{split}
&\kappa_{pq} = \overline{\varphi}_{p}\overline{\varphi}_{q} - \int\limits_{0}^{1} \frac{d\varphi_{p}(\tau)}{d\tau}\varphi_{q}(\tau)d\tau
			= \smash[b]{\underline{\varphi}}_{p}\smash[b]{\underline{\varphi}}_{q} + \int\limits_{0}^{1} \varphi_{p}(\tau)\frac{d\varphi_{q}(\tau)}{d\tau}d\tau,\\
&\mu_{pq} = \int\limits_{0}^{1} \varphi_{p}(\tau)\varphi_{q}(\tau)d\tau \equiv w_{p} \delta_{pq},\quad
\smash[b]{\underline{\varphi}}_{p} = \varphi_{p}(0),\quad
\overline{\varphi}_{p} = \varphi_{p}(1),
\end{split}
\end{equation}
where the last expression explicitly takes into account the orthogonality of the basis $\{\varphi_{p}\}$ and the expression for the weights $w_{p}$ of the quadrature formula (\ref{eq:gl_rule}). Multiplication by $\kappa^{-1}$ leads to the system of equations of the local DG predictor:
\begin{equation}\label{eq:snae_lstdg}
\hat{\mathbf{q}}_{n, p} = \mathbf{u}_{n} + {\Delta t}_{n}\sum\limits_{q = 0}^{N} a_{pq}\mathbf{F}(\hat{\mathbf{q}}_{n, q}, t(\tau_{q})),
\end{equation}
where $a_{pq}$ is elements of the matrix $a = \kappa^{-1}\mu$, and Lemma~\ref{lemma:ader_int_prop} is taken into account.
\begin{lemma}\label{lemma:ader_int_prop}
The parameters of the ADER-DG method satisfy the relations:
\begin{align}
&\smash[b]{\underline{\varphi}}_{p} = \sum\limits_{q = 0}^{N}\kappa_{pq},\quad
\sum\limits_{q = 0}^{N}\left[\kappa^{-1}\right]_{pq}\smash[b]{\underline{\varphi}}_{q} = 1,\label{eq:ader_int_prop}\\
&\overline{\varphi}_{q} = \sum\limits_{p = 0}^{N}\kappa_{pq},\quad
\sum\limits_{p = 0}^{N}\left[\kappa^{-1}\right]_{pq} \overline{\varphi}_{p} = 1,\quad
\sum\limits_{p = 0}^{N}a_{pq} \overline{\varphi}_{p} = w_{q},\\
&\mathrm{tr}(\kappa) = \sum\limits_{p = 0}^{N}\kappa_{pp} = 
\frac{1}{2}\sum\limits_{p = 0}^{N}\left[\smash[b]{\underline{\varphi}}_{p}^{2} + \overline{\varphi}_{p}^{2}\right] = 
\sum\limits_{p = 0}^{N}\smash[b]{\underline{\varphi}}_{p}^{2} = \sum\limits_{p = 0}^{N}\overline{\varphi}_{p}^{2}.
\end{align}
\end{lemma}
\begin{proof}
These relationships are proved by direct expansion of the expressions:
\begin{equation}\nonumber
\begin{split}
\sum\limits_{q = 0}^{N}\kappa_{pq} & = \Bigg|\sum\limits_{p = 0}^{N}\varphi_{p}(\tau) \equiv 1\Bigg|
	= \smash[b]{\underline{\varphi}}_{p} \sum\limits_{p = 0}^{N}\smash[b]{\underline{\varphi}}_{q} +
	\int\limits_{0}^{1} \varphi_{p}(\tau) \frac{d}{d\tau}\left[\sum\limits_{q = 0}^{N}\varphi_{q}(\tau)\right]d\tau = \smash[b]{\underline{\varphi}}_{p},\\
\sum\limits_{p = 0}^{N}\kappa_{pq} & = \overline{\varphi}_{q} \sum\limits_{p = 0}^{N}\overline{\varphi}_{p} - 
	\int\limits_{0}^{1} \varphi_{q}(\tau) \frac{d}{d\tau}\left[\sum\limits_{p = 0}^{N}\varphi_{p}(\tau)\right]d\tau = \overline{\varphi}_{p},
\end{split}
\end{equation}
\begin{equation}\nonumber
\begin{split}
\sum\limits_{q = 0}^{N}\kappa_{pq} & = \smash[b]{\underline{\varphi}}_{p},\quad \Leftrightarrow\quad
\sum\limits_{q = 0}^{N}\left[\kappa^{-1}\right]_{pq}\smash[b]{\underline{\varphi}}_{q} = 1,\\
\sum\limits_{p = 0}^{N}\kappa_{pq} & = \overline{\varphi}_{q},\quad \Leftrightarrow \quad
\sum\limits_{p = 0}^{N}\left[\kappa^{-1}\right]_{pq} \overline{\varphi}_{p} = 1,\quad \Leftrightarrow\quad
\sum\limits_{p = 0}^{N}a_{pq} \overline{\varphi}_{p} = w_{q},
\end{split}
\end{equation}
\begin{equation}\nonumber
\begin{split}
\mathrm{tr}(\kappa) &= \sum\limits_{p = 0}^{N}\kappa_{pp} =
\sum\limits_{p = 0}^{N}\left[\overline{\varphi}_{p}\overline{\varphi}_{p} - \int\limits_{0}^{1} \frac{d\varphi_{p}(\tau)}{d\tau}\varphi_{p}(\tau)d\tau\right]\\
& = \sum\limits_{p = 0}^{N}\left[
	\overline{\varphi}_{p}\overline{\varphi}_{p} - 
	\frac{1}{2}\left(\overline{\varphi}_{p}\overline{\varphi}_{p} - \smash[b]{\underline{\varphi}}_{p}\smash[b]{\underline{\varphi}}_{p}\right)
\right] = \frac{1}{2}\sum\limits_{p = 0}^{N}\left[\smash[b]{\underline{\varphi}}_{p}^{2} + \overline{\varphi}_{p}^{2}\right] =
\sum\limits_{p = 0}^{N}\smash[b]{\underline{\varphi}}_{p}^{2} = \sum\limits_{p = 0}^{N}\overline{\varphi}_{p}^{2},
\end{split}
\end{equation}
where in the last two identities the symmetry property $\tau_{p} = 1 - \tau_{N - p}$ of the zeros of shifted Legendre polynomials $\tilde{P}_{N+1}$ was used, resulting in $\varphi_{p}(\tau) = \varphi_{N - p}(1 - \tau)$, in particular, $\overline{\varphi}_{p} = \smash[b]{\underline{\varphi}}_{N - p}$, so the sums are the same. 
\end{proof}
\begin{remark}
\corrtext{Formula (\ref{eq:ader_int_prop}) is proved in~\cite{ader_improving_2024}, Proposition~2.4.}
\end{remark}
\begin{remark}
It can be noted that the second relation in (\ref{eq:ader_int_prop}) can be proved heuristically: in the case of $\mathbf{F} \equiv 0$, the solution $\mathbf{q}_{n}$ to the problem (\ref{eq:ivp_ode_weak}) is a constant, but the initial condition $\mathbf{q}_{n}(0) = \mathbf{u}_{n}$ assumes $\mathbf{q}_{n}(\tau) \equiv \mathbf{u}_{n}$, from which the relation (\ref{eq:ader_int_prop}) to be proven directly. Otherwise $\mathbf{F} \not\equiv 0$, the values $\kappa_{pq}$ and $\smash[b]{\underline{\varphi}}_{p}$ do not depend on $\mathbf{F}$ by construction, which reduces the proof to the case $\mathbf{F} \equiv 0$.
\end{remark}
\begin{remark}
The relations presented in the Lemma~\ref{lemma:ader_int_prop} do not depend on the explicit choice of Legendre polynomials $\tilde{P}_{N+1}$ and can also be used with other polynomials, such as the left and right Radau polynomials. However, the last transformation of the symmetrized sum into individual sums requires symmetry of the roots $\{\tau_{p}\}$.
\end{remark}
\begin{remark}
These relations can be useful not only in proving theorems, but also in software development, especially in testing and formulating assertions, or simplifying expressions, such as making the factor around $\mathbf{u}_{n}$ disappear in (\ref{eq:snae_lstdg}).
\end{remark}
The resulting system of equations (\ref{eq:snae_lstdg}) is generally nonlinear and can be solved by Newton's or Picard's iteration methods, as well as other methods for solving systems of nonlinear algebraic equations. The existence and uniqueness of a solution to the system of equations (\ref{eq:snae_lstdg}), which in this case can be obtained by iteration, is guaranteed under the conditions presented in Remark~\ref{remark:snae_lstdg_sol_exist_uniq}.

The expression (\ref{eq:ivp_ode_int_src}) of the solution $\mathbf{u}_{n+1}$ at the grid node $t_{n+1}$ as a result of using the quadrature formula (\ref{eq:gl_rule}) takes the form:
\begin{equation}\label{eq:ader_dg_node_sol}
\mathbf{u}_{n+1} = \mathbf{u}_{n} + {\Delta t}_{n}\sum\limits_{p = 0}^{N} w_{p} \mathbf{F}\left(\hat{\mathbf{q}}_{n, p}, t(\tau_{p})\right).
\end{equation}
The expressions (\ref{eq:snae_lstdg}), (\ref{eq:ader_dg_node_sol}) together constitute the ADER-DG method for solving the IVP for an ODE system (\ref{eq:ivp_ode_diff_src}).

The solution $\mathbf{u}_{n+1}$ at the grid node $t_{n+1}$ and the local solution $\mathbf{q}_{n}(1)$ coincide:
\begin{equation}\label{eq:ader_dg_node_sol_by_q}
\mathbf{q}_{n}(1)
= \sum\limits_{p = 0}^{N} \hat{\mathbf{q}}_{n, p} \overline{\varphi}_{p}
= \mathbf{u}_{n} + {\Delta t}_{n}\sum\limits_{q = 0}^{N}\left[
	\sum\limits_{p = 0}^{N} a_{pq} \overline{\varphi}_{p}
\right]\mathbf{F}(\hat{\mathbf{q}}_{n, q}, t(\tau_{q})) \equiv \mathbf{u}_{n+1},
\end{equation}
which is associated with a property of the parameters expressed by \corrtext{(\ref{eq:ader_int_prop})}.

\section{ADER-DG method as Runge-Kutta method}

The use of substitution
\begin{equation}\label{eq:k_to_q}
\hat{\mathbf{q}}_{n, p} = \mathbf{u}_{n} + {\Delta t}_{n}\sum\limits_{p = 0}^{N} a_{pq} \mathbf{k}_{q},
\end{equation}
known in the theory of Runge-Kutta (RK) methods~\cite{Butcher_book_2016, Hairer_book_1, Hairer_book_2, Dekker_Verwer_1984}, in expressions (\ref{eq:snae_lstdg}), (\ref{eq:ader_dg_node_sol}) allows to reformulate the ADER-DG method as the $s$-stage implicit RK method with $s = N+1$:
\begin{equation}\label{eq:ader_dg_rkm}
\begin{split}
\mathbf{k}_{p} = \mathbf{F}\left(
	\mathbf{u}_{n} + {\Delta t}_{n}\sum\limits_{p = 0}^{N} a_{pq} \mathbf{k}_{q},\ 
	t_{n} + \tau_{p}\cdot{\Delta t}_{n}
\right),\quad
\mathbf{u}_{n+1} = \mathbf{u}_{n} + {\Delta t}_{n}\sum\limits_{p = 0}^{N} w_{p} \mathbf{k}_{p},
\end{split}
\end{equation}
which in terms of weights $\{w_{p}\}$ and nodes $\{\tau_{p}\}$ corresponds to the implicit Gauss-Legendre (GL) method. The coefficients $a_{pq}$ do not coincide with any RK method known (at least to the author). The set of coefficients $\{a_{pq},\, w_{p},\, \tau_{p}\}$ can be exactly calculated in radicals for $N \leq 8$. In cases $N = 1$, $2$, $3$, the Butcher tables take the form:
\begin{equation}\label{eq:butcher_table}
\begin{split}
&\begin{array}{c|cc}
\frac{1}{2}-\frac{\sqrt{3}}{6} & \frac{1}{3} & -\frac{1-\sqrt{3}}{2}\\
\frac{1}{2}+\frac{\sqrt{3}}{6} & \frac{1-\sqrt{3}}{2} & \frac{1}{3}\\
\hline
& \frac{1}{2} & \frac{1}{2}
\end{array}\qquad
\begin{array}{c|ccc}
\frac{1}{2}-\frac{\sqrt{15}}{10}& \frac{29}{180} & \frac{8}{45}-\frac{\sqrt{15}}{15} & \frac{29}{180}-\frac{\sqrt{15}}{30}\\
\frac{1}{2} & \frac{1}{9}+\frac{\sqrt{15}}{24} & \frac{5}{18} & \frac{1}{9}-\frac{\sqrt{15}}{24}\\
\frac{1}{2}+\frac{\sqrt{15}}{10} & \frac{29}{180}+\frac{\sqrt{15}}{30} & \frac{8}{45}+\frac{\sqrt{15}}{15} & \frac{29}{180}\\[1mm]
\hline
& \frac{5}{18} & \frac{4}{9} & \frac{5}{18}
\end{array}\\
&\begin{array}{c|c}
\hspace{-4mm}
\begin{array}{c}
\frac{1}{2}-\frac{\sqrt{525+70 \sqrt{30}}}{70}\\
\frac{1}{2}-\frac{\sqrt{525-70 \sqrt{30}}}{70}\\
\frac{1}{2}+\frac{\sqrt{525-70 \sqrt{30}}}{70}\\
\frac{1}{2}+\frac{\sqrt{525+70 \sqrt{30}}}{70}
\end{array}\hspace{-2mm}&\hspace{-2mm}
\left[
\begin{array}{cccc}
\hspace{-2mm}\frac{5 \sqrt{30}}{72}+\frac{19}{24} & \hspace{-2mm}-\frac{43}{60}+\frac{17 \sqrt{30}}{72} & \hspace{-2mm}-\frac{13}{40} & \hspace{-2mm}\frac{13}{28}-\frac{53 \sqrt{30}}{630}\\
\hspace{-2mm}-\frac{43}{60}-\frac{17 \sqrt{30}}{72} & \hspace{-2mm}-\frac{5 \sqrt{30}}{72}+\frac{19}{24} & \hspace{-2mm}\frac{13}{28}+\frac{53 \sqrt{30}}{630} & \hspace{-2mm}-\frac{13}{40}\\
\hspace{-2mm}\frac{47}{40} & \hspace{-2mm}-\frac{27}{28}-\frac{47 \sqrt{30}}{420} & \hspace{-2mm}-\frac{5 \sqrt{30}}{72}+\frac{19}{24} & \hspace{-2mm}-\frac{43}{60}+\frac{17 \sqrt{30}}{72}\\
\hspace{-2mm}-\frac{27}{28}+\frac{47 \sqrt{30}}{420} & \hspace{-2mm}\frac{47}{40} & \hspace{-2mm}-\frac{43}{60}-\frac{17 \sqrt{30}}{72} & \hspace{-2mm}\frac{5 \sqrt{30}}{72}+\frac{19}{24}\hspace{-2mm}
\end{array}
\right]^{-1}\\
&\cdot\,\mathrm{diag}\left(
-\frac{\sqrt{30}}{72}+\frac{1}{4},\
\frac{\sqrt{30}}{72}+\frac{1}{4},\
\frac{\sqrt{30}}{72}+\frac{1}{4},\
-\frac{\sqrt{30}}{72}+\frac{1}{4} 
\right)\\
\hline
&\begin{array}{cccc}
-\frac{\sqrt{30}}{72}+\frac{1}{4}&
\frac{\sqrt{30}}{72}+\frac{1}{4}&
\frac{\sqrt{30}}{72}+\frac{1}{4}&
-\frac{\sqrt{30}}{72}+\frac{1}{4} 
\end{array}
\end{array}
\end{split}
\end{equation}
but already in the case $N = 3$, the expressions for the coefficients $a_{pq}$ become large, therefore in the case of $N \geq 3$, it is more convenient to choose the coefficients $a_{pq}$ numerically \corrtext{(see Supplementary materials~SM1).}
\begin{definition}[\!\!\cite{ader_proofs_2025}]
ADER-IWF-RK-GLG method is the implicit RK method (\ref{eq:ader_dg_rkm}) based on ADER-DG method (\ref{eq:snae_lstdg}), (\ref{eq:ader_dg_node_sol}) with nodal basis (\ref{eq:basis_funcs_def}) in roots of Legendre polynomial $\tilde{P}_{N+1}$ and using the GL quadrature formula (\ref{eq:gl_rule}).
\end{definition}
\begin{remark}
Calculation of Butcher tables $\{a_{pq},\, w_{p},\, \tau_{p}\}$ for left and right Radau polynomials showed that the results are completely consistent with the Butcher tables of Radau IA and IIA methods, respectively. Therefore, these results presented in the classical literature~\cite{Butcher_book_2016, Hairer_book_2, Dekker_Verwer_1984} are expected, and not considered in this paper. \corrtext{Using Gauss-Lobatto nodes in quadrature leads to RK method Lobatto IIIC~\cite{ader_improving_2024}.}
\end{remark}
\begin{remark}\label{remark:snae_lstdg_sol_exist_uniq}
The ADER-DG method (\ref{eq:snae_lstdg}), (\ref{eq:ader_dg_node_sol}) is equivalent to the implicit RK method, so to determine the existence and uniqueness of a solution $\{\hat{\mathbf{q}}_{p}\}$ one can use Theorem 7.2 in~\cite{Hairer_book_1}, according to which, if $\mathbf{F}(\mathbf{u}, t)\in\mathcal{C}_{0}(\Xi)$ satisfying the Lipschitz condition with constant $C > 0$ in neighborhood $\Xi\subseteq\mathcal{R}^{D}\times\Omega$ of the $(\mathbf{u}_{n}, t_{n})$, and if
\begin{equation}
{\Delta t}_{n} \cdot C \max\limits_{p}\sum\limits_{q = 0}^{N} \left|a_{pq}\right| < 1,
\end{equation}
then the system of equations (\ref{eq:snae_lstdg}) has a unique solution $\{\hat{\mathbf{q}}_{p}\}$, which can be obtained by iteration. If also $\mathbf{F}(\mathbf{u}, t)\in\mathcal{C}_{p}(\Xi\times\Omega)\subset\mathcal{C}_{0}(\Xi\times\Omega)$, then $\{\hat{\mathbf{q}}_{q}({\Delta t}_{n})\in\mathcal{C}_{p}(\mathcal{R})\}_{q = 0}^{N}$.
\end{remark}

\section{Superconvergence and orders}

The determination of the approximation order of the ADER-DG method can now be carried out on the basis of the well-known theory of RK methods~\cite{Butcher_book_2016, Hairer_book_1, Dekker_Verwer_1984}. In this paper, simplifying conditions (see~\cite{Dekker_Verwer_1984}) are used, which are formulated in terms of the parameters of the ADER-DG method:
\begin{equation}\label{eq:simp_cond}
\begin{split}
B(L): &\ \sum\limits_{q = 0}^{N} w_{q}\tau_{q}^{r} = \frac{1}{r+1},\quad 0 \leq r < L,\\
C(L): &\ \sum\limits_{q = 0}^{N} a_{pq} \tau_{q}^{r} = \frac{\tau_{p}^{r+1}}{r+1},\quad 0 \leq p \leq N,\ 0 \leq r < L,\\
D(L): &\ \sum\limits_{q = 0}^{N} w_{q} a_{qp} \tau_{q}^{r} = \frac{w_{p}}{r+1}\left(1 - \tau_{p}^{r+1}\right),\quad 0 \leq p \leq N,\ 0 \leq r < L.
\end{split}
\end{equation}
The nodes and weights of the quadrature formula (\ref{eq:gl_rule}) ensure (see also~\cite{Dekker_Verwer_1984}) that condition $B(2N+2)$ is satisfied, therefore, to determine the approximation order, it is necessary to study conditions $C(L)$ and $D(L)$.

The conditions (\ref{eq:simp_cond}) were verified numerically using module \texttt{mpmath} (with internal use of module \texttt{gmpy2}, based on the  \texttt{MPFR} numerical library) of the \texttt{python} programming language, with an arithmetic precision constant \texttt{mpmath.mp.dps = 1000}. Numerical investigation of cases $N = 1, \ldots, 75$ showed that conditions $C(N)$, $D(N)$ are satisfied with precision $10^{-978}$--$10^{-1006}$. However, conditions $C(N+1)$, $D(N+1)$, which are inherent in the classical GL method are not satisfied, but the deviations decrease with $N$, reaching $10^{-48}$--$10^{-49}$ for $C(N+1)$ and $10^{-50}$--$10^{-52}$ for $D(N+1)$ in case $N = 75$. Based on expected result, a rigorous proof is created, presented by Theorem~\ref{theorem:order_rk}. In the classical GL method, the proof is based on the property $a_{pq} = \int_{0}^{\tau_{p}} \varphi_{q}(\tau) d\tau$, however, in this case, it is obviously not satisfied.
\begin{lemma}\label{lemma:p_ode_exact}
The local solution $q_{n}$ (\ref{eq:local_sol}) is the exact solution of the ODE
\begin{equation}\label{eq:p_ode}
\frac{du}{dt} = f(t),\quad f(t)\in\mathcal{P}_{L}(\mathcal{R}),\quad L < N.
\end{equation}
\end{lemma}
\begin{proof}
The exact solution $u\in\mathcal{P}_{L+1}(\mathcal{R})$ and can therefore be represented exactly as $q_{n} \in \mathrm{span}(\{\varphi_{p}\}_{p = 0}^{N})$ (\ref{eq:local_sol}) if and only if $L < N$. The local solution $q_{n}$ (\ref{eq:local_sol}) is the solution of the weak form (\ref{eq:ivp_ode_weak}) of the ODE (\ref{eq:p_ode}). However, $f(t)\in\mathcal{P}_{L}(\mathcal{R})\subset\mathcal{C}_{1}(\mathcal{R})$, and $\{\varphi_{p}\}_{p = 0}^{N}$ is a complete set of functions in $\mathcal{P}_{L+1}(\mathcal{R})\subseteq\mathcal{P}_{N}(\mathcal{R})$. Therefore, the solution of the weak form coincides with the solution of the original ODE (\ref{eq:p_ode}). The provable statement follows directly from this.
\end{proof}
\begin{remark}\label{lemma:p_ode_exact:remark:radau_r}
Nowhere in Lemma~\ref{lemma:p_ode_exact} is it explicitly used that the nodes $\{\tau_{p}\}_{p = 0}^{N}$ must be chosen at the roots of Legendre polynomials. Therefore, the Lemma~\ref{lemma:p_ode_exact} is satisfied in the case of left Radau polynomials. In the case of right Radau polynomials, the right node $\tau_{N} = 1$, so this Lemma can be weakened. As a consequence of (\ref{eq:ader_dg_node_sol_by_q}), the local solution at $q_{n}(1)$ coincides with the solution $u_{n+1}$ at the node $t_{n+1}$ obtained by the Gauss-Radau quadrature formula, similar to (\ref{eq:gl_rule}), that is exact for $\mathcal{P}_{2N}(\mathcal{R})$, so the local solution will be exact for an ODE (\ref{eq:p_ode}) with $f(t)\in\mathcal{P}_{N}(\mathcal{R})$, i.e. $L \leq N$.
\end{remark}
\begin{remark}\label{lemma:p_ode_exact:remark:radau_l}
Following the logic of the previous Remark~\ref{lemma:p_ode_exact:remark:radau_r}, it may seem that in the case of left Radau polynomials the weakening of the Lemma~\ref{lemma:p_ode_exact} should also be fulfilled, since in this case there is always the left node $\tau_{0} = 0$, and the initial condition $\mathbf{q}_{n}(0) = \mathbf{u}_{n}$ is defined. However, the initial condition is defined precisely in the sense of the weak form (\ref{eq:ivp_ode_weak}) of the ODE system (\ref{eq:ivp_ode_diff_src}), and is used in the integration by parts in the derivation of (\ref{eq:ivp_ode_weak_rewr_dst}) from (\ref{eq:ivp_ode_weak}). A direct calculation shows:
\begin{equation}
\begin{split}
&\mathbf{q}_{n}(0) = \sum\limits_{p = 0}^{N} \hat{\mathbf{q}}_{n, p} \smash[b]{\underline{\varphi}}_{p} =
\mathbf{u}_{n} + {\Delta t}_{n}\sum\limits_{q = 0}^{N}\left[\sum\limits_{p = 0}^{N} a_{pq} \smash[b]{\underline{\varphi}}_{p}\right]\mathbf{F}(\hat{\mathbf{q}}_{n, q}, t(\tau_{q})),\\
&\tilde{w}_{q} = \sum\limits_{p = 0}^{N} a_{pq} \smash[b]{\underline{\varphi}}_{p} 
= \sum\limits_{p = 0}^{N}\sum\limits_{r = 0}^{N} \left[\kappa^{-1}\right]_{pr} w_{r} \delta_{rq} \smash[b]{\underline{\varphi}}_{p}
= w_{q} \sum\limits_{p = 0}^{N} \left[\kappa^{-1}\right]_{pq} \smash[b]{\underline{\varphi}}_{p},
\end{split}
\end{equation}
where vector $\tilde{\boldsymbol{w}} = \lbrace\tilde{w}_{q}\rbrace_{q} \neq 0$ due to $\smash[b]{\underline{\boldsymbol{\varphi}}} = \lbrace\smash[b]{\underline{\varphi}}_{p}\rbrace_{p} \neq 0$ and non-degeneracy of matrix $\kappa^{-1}$, therefore $\mathbf{q}_{n}(0) \neq \mathbf{u}_{n}$ in the general case $\mathbf{F} \not\equiv 0$. Therefore, in the case of left Radau polynomials, weakening Lemma~\ref{lemma:p_ode_exact} is not admissible.
\end{remark}
\begin{remark}\label{lemma:p_ode_exact:remark:gl_classic}
In the case of classical GL methods, the Lemma~\ref{lemma:p_ode_exact} can also be weakened to the condition $f(t)\in\mathcal{P}_{N}(\mathcal{R})$, due to condition $a_{pq} = \int_{0}^{\tau_{p}} \varphi_{q}(\tau) d\tau$.
\end{remark}
\begin{lemma}[\corrtext{\!\!\cite{ader_improving_2024}, Lemma 3.5}]\label{lemma:cn_cond}
$C(N)$ (\ref{eq:simp_cond}) is satisfied for the ADER-IWF-RK-GLG method.
\end{lemma}
\begin{proof}
The proof is based on the meaning of the simplifying condition $C(L)$ (see~\cite{Dekker_Verwer_1984}) and Lemma~\ref{lemma:p_ode_exact}. The simplifying condition $C(L)$ means that $q_{n, p}$ are the values of exact solution of the ODE (\ref{eq:p_ode}) at quadrature points $t_{n} + \tau_{p} {\Delta t}_{n}$ in discretization domain $\Omega_{n}$. Note that if and only if the $\mathbf{q}_{n, p}$ are exact solutions at quadrature points and $N > L$, then the local solution $\mathbf{q}$ is an exact solution. Lemma~\ref{lemma:p_ode_exact} shows that $\mathbf{q}$ is exact solution only in case $L < N$, and in case $L \geq N$ the local solution $\mathbf{q}$ is not an exact solution. Therefore, in case $L < N$, the values $\mathbf{q}_{n, p}$ are exact solutions at quadrature points, and in case $L \geq N$ they are not. This proves $C(N)$.
\end{proof}
\begin{proof}[\protect{Alternative proof of Lemma~\ref{lemma:cn_cond}}]
The proof of this lemma is based on the transformation of the simplifying condition $C(N)$:
\begin{equation}
\sum\limits_{q = 0}^{N} a_{pq} \tau_{q}^{r} = \frac{\tau_{p}^{r+1}}{r+1}, \quad 0 \leq p \leq N,\ 0 \leq r \leq N-1.
\end{equation}
Substituting the parameters of ADER-DG method into $C(N)$ leads to the expression:
\begin{equation}
w_{p}\tau_{p}^{r} = \sum\limits_{q = 0}^{N} \kappa_{pq} \frac{\tau_{q}^{r+1}}{r+1}
=\sum\limits_{q = 0}^{N} \left[
	\smash[b]{\underline{\varphi}}_{p}\smash[b]{\underline{\varphi}}_{q} + \int\limits_{0}^{1} \varphi_{p}(\tau)\frac{d\varphi_{q}(\tau)}{d\tau}d\tau
\right] \frac{\tau_{q}^{r+1}}{r+1},
\end{equation}
where for each term the following transformations are performed:
\begin{equation}
\begin{split}
&\sum\limits_{q = 0}^{N} \frac{\tau_{q}^{r+1}}{r+1}\smash[b]{\underline{\varphi}}_{p}\smash[b]{\underline{\varphi}}_{q} = 
\frac{\smash[b]{\underline{\varphi}}_{p}}{r+1} \sum\limits_{q = 0}^{N} \tau_{q}^{r+1}\smash[b]{\underline{\varphi}}_{q} = \frac{\smash[b]{\underline{\varphi}}_{p}}{r+1}\cdot(0)^{r+1} = 0,
\end{split}
\end{equation}
\begin{equation}\nonumber
\begin{split}
&\sum\limits_{q = 0}^{N} \frac{\tau_{q}^{r+1}}{r+1}\int\limits_{0}^{1} \varphi_{p}(\tau)\frac{d\varphi_{q}(\tau)}{d\tau}d\tau =
\int\limits_{0}^{1} \frac{\varphi_{p}(\tau)}{r+1} \frac{d}{d\tau}\left[\sum\limits_{q = 0}^{N}\tau_{q}^{r+1}\varphi_{q}(\tau)\right]d\tau\\
&= \int\limits_{0}^{1} \frac{\varphi_{p}(\tau)}{r+1} \frac{d\left(\tau^{r+1}\right)}{d\tau}d\tau =
\int\limits_{0}^{1} \tau^{r}\varphi_{p}(\tau) d\tau = w_{p}\tau_{p}^{r},
\end{split}
\end{equation}
\corrtext{where the power function is represent as a basis expansion in $\{\varphi_{p}\}$.}
Substituting the transformed terms into the original expression turns it into an identity, which proves the fulfillment of the simplifying condition $C(N)$.
\end{proof}
\begin{remark}
In the case of right Radau polynomials, following Remark~\ref{lemma:p_ode_exact:remark:radau_r}, the Lemma~\ref{lemma:p_ode_exact} can be weakened, which leads to the fulfillment of the simplifying condition $C(N+1)$, which is known~\cite{Hairer_book_1, Dekker_Verwer_1984}.
\end{remark}
\begin{remark}
The $r = 0$ in $C(N)$ corresponds to the relation: $\sum_{q = 0}^{N} a_{pq} = \tau_{p}$.
\end{remark}
\begin{lemma}[\corrtext{\!\!\cite{ader_improving_2024}, Lemma 3.6}]\label{lemma:dn_cond}
$D(N)$ (\ref{eq:simp_cond}) is satisfied for the ADER-IWF-RK-GLG method.
\end{lemma}
\begin{proof}
The proof of this lemma is based on the transformation of the simplifying condition $D(N)$:
\begin{equation}
\sum\limits_{q = 0}^{N} w_{q} a_{qp} \tau_{q}^{r} = \frac{w_{p}}{r+1}\left(1 - \tau_{p}^{r+1}\right),\quad 0 \leq p \leq N,\ 0 \leq r \leq N-1.
\end{equation}
Substituting the parameters of ADER-DG method (\ref{eq:snae_lstdg}), (\ref{eq:ader_dg_node_sol}) into $D(N)$ leads to the following expression:
\begin{equation}
\begin{split}
w_{p}\tau_{p}^{r} = \sum\limits_{q = 0}^{N} \frac{\kappa_{qp}}{w_{q}}\frac{w_{q}}{r+1}\left(1 - \tau_{q}^{r+1}\right)
= \sum\limits_{q = 0}^{N}\hspace{-1mm}\left[
	\smash[b]{\underline{\varphi}}_{q}\smash[b]{\underline{\varphi}}_{p} + \int\limits_{0}^{1} \varphi_{q}(\tau)\frac{d\varphi_{p}(\tau)}{d\tau}d\tau
\right]\hspace{-1mm}\frac{1 - \tau_{q}^{r+1}}{r+1},
\end{split}
\end{equation}
where for each term the following transformations are performed:
\begin{equation}
\begin{split}
&\sum\limits_{q = 0}^{N} \frac{\smash[b]{\underline{\varphi}}_{q}\smash[b]{\underline{\varphi}}_{p}}{r+1} =
\frac{\smash[b]{\underline{\varphi}}_{p}}{r+1} \sum\limits_{q = 0}^{N} \smash[b]{\underline{\varphi}}_{q} = \frac{\smash[b]{\underline{\varphi}}_{p}}{r+1},\\
&\sum\limits_{q = 0}^{N} \frac{\tau_{q}^{r+1}}{r+1}\smash[b]{\underline{\varphi}}_{q}\smash[b]{\underline{\varphi}}_{p} =
\frac{\smash[b]{\underline{\varphi}}_{p}}{r+1} \sum\limits_{q = 0}^{N} \tau_{q}^{r+1}\smash[b]{\underline{\varphi}}_{q} = \frac{\smash[b]{\underline{\varphi}}_{p}}{r+1} \cdot (0)^{r+1} = 0,\\
&\sum\limits_{q = 0}^{N} \int\limits_{0}^{1} \frac{\varphi_{q}(\tau)}{r+1}\frac{d\varphi_{p}(\tau)}{d\tau}d\tau =
\frac{1}{r+1}\int\limits_{0}^{1} \left[\sum\limits_{q = 0}^{N}\varphi_{q}(\tau)\right] \frac{d\varphi_{p}(\tau)}{d\tau}d\tau =
\frac{\overline{\varphi}_{p} - \smash[b]{\underline{\varphi}}_{p}}{r+1},
\end{split}
\end{equation}
\begin{equation}
\begin{split}
&\sum\limits_{q = 0}^{N} \frac{\tau_{q}^{r+1}}{r+1}\int\limits_{0}^{1} \varphi_{q}(\tau)\frac{d\varphi_{p}(\tau)}{d\tau}d\tau =
\frac{1}{r+1}\int\limits_{0}^{1} \left[\sum\limits_{q = 0}^{N}\tau_{q}^{r+1}\varphi_{q}(\tau)\right] \frac{d\varphi_{p}(\tau)}{d\tau}d\tau\\
& = \frac{1}{r+1} \int\limits_{0}^{1} \tau^{r+1} \frac{d\varphi_{p}(\tau)}{d\tau}d\tau =
\frac{\tau^{r+1}\varphi_{p}(\tau)}{r+1}\Bigg|_{0}^{1} - \int\limits_{0}^{1} \tau^{r} \varphi_{p}(\tau)d\tau =
\frac{\overline{\varphi}_{p}}{r+1} + w_{p}\tau_{p}^{r}.
\end{split}
\end{equation}
Substituting the transformed terms into the original expression turns it into identity:
\begin{equation}
w_{p}\tau_{p}^{r} = 
\frac{\smash[b]{\underline{\varphi}}_{p}}{r+1} - 0 + \frac{\overline{\varphi}_{p} - \smash[b]{\underline{\varphi}}_{p}}{r+1} - \frac{\overline{\varphi}_{p}}{r+1} + w_{p}\tau_{p}^{r} = 
w_{p}\tau_{p}^{r},
\end{equation}
which proves the fulfillment of the simplifying condition $D(N)$.
\end{proof}
\begin{theorem}[\corrtext{\!\!\cite{ader_improving_2024}, Theorem 3.9}]\label{theorem:order_rk}
The approximation order of the ADER-IWF-RK-GLG method is $p_{\rm RK} = 2N + 1$.
\end{theorem}
The proof was based on the use of the famous Butcher's theorem and properties $B(2N+2)$ (Lemma~3.3.1 in~\cite{Dekker_Verwer_1984}), $C(N)$ and $D(N)$ (Lemmas~\ref{lemma:cn_cond} and~\ref{lemma:dn_cond}).
\begin{corollary}\label{coroll:order_ader_dg_nodes}
The ADER-IWF-RK-GLG method is equivalent to the original ADER-DG method, and the numerical solution $\{\mathbf{u}_{n}\}$ of the RK method corresponds to the numerical solution $\{\mathbf{u}_{n}\}$ (\ref{eq:ader_dg_node_sol}) obtained by the ADER-DG method, therefore the approximation order for the solution at the nodes $\{\mathbf{u}_{n}\}$ is also equal $p_{\rm G} = 2N+1$ --- superconvergence is observed.
\end{corollary}
The approximation order $p_{\rm RK} = 2N+1$ of the ADER-IWF-RK-GLG method is one unit less than the approximation order $2N+2$ of the classical GL method.
\begin{corollary}[\corrtext{\!\!\cite{ader_improving_2024}, Remark 3.7}]\label{coroll:rk_is_not_coll}
The ADER-IWF-RK-GLG method is not a collocation method according to Theorem 7.7 in~\cite{Hairer_book_1}, since it has the convergence order $p_{\rm RK} > s$ and the simplifying condition $C(N+1)$ is not satisfied.
\end{corollary}

\section{Local solution and subgrid resolution}

DG methods, and ADER-DG in particular, have the ability to obtain a numerical solution with subgrid resolution, which was noted in~\cite{dg_ivp_ode_1, dg_ivp_ode_2, dg_ivp_ode_3, dg_ivp_ode_4, dg_ivp_ode_5, dg_ivp_ode_6, ader_dg_ivp_ode, ader_dg_ivp_dae}. The ADER-DG method is equivalent to the ADER-IWF-RK-GLG method, which is not a collocation method, so Theorem 7.9 in~\cite{Hairer_book_1} cannot be used directly for it.
\begin{lemma}[\corrtext{\!\!\cite{ader_improving_2024}, Proposition~3.3}]\label{lemma:local_sol_order_L_inf}
The local solution $\mathbf{q}_{n}$ (\ref{eq:local_sol}) of the ADER-DG method has an approximation order $p = N+1$ in $\mathcal{L}_{\infty}$-norm on each discretization domain $\Omega_{n}$:
\begin{equation}
\exists C_{n}\in\mathcal{R}_{+}:\ 
\left|\left|\mathbf{q}_{n} - \mathbf{u}\right|\right|_{\mathcal{L}_{\infty}} \triangleq 
\sup\limits_{t\in\Omega_{n}}\left|\mathbf{q}_{n}(t) - \mathbf{u}(t)\right| < C_{n} \cdot {\Delta t}_{n}^{N+1}.
\end{equation}
\end{lemma}
\begin{proof}
The local solution $\mathbf{q}_{n}$ is the Lagrange interpolation polynomial of degree $N$ (\ref{eq:local_sol}) whose coefficients $\{\hat{\mathbf{q}}_{p}\}_{p = 0}^{N}$ are the solution of the system of equations (\ref{eq:snae_lstdg}), where $\mathbf{u}_{n}$ is known with the approximation order $p_{\rm RK} = 2N+1$. The exact solution of the ODE system (\ref{eq:ivp_ode_diff_src}) in domain $\Omega_{n}$ with the initial condition $\tilde{\mathbf{u}}(t_{n}) = \mathbf{u}_{n}$ is denoted as $\tilde{\mathbf{u}}: \Omega_{n}\rightarrow\mathcal{R}^{D}$. If $\mathbf{F}(\mathbf{u}, t)$ is a sufficiently smooth function, then
\begin{equation}
\exists A_{n}\in\mathcal{R}_{+}:\ \left|\tilde{\mathbf{u}}(t) - \mathbf{u}(t)\right| < 
A_{n} \left|\tilde{\mathbf{u}}(t_{n}) - \mathbf{u}(t_{n})\right| <
A_{n} C_{\rm app} {\Delta t}_{n}^{2N+2},
\end{equation}
where $C_{\rm app}$ is the approximation constant of the ADER-IWF-RK-GLG method, and $2N+2$ in the exponent is associated with the approximation order $p_{\rm RK} = 2N+1$ at one step $\varepsilon = O({\Delta t}_{n}^{p_{\rm RK}+1})$, expressed through the local truncation error $\varepsilon$. The error estimate for the Lagrange interpolation polynomial can be presented in the form:
\begin{equation}
\left|\mathbf{q}_{n}(t) - \tilde{\mathbf{u}}(t)\right| \leq 
\frac{{\Delta t}_{n}^{N+1}}{(N+1)!}\max\limits_{t\in\Omega_{n}}\left|\frac{d^{N+1}\tilde{\mathbf{u}}(t)}{dt^{N+1}}\right| \equiv B_{n} {\Delta t}_{n}^{N+1},
\end{equation}
which is also applicable for sufficiently smooth functions $\mathbf{F}(\mathbf{u}, t)$ due to Lemma~\ref{lemma:p_ode_exact}. As a result of applying the norm property, the following expression is obtained:
\begin{equation}
\begin{split}
\left|\mathbf{q}_{n}(t) - \mathbf{u}(t)\right| &\leq 
\left|\mathbf{q}_{n}(t) - \tilde{\mathbf{u}}(t)\right| + 
\left|\tilde{\mathbf{u}}(t) - \mathbf{u}(t)\right|\\
&\leq B_{n} {\Delta t}_{n}^{N+1} + A_{n} C_{\rm app} {\Delta t}_{n}^{2N+2} \equiv O\left({\Delta t}_{n}^{N+1}\right),
\end{split}
\end{equation}
where the initial provable assumption follows.
\end{proof}
\begin{theorem}\label{theorem:local_sol_order_L_inf}
The local solution $\mathbf{u}_{L}$ (\ref{eq:local_sol_assembly}) of the ADER-DG method has an approximation order $p_{\rm L} = N+1$ in $\mathcal{L}_{\infty}$-norm:
\begin{equation}
\exists C\in\mathcal{R}_{+}:\ 
\left|\left|\mathbf{u}_{L} - \mathbf{u}\right|\right|_{\mathcal{L}_{\infty}} \hspace{-2mm}\triangleq 
\operatorname{ess}\sup\limits_{\hspace{-5mm}t\in\Omega}\left|\mathbf{u}_{L}(t) - \mathbf{u}(t)\right| < C \cdot {\Delta t}^{N+1},\ 
{\Delta t} = \max\limits_{n} {\Delta t}_{n}.
\end{equation}
\end{theorem}
\begin{proof}
Taking into account the fact of discretization of the domain of definition $\Omega$ by a finite set of $\{\Omega_{n}\}$ and in the case $\mathbf{u}\in\mathcal{C}_{0}(\Omega)$, the following expression is obtained:
\begin{equation}
\operatorname{ess}\sup\limits_{\hspace{-5mm}t\in\Omega}\left|\mathbf{u}_{L}(t) - \mathbf{u}(t)\right| =
\max\limits_{n}\sup\limits_{t\in\Omega_{n}}\left|\mathbf{u}_{L}(t) - \mathbf{u}(t)\right| = 
\max\limits_{n}\sup\limits_{t\in\Omega_{n}}\left|\mathbf{q}_{n}(t) - \mathbf{u}(t)\right|.
\end{equation}
Following Lemma~\ref{lemma:local_sol_order_L_inf}, the following estimate is obtained:
\begin{equation}
\begin{split}
\max\limits_{n} \sup\limits_{t\in\Omega_{n}}\left|\mathbf{q}_{n}(t) - \mathbf{u}(t)\right| < 
\max\limits_{n} C_{n} \cdot {\Delta t}_{n}^{N+1} < C \cdot {\Delta t}^{N+1},\ C > \max\limits_{n} C_{n},
\end{split}
\end{equation}
where the initial provable assumption follows.
\end{proof}
\begin{corollary}\label{coroll:local_sol_order_qpoints}
The local solution $\mathbf{u}_{L}$ at the nodal points $t_{n}+\tau_{p}{\Delta t}_{n}$ coincides with the coefficients $\{\hat{\mathbf{q}}_{n, p}\}$, so the approximation order of the local solution at the nodal points $\{\hat{\mathbf{q}}_{n, p}\}$ is also $p_{\rm L} = N+1$.
\end{corollary}
\begin{theorem}\label{theorem:local_sol_order_L_p}
The local solution $\mathbf{u}_{L}$ (\ref{eq:local_sol_assembly}) of the ADER-DG method has an approximation order $p_{\rm L} = N+1$ in $\mathcal{L}_{q}$-norm:
\begin{equation}
\exists C_{q}\in\mathcal{R}_{+}:\ 
\left|\left|\mathbf{u}_{L} - \mathbf{u}\right|\right|_{\mathcal{L}_{q}} \hspace{-2mm}\triangleq 
\left[\int\limits_{\Omega}\left|\mathbf{u}_{L}(t) - \mathbf{u}(t)\right|^{q} dt\right]^{1/q} < C_{q} \cdot {\Delta t}^{N+1},\ 
{\Delta t} = \max\limits_{n} {\Delta t}_{n}.
\end{equation}
\end{theorem}
\begin{proof}
Using the inequality $|\mathbf{f}| < ||\mathbf{f}||_{\mathcal{L}_{\infty}}$ in the assumption:
\begin{equation}
\int\limits_{\Omega}\left|\mathbf{u}_{L}(t) - \mathbf{u}(t)\right|^{q} dt \leq
\int\limits_{\Omega} \corrtext{\left|\left|\mathbf{u}_{L} - \mathbf{u}\right|\right|_{\mathcal{L}_{\infty}}^{q}} dt =
\left|\left|\mathbf{u}_{L}(t) - \mathbf{u}(t)\right|\right|_{\mathcal{L}_{\infty}}^{q}\cdot|\Omega|,
\end{equation}
where $|\Omega| = t_{f} - t_{0}$, and according to Theorem~\ref{theorem:local_sol_order_L_inf}, the estimate is obtained:
\begin{equation}
\begin{split}
\left|\left|\mathbf{u}_{L} - \mathbf{u}\right|\right|_{\mathcal{L}_{q}} \leq
\left[|\Omega|\cdot\left|\left|\mathbf{u}_{L}(t) - \mathbf{u}(t)\right|\right|_{\mathcal{L}_{\infty}}^{q}\right]^{1/q} \leq
C|\Omega|^{1/q} \cdot {\Delta t}^{N+1} \equiv C_{q} \cdot {\Delta t}^{N+1},
\end{split}
\end{equation}
where the initial provable assumption follows.
\end{proof}

\section{Linear and nonlinear stability}

The analysis of nonlinear stability of the ADER-IWF-RK-GLG method is based on the study of the properties of non-negative definiteness of matrices $\mu$, $\mathrm{Q}$ and $\mathrm{M}$~\cite{Dekker_Verwer_1984}:
\begin{equation}\label{eq:q_m_matrices_def}
\mathrm{Q} = \mu\alpha + \alpha^{T}\mu - \alpha^{T}\mathbf{w}\mathbf{w}^{T}\alpha,\quad
\mathrm{M} = a^{T}\mathrm{Q}a = \mu a + a^{T} \mu - \mathbf{w}\mathbf{w}^{T},
\end{equation}
where $\mathbf{w} = [w_{0} \ldots w_{N}]^{T}$ and $\alpha = a^{-1} = [\kappa^{-1}\mu]^{-1} = \mu^{-1}\kappa$. Matrix $\mu$ is positive definite, since it is diagonal and contains $w_{p} > 0$ on the diagonal.
\begin{lemma}\label{lemma:q_matrix_comp}
Matrix $\mathrm{Q} = \smash[b]{\underline{\boldsymbol{\varphi}}}\smash[b]{\underline{\boldsymbol{\varphi}}}^{T}$, where $\smash[b]{\underline{\boldsymbol{\varphi}}} = \lbrace\smash[b]{\underline{\varphi}}_{p}\rbrace_{p} \neq 0$.
\end{lemma}
\begin{proof}
The proof is based on the element-by-element description of the matrix $\mathrm{Q} = \lbrace Q_{pq}\rbrace_{p,q}$ and the use of the definitions of matrices $\kappa$ and $\mu$ (\ref{eq:kappa_mu_defs}):
\begin{equation}
Q_{pq} = w_{p}\alpha_{pq} + w_{q}\alpha_{qp} - \sum\limits_{i = 0}^{N}\sum\limits_{j = 0}^{N} \alpha_{ip}w_{i}w_{j}\alpha_{jq} =
\kappa_{pq} + \kappa_{qp} - \sum\limits_{i = 0}^{N}\sum\limits_{j = 0}^{N} \kappa_{ip}\kappa_{jq}.
\end{equation}
Then each term is studied and folded separately:
\begin{equation}\nonumber
\begin{split}
&\kappa_{pq} + \kappa_{qp} = 
\overline{\varphi}_{p}\overline{\varphi}_{q} + \overline{\varphi}_{q}\overline{\varphi}_{p} - 
\hspace{-1mm}\int\limits_{0}^{1} \left[\frac{d\varphi_{p}(\tau)}{d\tau}\varphi_{q}(\tau) + \frac{d\varphi_{q}(\tau)}{d\tau}\varphi_{p}(\tau)\right] d\tau =
\overline{\varphi}_{p}\overline{\varphi}_{q} + \smash[b]{\underline{\varphi}}_{p}\smash[b]{\underline{\varphi}}_{q},\\
&\sum\limits_{i = 0}^{N}\sum\limits_{j = 0}^{N} \kappa_{ip}\kappa_{jq} =
\corrtext{\Bigg[\sum\limits_{i = 0}^{N}\kappa_{ip}\Bigg]\Bigg[\sum\limits_{j = 0}^{N}\kappa_{jq}\Bigg]}
= \overline{\varphi}_{p}\overline{\varphi}_{q}.
\end{split}
\end{equation}
The final substitution of the folded expressions of the terms into the original expression leads to the result $Q_{pq} = \corrtext{\smash[b]{\underline{\varphi}}_{p}\smash[b]{\underline{\varphi}}_{q}}$, which is the result for matrix $\mathrm{Q}$ being proven.
\end{proof}
\begin{lemma}\label{lemma:q_m_matrices_posdef}
Matrices $\mathrm{Q}$ and $\mathrm{M}$ are non-negative definite.
\end{lemma}
\begin{proof}
The dyadic matrix $\mathbf{a}\mathbf{b}^{T}$, where $\mathbf{a}\neq\mathbf{0}$, $\mathbf{b}\neq\mathbf{0}$, has all zero eigenvalues except one $\lambda = \mathbf{a}^{T}\mathbf{b}$. In the case of dyadic square $\mathbf{a}\mathbf{a}^{T}$, where $\mathbf{a}\neq\mathbf{0}$, $\lambda = \mathbf{a}^{T}\mathbf{a} = |\mathbf{a}|^{2} > 0$, so dyadic squares are non-negative definite. Matrix $\mathrm{Q}$ is dyadic square by Lemma~\ref{lemma:q_matrix_comp}, so it is non-negative definite. Matrix $\mathrm{M} = a^{T}\mathrm{Q}a = a^{T}\smash[b]{\underline{\boldsymbol{\varphi}}}\smash[b]{\underline{\boldsymbol{\varphi}}}^{T}a = (a^{T}\smash[b]{\underline{\boldsymbol{\varphi}}})(a^{T}\smash[b]{\underline{\boldsymbol{\varphi}}})^{T}$, and is also dyadic square, so it is non-negative definite, since $a$ is not degenerate and $\smash[b]{\underline{\boldsymbol{\varphi}}} = \lbrace\smash[b]{\underline{\varphi}}_{p}\rbrace_{p} \neq 0$.
\end{proof}
\begin{corollary}
$\lambda_{Q} = |\smash[b]{\underline{\boldsymbol{\varphi}}}|^{2} = \sum_{p = 0}^{N} \overline{\varphi}_{p}^{2} > 0$ and $\lambda_{M} = |a^{T} \smash[b]{\underline{\boldsymbol{\varphi}}}|^{2} > 0$ are the only non-negative eigenvalues of $\mathrm{Q}$ and $\mathrm{M}$, respectively.
\end{corollary}
\begin{theorem}\label{theorem:alg_stab_rk}
The ADER-IWF-RK-GLG method is algebraically stable.
\end{theorem}
\begin{proof}
The matrix $\mathrm{M}$ (\ref{eq:q_m_matrices_def}) is non-negative definite by Lemma~\ref{lemma:q_m_matrices_posdef}, and matrix $\mu$ (\ref{eq:kappa_mu_defs}) is positive definite, so the ADER-IWF-RK-GLG method is algebraically stable, by Definition~4.2.1 in~\cite{Dekker_Verwer_1984}.
\end{proof}
\begin{corollary}\label{coroll:b_bn_stab_rk}
The ADER-IWF-RK-GLG method is $BN$-stable by the Theorem~4.2.2 in~\cite{Dekker_Verwer_1984}, and $B$-stable as a special case.
\end{corollary}
RK methods are usually defined by matrix $a$, so the stability analysis is carried out not by $\mathrm{Q}$, but by $\mathrm{M}$. It is more convenient to use matrix $\mathrm{Q}$ for the ADER-IWF-RK-GLG method, since it is more convenient to use $a^{-1}$ containing $\kappa$, not $\kappa^{-1}$.
\begin{remark}
Theorem~4.6.1 in~\cite{Dekker_Verwer_1984} shows $\mathrm{rank}(\mathrm{M}) \leq 2s - p_{\rm RK} = 1$, and the Lemma~\ref{lemma:q_m_matrices_posdef} leads to result $\mathrm{rank}(\mathrm{M}) = 1$, since the matrix $\mathrm{M}$ is dyadic.
\end{remark}
\begin{remark}
According to Lemma~\ref{lemma:q_m_matrices_posdef}, by Theorem 4.1.3 in~\cite{Dekker_Verwer_1984}, $B$-stability of the ADER-IWF-RK-GLG method can be proven by the non-negative definiteness of $\mathrm{Q}$.
\end{remark}
\begin{theorem}\label{theorem:an_stab_rk}
The ADER-IWF-RK-GLG method is $AN$-stable.
\end{theorem}
\begin{proof}
Legendre polynomials $\tilde{P}_{N+1}$ do not have multiple roots $\tau_{p}$, so all abscissas $\tau_{p}$ are different: $\tau_{p} = \tau_{q}$ if and only if $p = q$, hence the ADER-IWF-RK-GLG method is non-confluent. By Theorem~4.3.5 in~\cite{Dekker_Verwer_1984} any non-confluent algebraically stable method is $AN$-stable. This proves the $AN$-stability of the ADER-IWF-RK-GLG method.
\end{proof}
\begin{corollary}\label{coroll:a_stab_rk}
The ADER-IWF-RK-GLG method is $A$-stable as a special case.
\end{corollary}
\begin{theorem}[\corrtext{\!\!\cite{ader_proofs_2025}, Theorem 3.7}]\label{theorem:pade_approx}
The stability function $R(z)$ of the ADER-IWF-RK-GLG method is the $(N, N+1)$-Pad\'{e} approximation $R_{N, N+1}(z)$ of $\exp(z)$.
\end{theorem}
Symbolic calculations in the \texttt{Maple} for degrees $N = 1,\, \ldots,\, 6$ showed that the function $R(z)$ coincides with the $(N, N+1)$-Pad\'{e} approximation of $\exp(z)$. Numerical calculations using module \texttt{mpmath} of the \texttt{python} programming language, with an arithmetic precision constant \texttt{mpmath.mp.dps = 1000}, also confirmed this result for degrees $N = 1,\, \ldots,\, 75$ with accuracy $\approx 10^{-946}$--$10^{-1002}$.
\begin{corollary}[\corrtext{\!\!\cite{ader_proofs_2025}, Corollary 3.8}]
The $(N, N+1)$-Pad\'{e} approximation of $\exp(z)$ is $A$-admissible function by Theorem~3.4.8 in~\cite{Dekker_Verwer_1984}, so the Theorem~\ref{theorem:pade_approx} provides an alternative proof of the $A$-stability of the ADER-IWF-RK-GLG method.
\end{corollary}
\begin{corollary}
By Theorem~3.11 in~\cite{Hairer_book_2}, the residual term of the $R_{N, N+1}(z)$ is $\exp(z) - R_{N, N+1}(z) = O(z^{2N+2})$, so the approximation order of the ADER-IWF-RK-GLG method $p_{\rm RK} = 2N+1$, which is an alternative proof of Theorem~\ref{theorem:order_rk}.
\end{corollary}
\begin{theorem}\label{theorem:l2_stab_rk}
The ADER-IWF-RK-GLG method is $L$-stable.
\end{theorem}
\begin{proof}
The ADER-IWF-RK-GLG method is $A$-stable, and the stability function $R(z)$ has asymptotics $|z|^{-1}$, $z\in\mathcal{C}$, in domain $z\rightarrow\infty$ as $(N, N+1)$-Pad\'{e} approximation, so the ADER-IWF-RK-GLG method is $L$-stable.
\end{proof}
\begin{remark}
The classical GL method is only $A$-stable, but not $L$-stable. Decreasing the approximation order ``led to'' an increase in the stability.
\end{remark}
\begin{corollary}\label{coroll:ader_dg_stabs}
The ADER-IWF-RK-GLG method is equivalent to the original ADER-DG method, therefore, the ADER-DG method is $A$-stable, $AN$-stable, $L$-stable, algebraically stable, $B$-stable and $BN$-stable.
\end{corollary}
$A$-stability and $L$-stability of the ADER-DG method were numerically demonstrated in~\cite{ader_dg_ivp_ode}, and in this paper, these and other types of numerical stability are proven.

\section{Computational results}
The Section presents examples of using the ADER-DG to solve a system of ODEs, which demonstrate the convergence and stability properties of the method. The software implementation of the ADER-DG numerical method is developed in the \texttt{python} programming language. The error of the numerical solution becomes very small even on fairly coarse grids. Therefore, floating-point numbers of arbitrary precision are used in the \texttt{mpmath} module, and \texttt{mpmath.mp.dps = 500} is chosen. In the implementation of the ADER-DG method, the system of equations (\ref{eq:ader_dg_rkm}) was solved to calculate $\{\mathbf{k}_{p}\}$ (for which numerical methods are developed~\cite{Butcher_book_2016, Hairer_book_1, Hairer_book_2, Dekker_Verwer_1984}), from which $\{\mathbf{q}_{p}\}$ are calculated (\ref{eq:k_to_q}), unlike~\cite{ader_dg_ivp_ode}, where direct system (\ref{eq:snae_lstdg}) of the local DG predictor was solved, although of course the result did not change. \corrtext{The system of nonlinear algebraic equations (\ref{eq:snae_lstdg}) is solved using Newton's method until an absolute error $\epsilon_{0} = 10^{-490}$ is reached, which required no more than $11$-$20$ iterations. Relative error $10^{-12}$ was achieved in no more than $4$-$5$ iterations.}

Convergence is verified by solving two ODE systems. The first example is a classical one-dimensional harmonic oscillator described by a linear ODE:
\begin{equation}\label{eq:harm_osc}
\ddot{x} + x = 0,\quad x(0) = 1,\quad \dot{x}(0) = 0,\quad t\in[0,\, 4\pi],
\end{equation}
and the second example is a mathematical pendulum described by a nonlinear ODE:
\begin{equation}\label{eq:pend}
\ddot{\phi} + \sin(\phi) = 0,\quad \phi(0) = \frac{\pi}{2},\quad \dot{\phi}(0) = 0,\quad t\in[0,\, 10].
\end{equation}
Exact analytical solutions of these ODE systems are obtained by trivial integration.

The solution at grid nodes $\mathbf{u}_{n}$, the local solution $\mathbf{u}_{L}$ (\ref{eq:local_sol_assembly}) in space between nodes $\{\Omega_{n}\}$, and the local solution $\mathbf{u}_{L}(t_{n, p})$ at nodal points $\{t_{n}+\tau_{p}{\Delta t}\}$ are calculated and analyzed. The local error $\varepsilon:\, \Omega\rightarrow\mathcal{R}_{+}$ is calculated using the uniform norm of the solution vector $\mathbf{u}$: $\varepsilon(t) = |\mathbf{u}(t) - \mathbf{u}^{\rm ex}(t)|$, where $\mathbf{u}^{\rm ex}$ is the exact solution. The calculation of empirical convergence orders $p$ is carried out separately for each presented type of numerical solution, and is based on a power approximation of the dependence of global errors $e({\Delta t}) \propto {\Delta t}^{p}$ on the discretization step ${\Delta t}$. The domain of definition $\Omega$ was discretized into $M$ discretization domains $\{\Omega_{n}\}$, with a discretization step ${\Delta t} = (t_{f} - t_{0})/M$. The following \corrtext{$10$} types of global errors $e$:
\begin{equation}\nonumber
\begin{array}{ll}
e^{n}_{L_{1}} = {\displaystyle\sum\limits_{n = 0}^{M}} {\Delta t}_{n}\left|\mathbf{u}_{n} - \mathbf{u}^{\rm ex}(t_{n})\right|,&
e^{l,q}_{L_{1}} = {\displaystyle\sum\limits_{n = 0}^{M}\sum\limits_{p = 0}^{N}} {\Delta t}_{n, p}\left|\mathbf{u}_{L}(t_{n, p}) - \mathbf{u}^{\rm ex}(t_{n, p})\right|,\\[5mm]
e^{n}_{L_{2}} = \left[{\displaystyle\sum\limits_{n = 0}^{M}} {\Delta t}_{n}\left|\mathbf{u}_{n} - \mathbf{u}^{\rm ex}(t_{n})\right|^{2}\right]^{\frac{1}{2}},&
e^{l,q}_{L_{2}} = \left[{\displaystyle\sum\limits_{n = 0}^{M}\sum\limits_{p = 0}^{N}} {\Delta t}_{n, p}\left|\mathbf{u}_{L}(t_{n, p}) - \mathbf{u}^{\rm ex}(t_{n, p})\right|^{2}\right]^{\frac{1}{2}}\hspace{-3mm},\\[7.5mm]
e^{n}_{L_{\infty}} = \max\limits_{0 \leqslant n \leqslant M}\left|\mathbf{u}_{n} - \mathbf{u}^{\rm ex}(t_{n})\right|,&
e^{l,q}_{L_{\infty}} = \max\limits_{0 \leqslant n \leqslant M}\max\limits_{0 \leqslant p \leqslant N}\left|\mathbf{u}_{L}(t_{n, p}) - \mathbf{u}^{\rm ex}(t_{n, p})\right|,\\[7.5mm]
e^{n}_{f} = \left|\mathbf{u}_{M} - \mathbf{u}^{\rm ex}(t_{f})\right|,&
e^{l}_{L_{\infty}} = \max\limits_{0 \leqslant n \leqslant M}\sup\limits_{t\in\Omega_{n}}\left|\mathbf{u}_{L}(t) - \mathbf{u}^{\rm ex}(t)\right|,\\[3mm]
e^{l}_{L_{1}} = {\displaystyle\sum\limits_{n = 0}^{M-1}}\ {\displaystyle\int\limits_{\Omega_{n}}}\left|\mathbf{u}_{L}(t) - \mathbf{u}^{\rm ex}(t)\right| dt,&
e^{l}_{L_{2}} = \left[{\displaystyle\sum\limits_{n = 0}^{M-1}}\ {\displaystyle\int\limits_{\Omega_{n}}}\left|\mathbf{u}_{L}(t) - \mathbf{u}^{\rm ex}(t)\right|^{2} dt\right]^{\frac{1}{2}},
\end{array}
\end{equation}
where $t_{n, p} = t_{n}+\tau_{p}{\Delta t}$ is nodal points, ${\Delta t}_{n, p}$ is distance to the next nodal point or $t_{f}$, and \corrtext{$10$} corresponding empirical convergence orders $p$ in different norms are selected.

\begin{figure}[h!]
\centering
\includegraphics[width=0.28\textwidth]{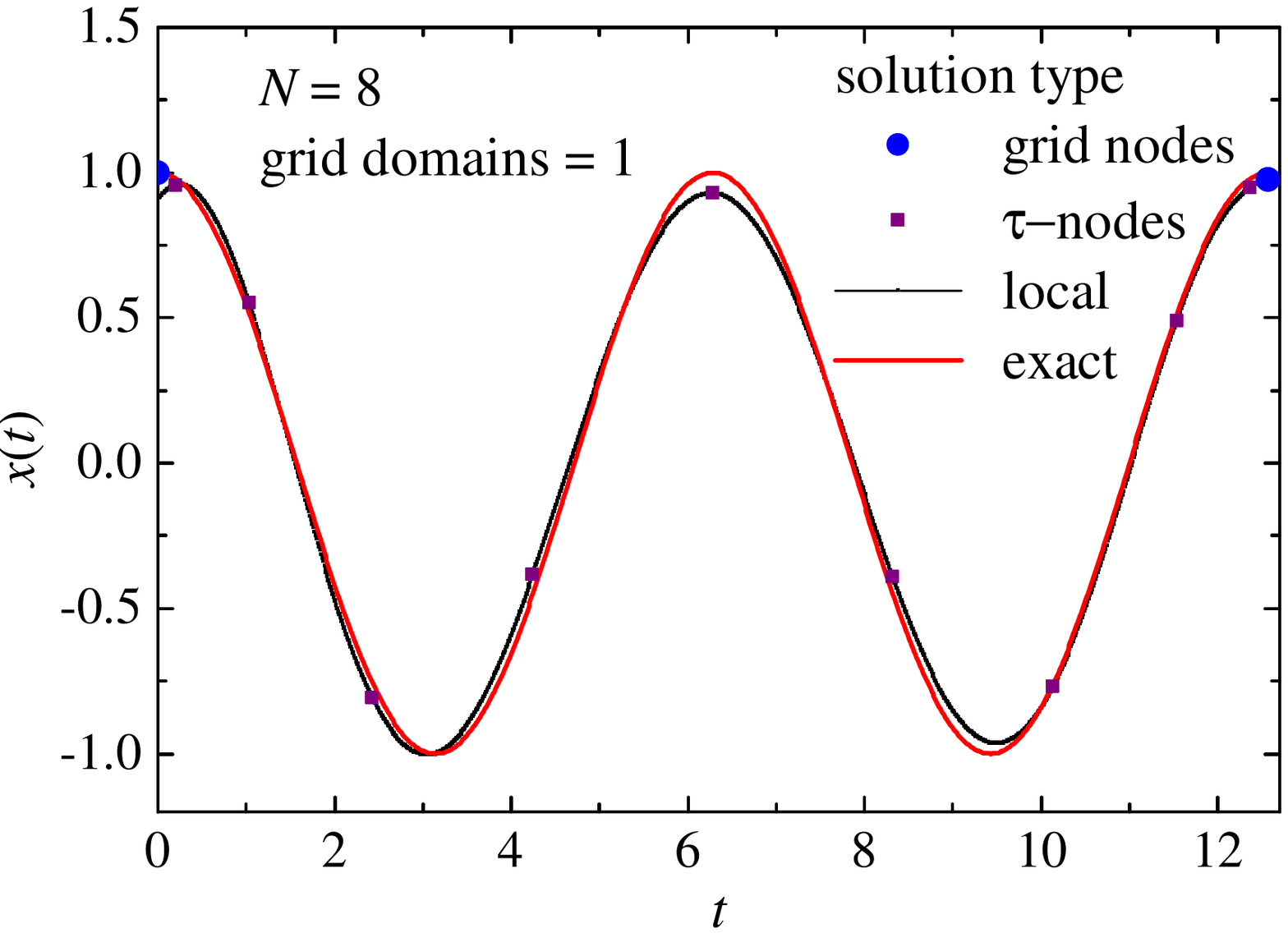}\hspace{4mm}
\includegraphics[width=0.28\textwidth]{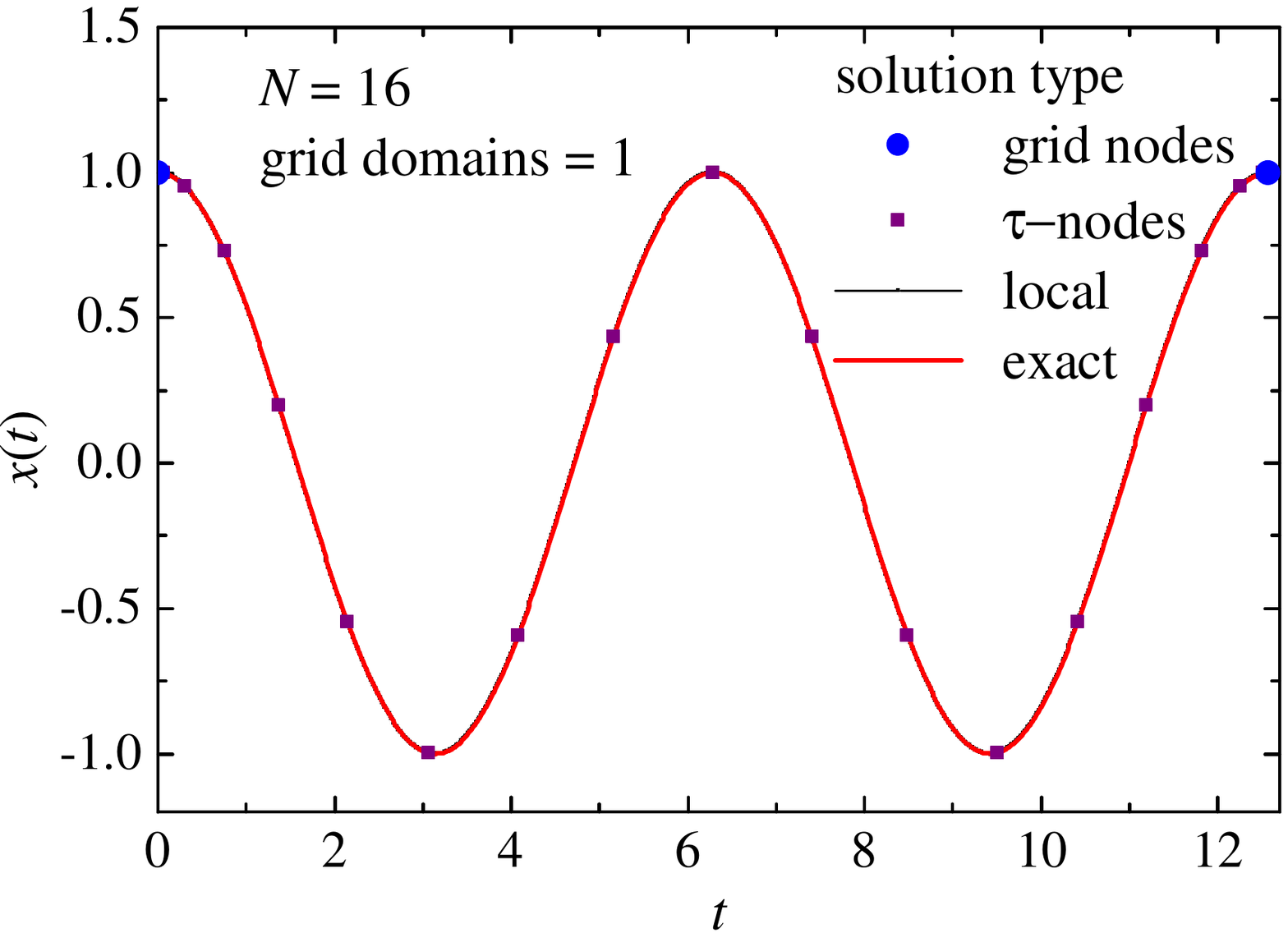}\hspace{4mm}
\includegraphics[width=0.28\textwidth]{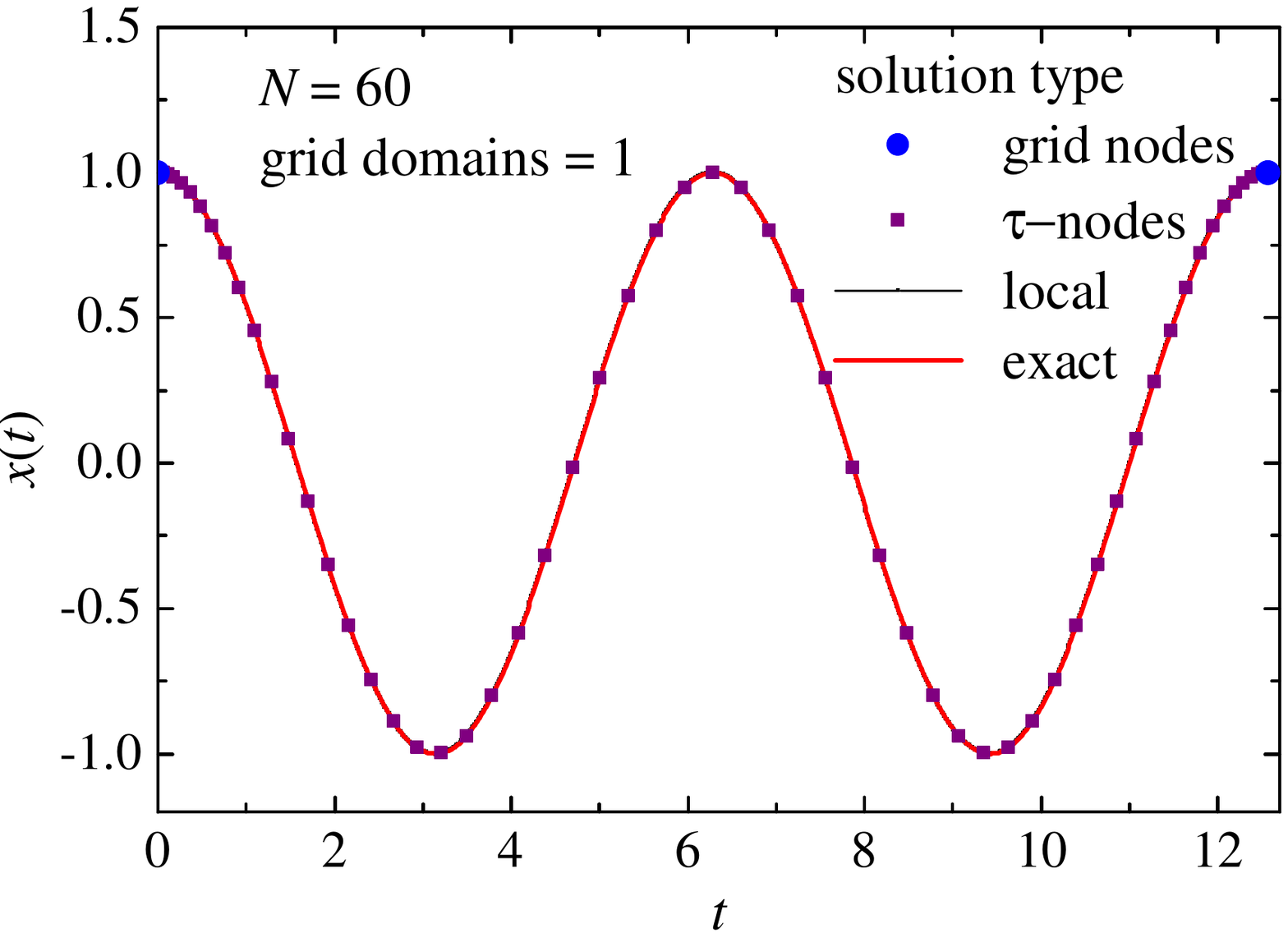}\\
\includegraphics[width=0.28\textwidth]{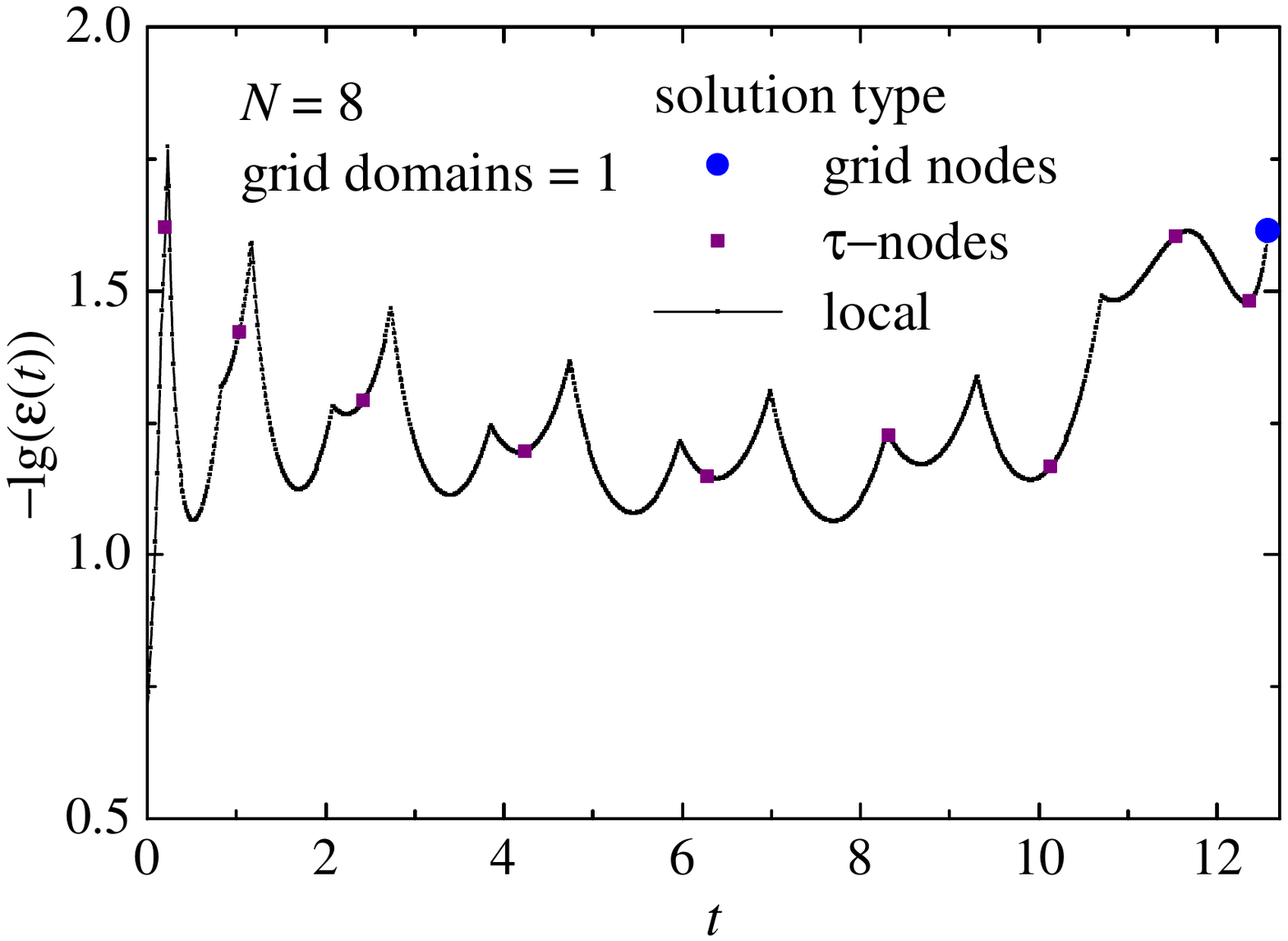}\hspace{4mm}
\includegraphics[width=0.28\textwidth]{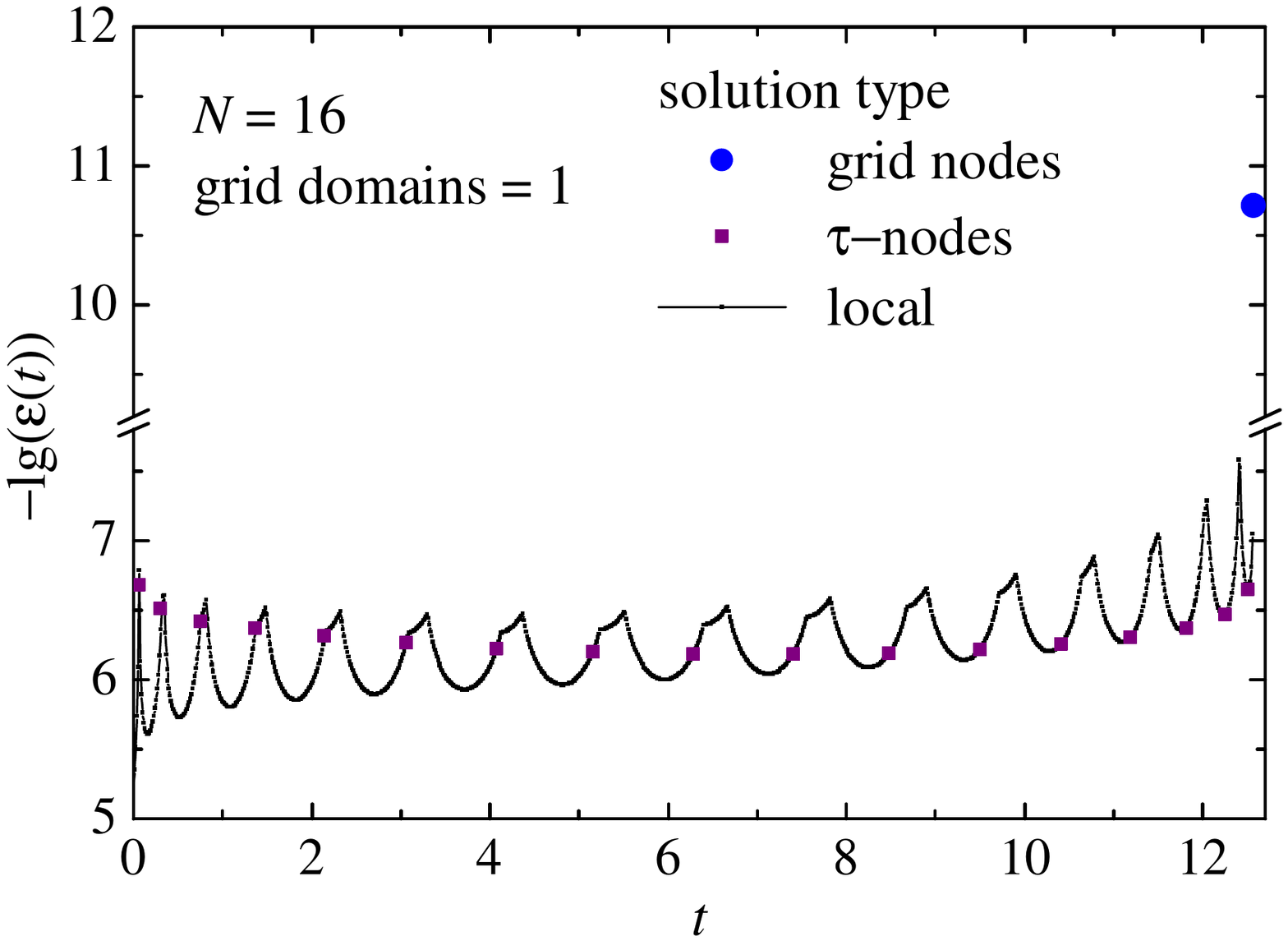}\hspace{4mm}
\includegraphics[width=0.28\textwidth]{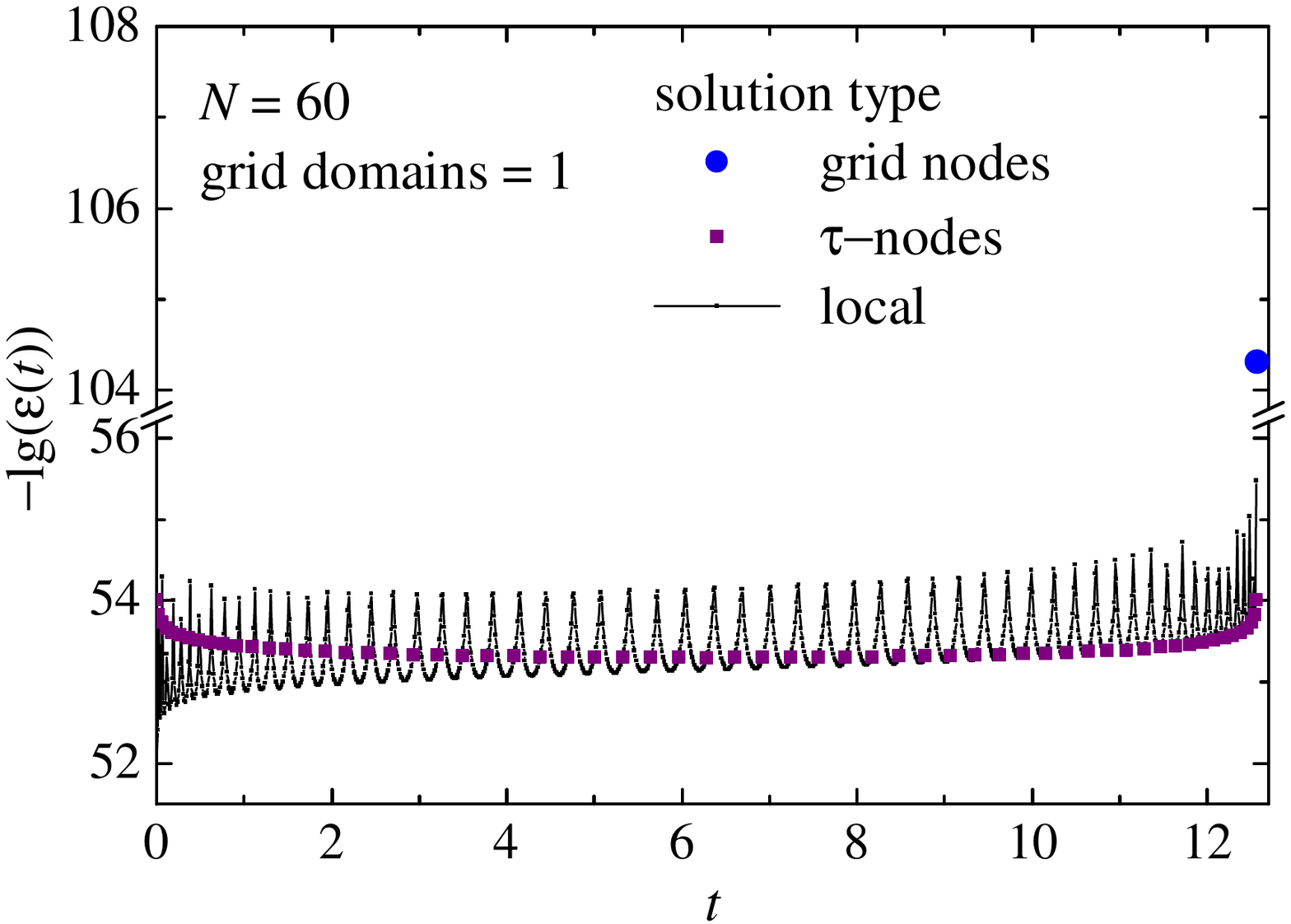}\\
\caption{\label{fig:harm_osc_sols}%
Numerical solution $x(t)$ and its error $\varepsilon(t)$, as the negative common logarithm $-\lg(\varepsilon(t))$, obtained by the ADER-DG method with $N = 8$, $16$, $60$ on a grid with one discretization domain for problem (\ref{eq:harm_osc}). The red line shows the exact analytical solution for comparison.
}
\end{figure}

\begin{figure}[h!]
\centering
\includegraphics[width=0.28\textwidth]{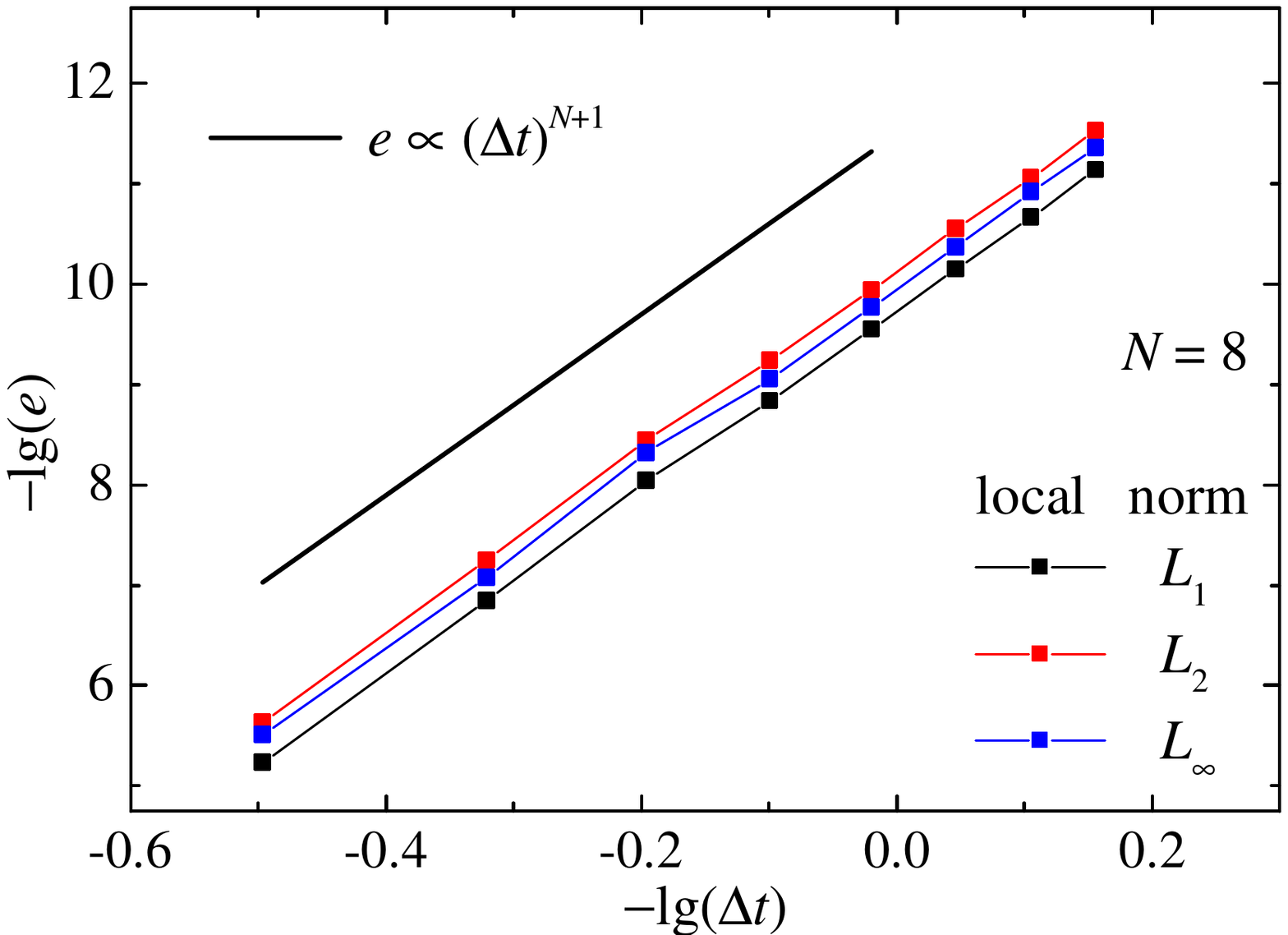}\hspace{4mm}
\includegraphics[width=0.28\textwidth]{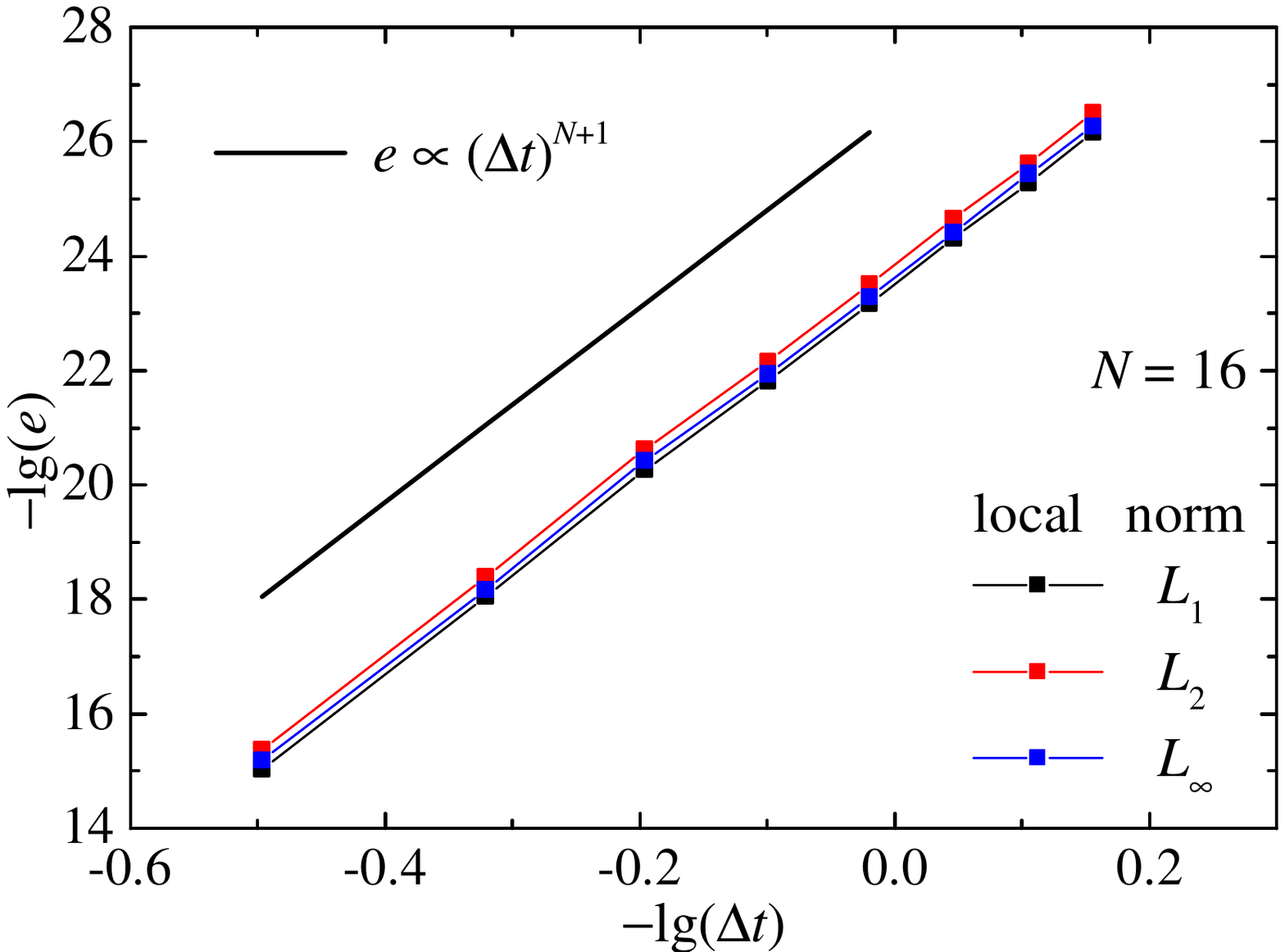}\hspace{4mm}
\includegraphics[width=0.28\textwidth]{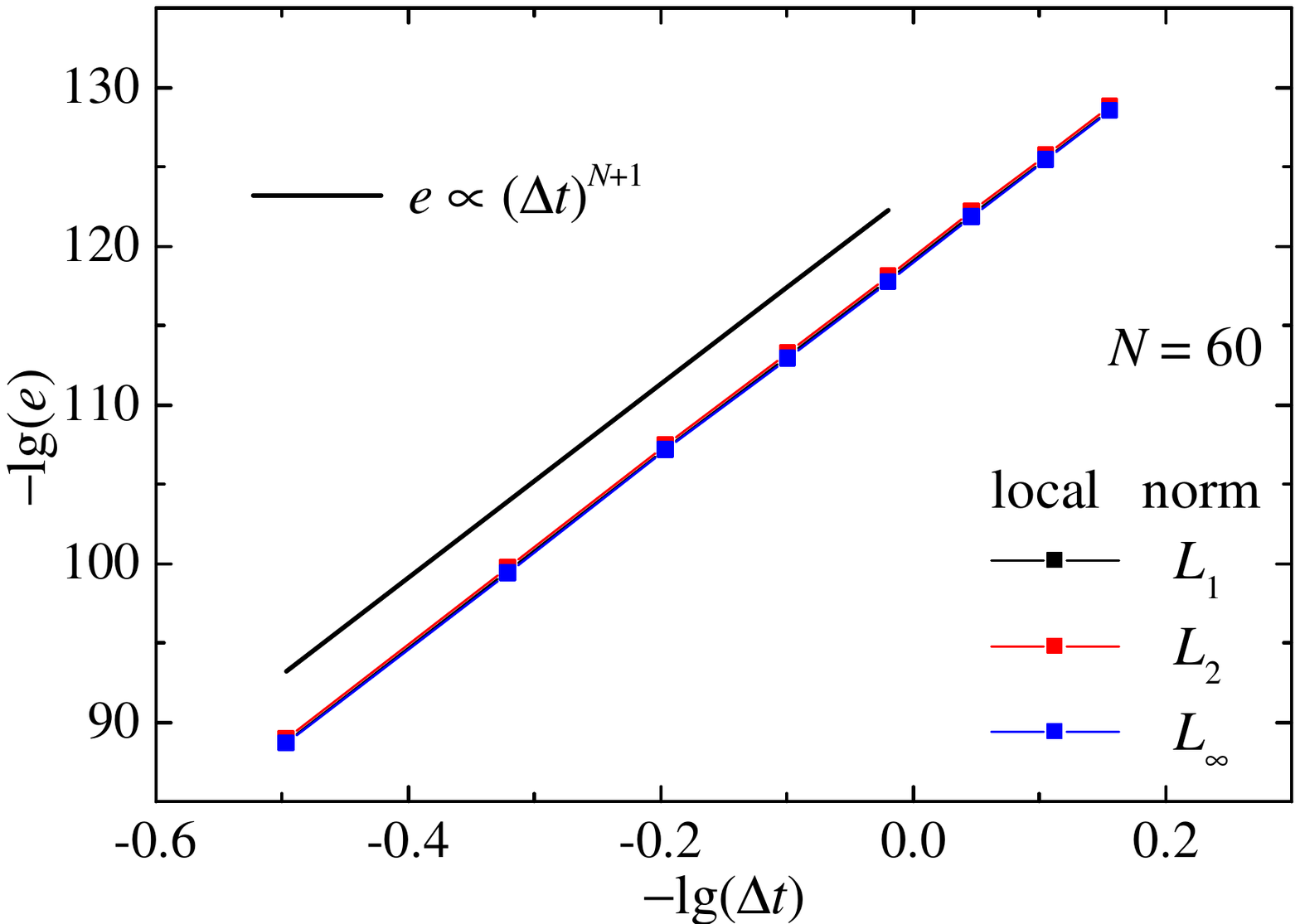}\\
\caption{\label{fig:errs_harm_osc_local}%
Dependence of error $e^{l}$ on discretization step ${\Delta t}$ for local solution $\mathbf{u}_{L}$ in space between nodes $\{\Omega_{n}\}$, calculated in norms $\mathcal{L}_{1}$, $\mathcal{L}_{2}$, $\mathcal{L}_{\infty}$, obtained by the ADER-DG method with $N = 8$, $16$, $60$. Lines represent reference slopes for convergence order $p_{\rm L} = N+1$.
}
\end{figure}

\begin{figure}[h!]
\centering
\includegraphics[width=0.28\textwidth]{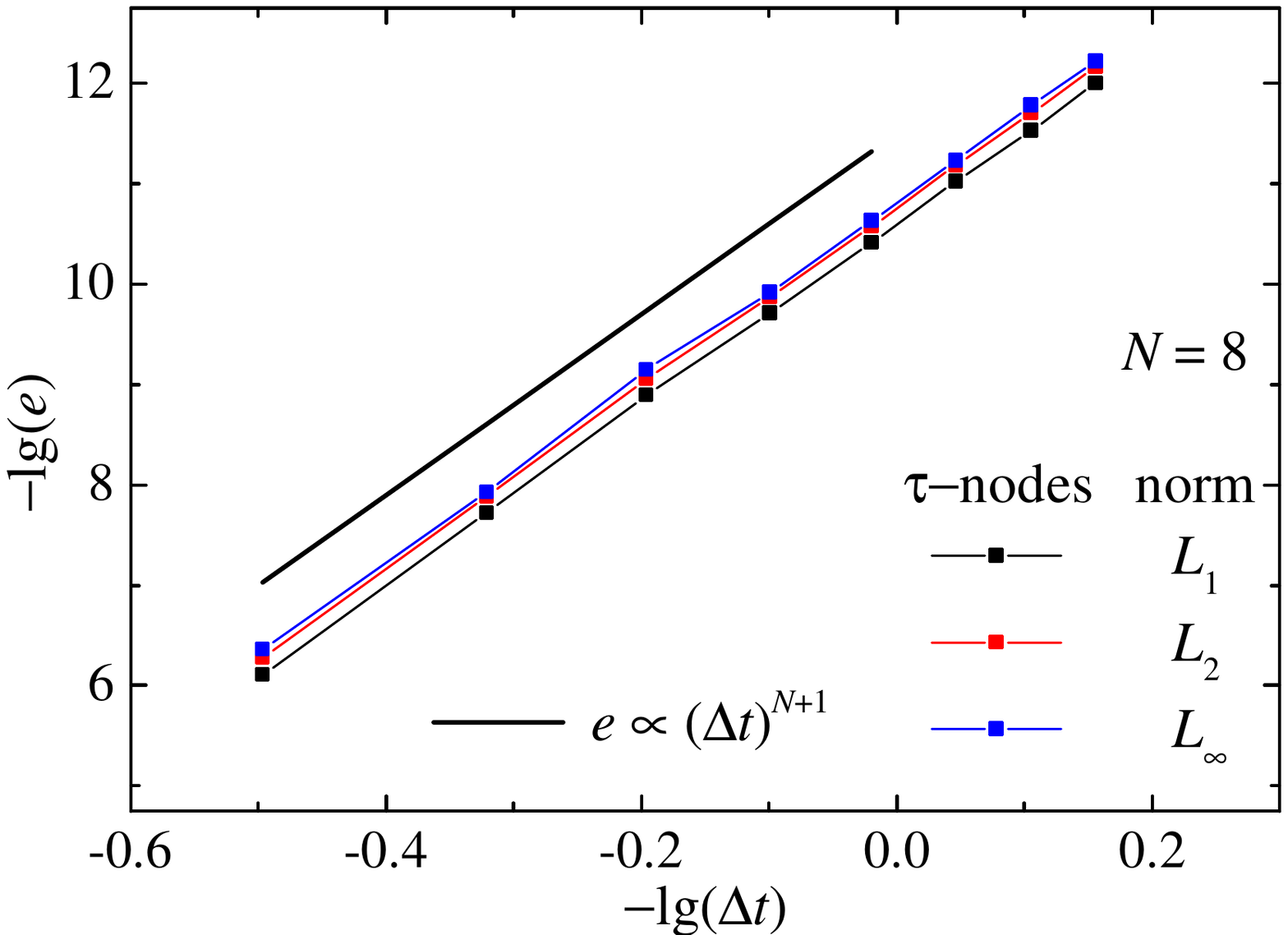}\hspace{4mm}
\includegraphics[width=0.28\textwidth]{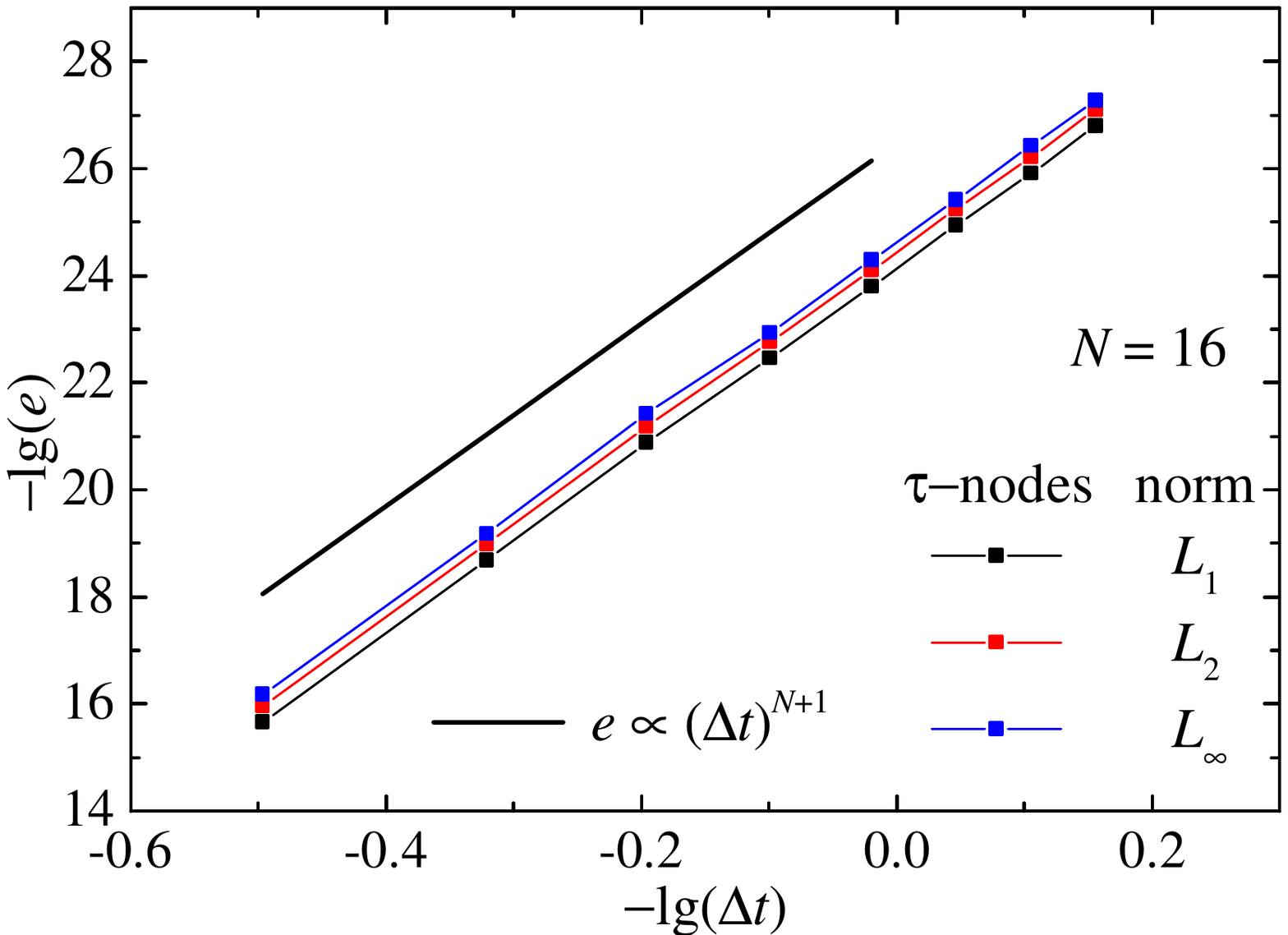}\hspace{4mm}
\includegraphics[width=0.28\textwidth]{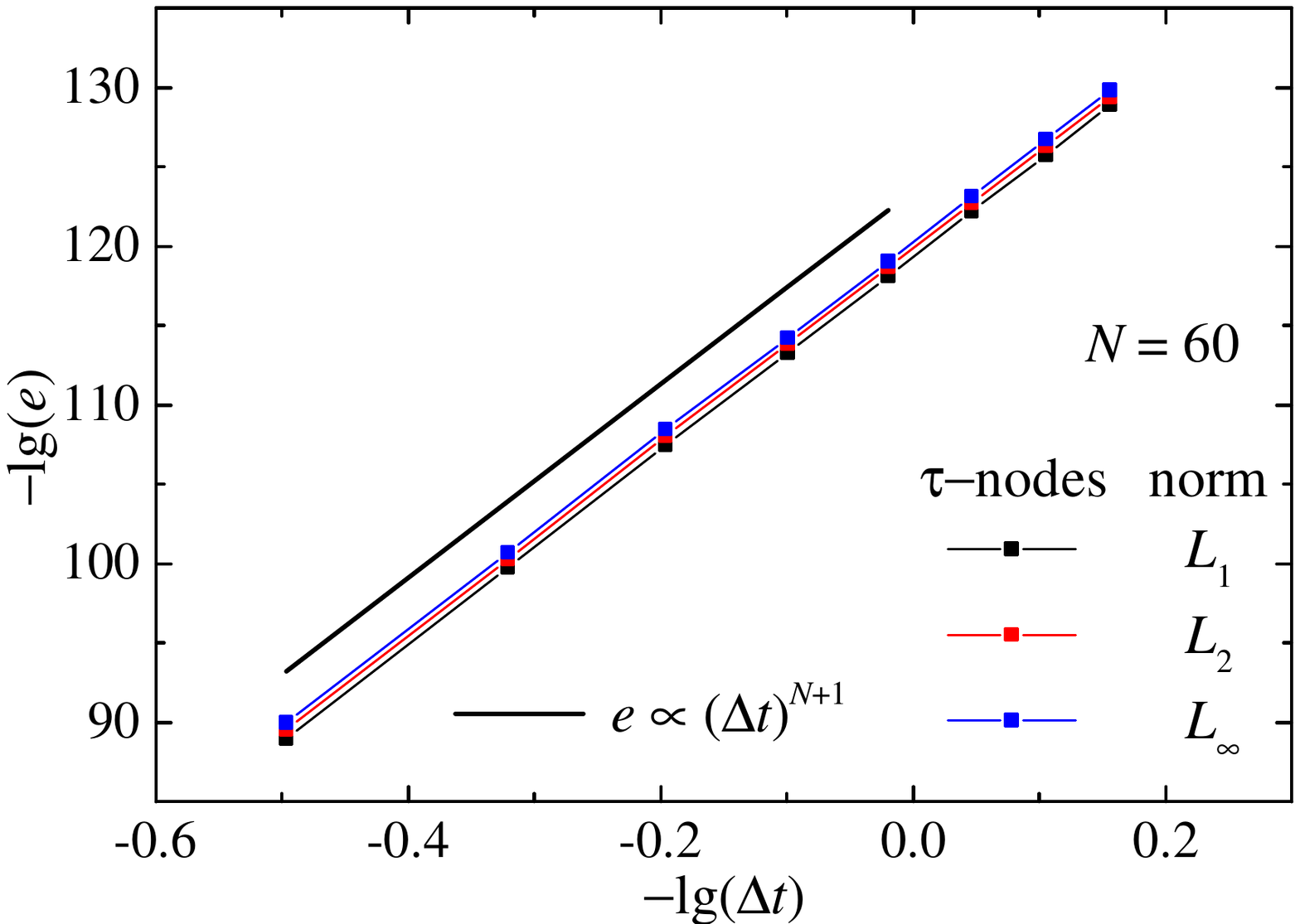}\\
\caption{\label{fig:errs_harm_osc_qnodes}%
Dependence of error $e^{l,q}$ on discretization step ${\Delta t}$ for local solution $\mathbf{u}_{L}(t_{n, p})$ at nodal points $\{t_{n}+\tau_{p}{\Delta t}\}$, calculated in norms $\mathcal{L}_{1}$, $\mathcal{L}_{2}$, $\mathcal{L}_{\infty}$, obtained by the ADER-DG method with $N = 8$, $16$, $60$. Lines represent reference slopes for convergence order $p_{\rm L} = N+1$.
}
\end{figure}

\begin{figure}[h!]
\centering
\includegraphics[width=0.28\textwidth]{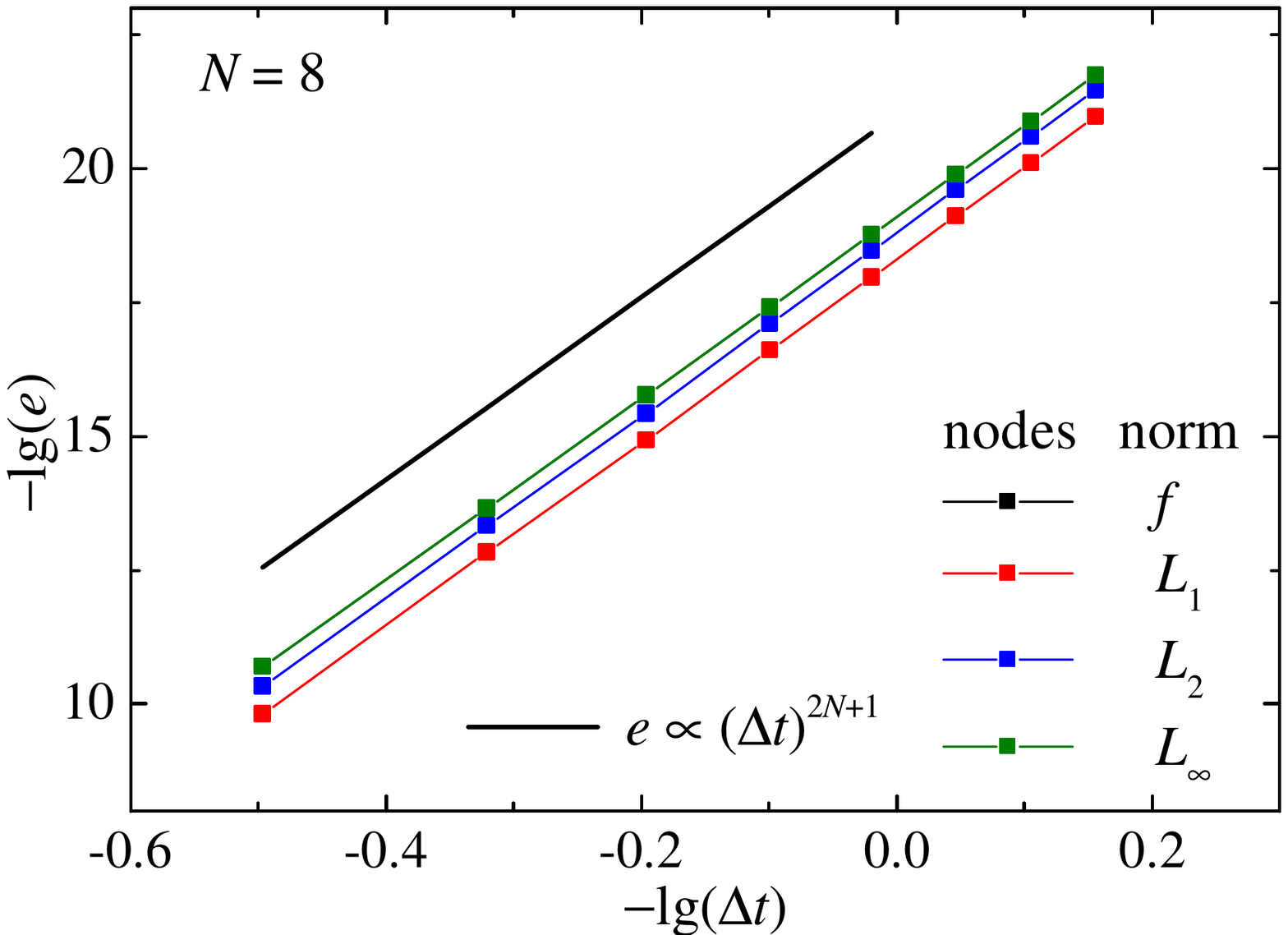}\hspace{4mm}
\includegraphics[width=0.28\textwidth]{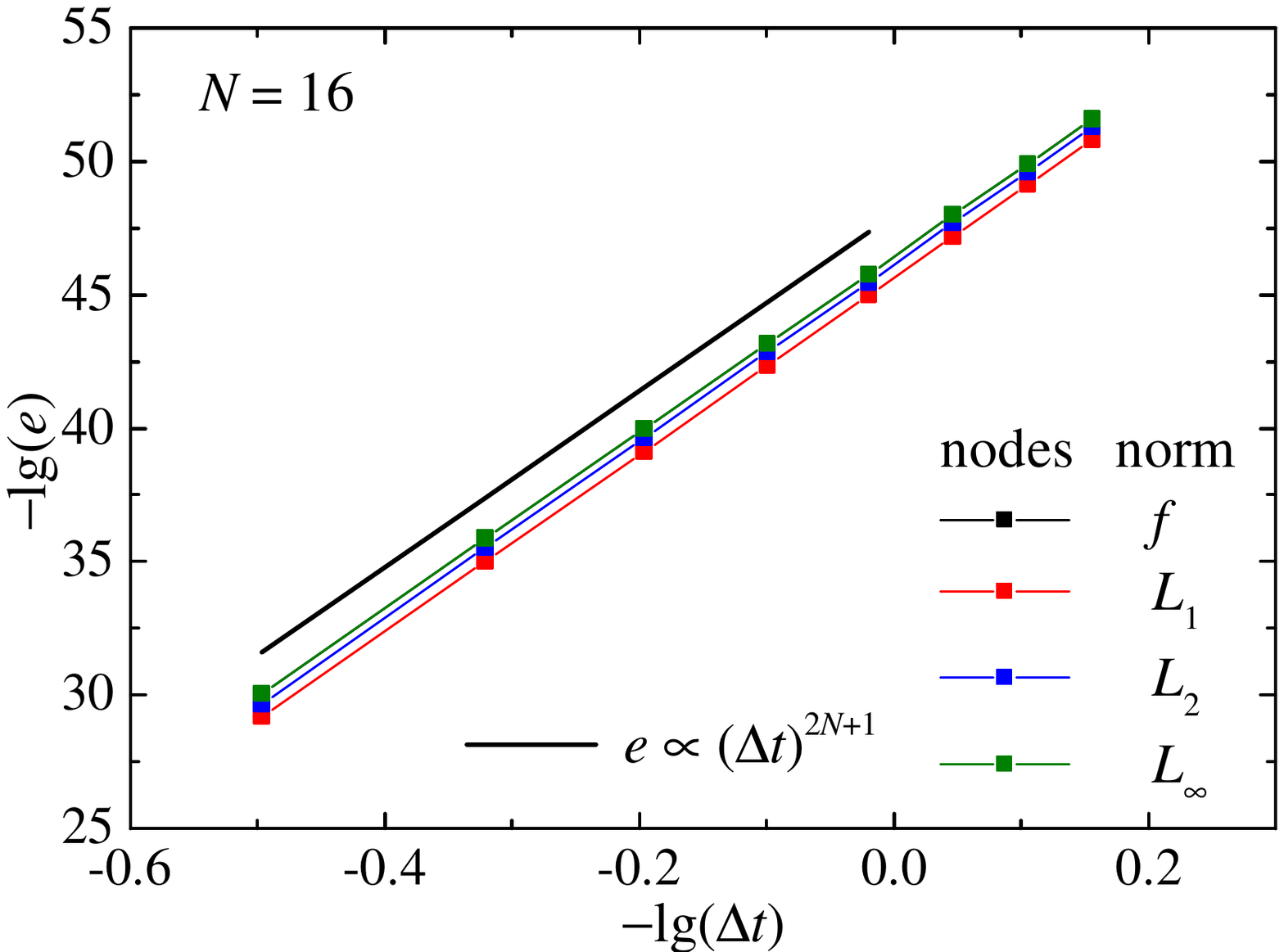}\hspace{4mm}
\includegraphics[width=0.28\textwidth]{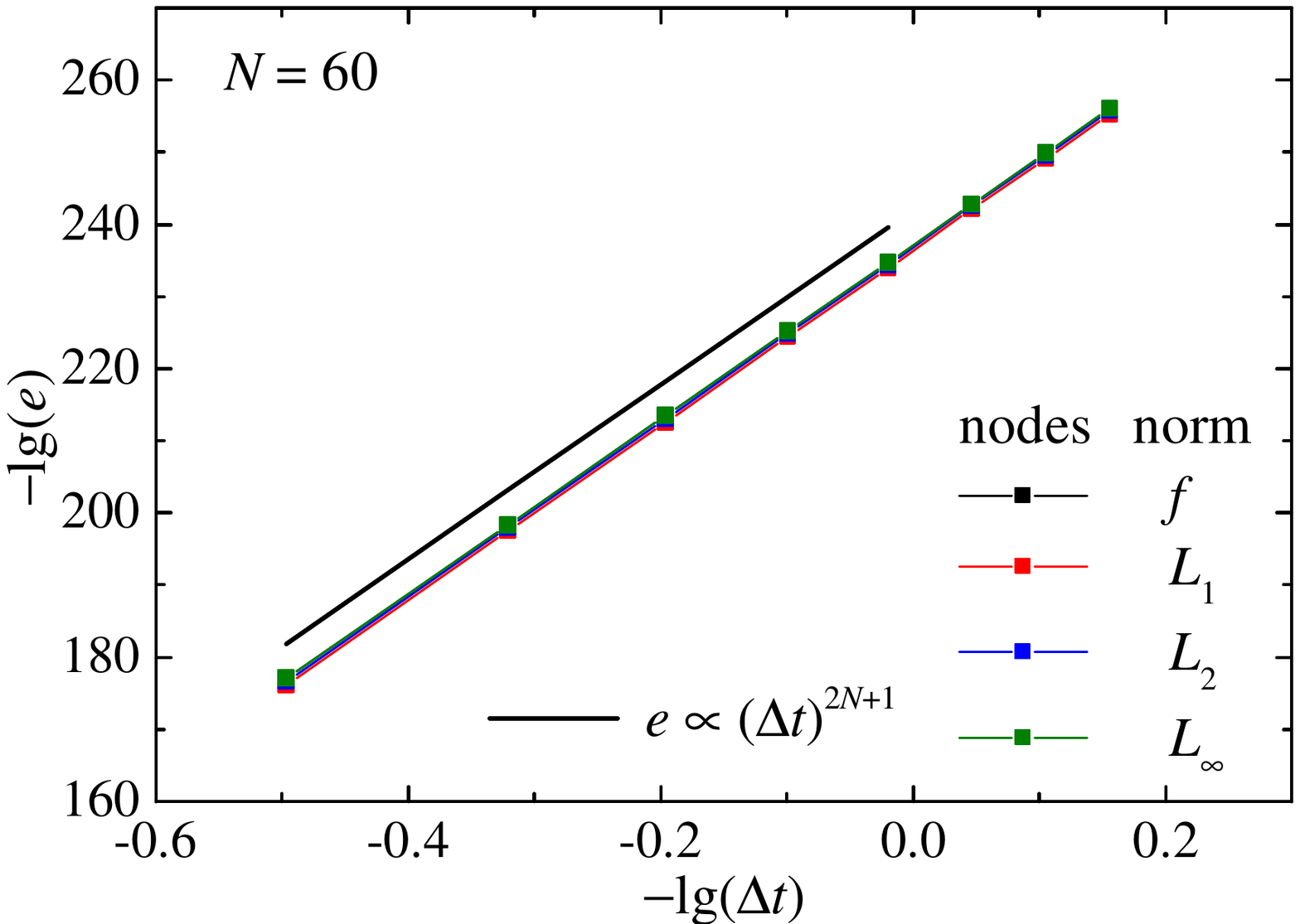}\\
\caption{\label{fig:errs_harm_osc_nodes}%
Dependence of error $e^{n}$ on discretization step ${\Delta t}$ for solution at grid nodes $\mathbf{u}_{n}$, calculated in norms $\mathcal{L}_{1}$, $\mathcal{L}_{2}$, $\mathcal{L}_{\infty}$ and ``final'' norm, obtained by the ADER-DG method with $N = 8$, $16$, $60$. Lines represent reference slopes for convergence order $p_{\rm G} = 2N+1$.
}
\end{figure}

Fig.~\ref{fig:harm_osc_sols} shows a demonstration example of numerical solution of the problem (\ref{eq:harm_osc}) obtained on a grid with only one discretization domain, which demonstrates a very high accuracy of the numerical solution. Figs.~\ref{fig:errs_harm_osc_local},~\ref{fig:errs_harm_osc_qnodes} and~\ref{fig:errs_harm_osc_nodes} show examples of the dependencies of errors $e$ on the discretization step ${\Delta t}$ for $N = 8$, $16$, $60$, with number $M = 4$, $6$, $8$, $10$, $12$, $14$, $16$, $18$ of discretization domains. It is evident that the curves $e({\Delta t})$ in the double logarithmic scale are lines $\lg(e({\Delta t})) \propto \lg({\Delta t})$, and their slopes corresponds well to the expected reference values $p_{\rm G}$ and $p_{\rm L}$.


\begin{table}[htbp]
\centering
\footnotesize
\caption{%
Empirical orders of convergence $p$ of the Runge--Kutta method based on ADER-DG method for the ODE system (\ref{eq:harm_osc}).
}\vspace{-2mm}
\label{tab:conv_orders_harm_osc}
\begin{tabular}{|c|l|lll|l|lll|lll|l|}
\hline
& \multicolumn{5}{c|}{Numerical solution $\mathbf{u}_{n}$ at grid nodes} & \multicolumn{7}{c|}{Local solution $\mathbf{u}_{L}(t)$} \\
\hline
$N$ & $p^{n}_{f}$ & $p^{n}_{L_{1}}$ & $p^{n}_{L_{2}}$ & $p^{n}_{L_{\infty}}$ & $p_{\rm G}$ &
$p^{l}_{L_{1}}$ & $p^{l}_{L_{2}}$ & $p^{l}_{L_{\infty}}$ & 
$p^{l, q}_{L_{1}}$ & $p^{l, q}_{L_{2}}$ & $p^{l, q}_{L_{\infty}}$ & $p_{\rm L}$ \\
\hline
$1$	&	$2.0$	&	$2.3$	&	$2.2$	&	$2.0$	&	$3$	&	$2.0$	&	$2.0$	&	$1.8$	&	$2.1$	&	$2.0$	&	$1.7$	&	$2$\\
$2$	&	$4.6$	&	$4.8$	&	$4.7$	&	$4.6$	&	$5$	&	$3.3$	&	$3.2$	&	$2.8$	&	$3.3$	&	$3.4$	&	$3.6$	&	$3$\\
$3$	&	$6.8$	&	$6.9$	&	$6.9$	&	$6.8$	&	$7$	&	$4.1$	&	$4.0$	&	$3.9$	&	$4.0$	&	$4.0$	&	$4.2$	&	$4$\\
$4$	&	$8.8$	&	$9.0$	&	$9.0$	&	$8.8$	&	$9$	&	$5.0$	&	$5.0$	&	$4.9$	&	$5.0$	&	$5.0$	&	$5.0$	&	$5$\\
$5$	&	$10.9$	&	$11.0$	&	$11.0$	&	$10.9$	&	$11$	&	$6.0$	&	$6.0$	&	$5.9$	&	$6.0$	&	$6.0$	&	$6.0$	&	$6$\\
$6$	&	$12.9$	&	$13.1$	&	$13.0$	&	$12.9$	&	$13$	&	$7.0$	&	$7.0$	&	$6.9$	&	$7.0$	&	$7.0$	&	$6.9$	&	$7$\\
$7$	&	$14.9$	&	$15.1$	&	$15.0$	&	$14.9$	&	$15$	&	$8.0$	&	$8.0$	&	$7.9$	&	$8.0$	&	$8.0$	&	$8.0$	&	$8$\\
$8$	&	$16.9$	&	$17.1$	&	$17.1$	&	$16.9$	&	$17$	&	$9.0$	&	$9.0$	&	$8.9$	&	$9.0$	&	$9.0$	&	$9.0$	&	$9$\\
$9$	&	$18.9$	&	$19.1$	&	$19.1$	&	$18.9$	&	$19$	&	$10.0$	&	$10.0$	&	$9.9$	&	$10.0$	&	$10.0$	&	$10.0$	&	$10$\\
$10$	&	$20.9$	&	$21.1$	&	$21.1$	&	$20.9$	&	$21$	&	$11.0$	&	$11.0$	&	$10.9$	&	$11.0$	&	$11.0$	&	$11.0$	&	$11$\\
$11$	&	$22.9$	&	$23.1$	&	$23.1$	&	$22.9$	&	$23$	&	$12.0$	&	$12.0$	&	$11.9$	&	$12.0$	&	$12.0$	&	$12.0$	&	$12$\\
$12$	&	$24.9$	&	$25.1$	&	$25.1$	&	$24.9$	&	$25$	&	$13.0$	&	$13.0$	&	$13.0$	&	$13.0$	&	$13.0$	&	$13.0$	&	$13$\\
$13$	&	$26.9$	&	$27.1$	&	$27.1$	&	$26.9$	&	$27$	&	$14.0$	&	$14.0$	&	$14.0$	&	$14.0$	&	$14.0$	&	$14.0$	&	$14$\\
$14$	&	$28.9$	&	$29.1$	&	$29.1$	&	$28.9$	&	$29$	&	$15.0$	&	$15.0$	&	$15.0$	&	$15.0$	&	$15.0$	&	$15.0$	&	$15$\\
$15$	&	$31.0$	&	$31.1$	&	$31.1$	&	$31.0$	&	$31$	&	$16.0$	&	$16.0$	&	$16.0$	&	$16.0$	&	$16.0$	&	$16.0$	&	$16$\\
$16$	&	$33.0$	&	$33.1$	&	$33.1$	&	$33.0$	&	$33$	&	$17.0$	&	$17.0$	&	$17.0$	&	$17.0$	&	$17.0$	&	$17.0$	&	$17$\\
$17$	&	$35.0$	&	$35.1$	&	$35.1$	&	$35.0$	&	$35$	&	$18.0$	&	$18.0$	&	$18.0$	&	$18.0$	&	$18.0$	&	$18.0$	&	$18$\\
$18$	&	$37.0$	&	$37.1$	&	$37.1$	&	$37.0$	&	$37$	&	$19.0$	&	$19.0$	&	$19.0$	&	$19.0$	&	$19.0$	&	$19.0$	&	$19$\\
$19$	&	$39.0$	&	$39.1$	&	$39.1$	&	$39.0$	&	$39$	&	$20.0$	&	$20.0$	&	$20.0$	&	$20.0$	&	$20.0$	&	$20.0$	&	$20$\\
$20$	&	$41.0$	&	$41.1$	&	$41.1$	&	$41.0$	&	$41$	&	$21.0$	&	$21.0$	&	$21.0$	&	$21.0$	&	$21.0$	&	$21.0$	&	$21$\\
$21$	&	$43.0$	&	$43.1$	&	$43.1$	&	$43.0$	&	$43$	&	$22.0$	&	$22.0$	&	$22.0$	&	$22.0$	&	$22.0$	&	$22.0$	&	$22$\\
$22$	&	$45.0$	&	$45.2$	&	$45.1$	&	$45.0$	&	$45$	&	$23.0$	&	$23.0$	&	$23.0$	&	$23.0$	&	$23.0$	&	$23.0$	&	$23$\\
$23$	&	$47.0$	&	$47.2$	&	$47.1$	&	$47.0$	&	$47$	&	$24.0$	&	$24.0$	&	$24.0$	&	$24.0$	&	$24.0$	&	$24.0$	&	$24$\\
$24$	&	$49.0$	&	$49.2$	&	$49.1$	&	$49.0$	&	$49$	&	$25.0$	&	$25.0$	&	$25.0$	&	$25.0$	&	$25.0$	&	$25.0$	&	$25$\\
$25$	&	$51.0$	&	$51.2$	&	$51.1$	&	$51.0$	&	$51$	&	$26.0$	&	$26.0$	&	$26.0$	&	$26.0$	&	$26.0$	&	$26.0$	&	$26$\\
$30$	&	$61.0$	&	$61.2$	&	$61.1$	&	$61.0$	&	$61$	&	$31.0$	&	$31.0$	&	$31.0$	&	$31.0$	&	$31.0$	&	$31.0$	&	$31$\\
$35$	&	$71.0$	&	$71.2$	&	$71.1$	&	$71.0$	&	$71$	&	$36.0$	&	$36.0$	&	$36.0$	&	$36.0$	&	$36.0$	&	$36.0$	&	$36$\\
$40$	&	$81.0$	&	$81.2$	&	$81.1$	&	$81.0$	&	$81$	&	$41.0$	&	$41.0$	&	$41.0$	&	$41.0$	&	$41.0$	&	$41.0$	&	$41$\\
$45$	&	$91.0$	&	$91.2$	&	$91.1$	&	$91.0$	&	$91$	&	$46.0$	&	$46.0$	&	$46.0$	&	$46.0$	&	$46.0$	&	$46.0$	&	$46$\\
$50$	&	$101.0$	&	$101.2$	&	$101.1$	&	$101.0$	&	$101$	&	$51.0$	&	$51.0$	&	$51.0$	&	$51.0$	&	$51.0$	&	$51.0$	&	$51$\\
$55$	&	$111.0$	&	$111.2$	&	$111.1$	&	$111.0$	&	$111$	&	$56.0$	&	$56.0$	&	$56.0$	&	$56.0$	&	$56.0$	&	$56.0$	&	$56$\\
$60$	&	$121.0$	&	$121.2$	&	$121.1$	&	$121.0$	&	$121$	&	$61.0$	&	$61.0$	&	$61.0$	&	$61.0$	&	$61.0$	&	$61.0$	&	$61$\\
\hline
\end{tabular}
\end{table}

\begin{table}[htbp]
\centering
\footnotesize
\caption{%
Empirical orders of convergence $p$ of the Runge--Kutta method based on ADER-DG method for the ODE system (\ref{eq:pend}).
}\vspace{-2mm}
\label{tab:conv_orders_pend}
\begin{tabular}{|c|l|lll|l|lll|lll|l|}
\hline
& \multicolumn{5}{c|}{Numerical solution $\mathbf{u}_{n}$ at grid nodes} & \multicolumn{7}{c|}{Local solution $\mathbf{u}_{L}(t)$} \\
\hline
$N$ & $p^{n}_{f}$ & $p^{n}_{L_{1}}$ & $p^{n}_{L_{2}}$ & $p^{n}_{L_{\infty}}$ & $p_{\rm G}$ &
$p^{l}_{L_{1}}$ & $p^{l}_{L_{2}}$ & $p^{l}_{L_{\infty}}$ & 
$p^{l, q}_{L_{1}}$ & $p^{l, q}_{L_{2}}$ & $p^{l, q}_{L_{\infty}}$ & $p_{\rm L}$ \\
\hline
$1$	&	$2.8$	&	$2.9$	&	$2.9$	&	$2.7$	&	$3$	&	$2.4$	&	$2.4$	&	$2.1$	&	$2.5$	&	$2.5$	&	$2.5$	&	$2$\\
$2$	&	$4.8$	&	$4.9$	&	$4.8$	&	$4.8$	&	$5$	&	$3.1$	&	$3.0$	&	$2.9$	&	$3.1$	&	$3.0$	&	$3.1$	&	$3$\\
$3$	&	$6.8$	&	$7.0$	&	$6.9$	&	$6.8$	&	$7$	&	$4.0$	&	$4.0$	&	$3.9$	&	$4.0$	&	$4.0$	&	$4.0$	&	$4$\\
$4$	&	$8.7$	&	$8.7$	&	$8.7$	&	$8.6$	&	$9$	&	$4.9$	&	$4.9$	&	$4.8$	&	$4.9$	&	$4.9$	&	$4.9$	&	$5$\\
$5$	&	$10.8$	&	$11.0$	&	$11.0$	&	$10.9$	&	$11$	&	$6.0$	&	$5.9$	&	$5.8$	&	$6.0$	&	$6.0$	&	$5.9$	&	$6$\\
$6$	&	$12.6$	&	$12.7$	&	$12.6$	&	$12.6$	&	$13$	&	$6.9$	&	$6.9$	&	$6.9$	&	$6.9$	&	$6.9$	&	$6.9$	&	$7$\\
$7$	&	$14.7$	&	$14.8$	&	$14.8$	&	$14.7$	&	$15$	&	$7.8$	&	$7.8$	&	$7.7$	&	$7.8$	&	$7.9$	&	$7.9$	&	$8$\\
$8$	&	$16.6$	&	$16.7$	&	$16.7$	&	$16.6$	&	$17$	&	$9.0$	&	$9.0$	&	$8.8$	&	$9.0$	&	$9.0$	&	$8.9$	&	$9$\\
$9$	&	$18.6$	&	$18.7$	&	$18.7$	&	$18.6$	&	$19$	&	$9.7$	&	$9.7$	&	$9.6$	&	$9.7$	&	$9.8$	&	$9.8$	&	$10$\\
$10$	&	$20.4$	&	$20.6$	&	$20.6$	&	$20.4$	&	$21$	&	$11.1$	&	$11.1$	&	$10.9$	&	$11.0$	&	$11.0$	&	$10.9$	&	$11$\\
$11$	&	$22.8$	&	$22.8$	&	$22.9$	&	$22.8$	&	$23$	&	$11.6$	&	$11.6$	&	$11.5$	&	$11.7$	&	$11.7$	&	$11.7$	&	$12$\\
$12$	&	$24.0$	&	$24.3$	&	$24.2$	&	$24.1$	&	$25$	&	$13.1$	&	$13.1$	&	$12.9$	&	$13.0$	&	$12.9$	&	$12.8$	&	$13$\\
$13$	&	$27.1$	&	$27.1$	&	$27.2$	&	$27.1$	&	$27$	&	$13.6$	&	$13.6$	&	$13.5$	&	$13.8$	&	$13.8$	&	$13.7$	&	$14$\\
$14$	&	$27.3$	&	$27.6$	&	$27.5$	&	$27.3$	&	$29$	&	$15.0$	&	$14.9$	&	$14.8$	&	$14.8$	&	$14.7$	&	$14.6$	&	$15$\\
$15$	&	$31.3$	&	$31.3$	&	$31.3$	&	$31.2$	&	$31$	&	$15.7$	&	$15.7$	&	$15.5$	&	$15.9$	&	$15.9$	&	$15.7$	&	$16$\\
$16$	&	$30.5$	&	$31.2$	&	$31.1$	&	$30.9$	&	$33$	&	$16.6$	&	$16.7$	&	$16.6$	&	$16.4$	&	$16.4$	&	$16.4$	&	$17$\\
$17$	&	$35.2$	&	$35.4$	&	$35.3$	&	$35.2$	&	$35$	&	$18.0$	&	$17.8$	&	$17.5$	&	$18.1$	&	$17.9$	&	$17.7$	&	$18$\\
$18$	&	$35.4$	&	$35.3$	&	$35.4$	&	$35.4$	&	$37$	&	$18.5$	&	$18.6$	&	$18.4$	&	$18.4$	&	$18.3$	&	$18.1$	&	$19$\\
$19$	&	$39.1$	&	$39.3$	&	$39.3$	&	$39.2$	&	$39$	&	$19.9$	&	$19.8$	&	$19.4$	&	$19.9$	&	$19.9$	&	$19.6$	&	$20$\\
$20$	&	$39.7$	&	$39.4$	&	$39.5$	&	$39.6$	&	$41$	&	$20.7$	&	$20.6$	&	$20.3$	&	$20.6$	&	$20.4$	&	$20.1$	&	$21$\\
$21$	&	$43.0$	&	$43.2$	&	$43.1$	&	$43.0$	&	$43$	&	$21.8$	&	$21.7$	&	$21.4$	&	$21.8$	&	$21.8$	&	$21.6$	&	$22$\\
$22$	&	$43.7$	&	$43.7$	&	$43.6$	&	$43.6$	&	$45$	&	$22.6$	&	$22.6$	&	$22.3$	&	$22.5$	&	$22.5$	&	$22.2$	&	$23$\\
$23$	&	$46.7$	&	$46.9$	&	$46.9$	&	$46.8$	&	$47$	&	$23.7$	&	$23.7$	&	$23.5$	&	$23.8$	&	$23.8$	&	$23.6$	&	$24$\\
$24$	&	$46.9$	&	$47.9$	&	$47.8$	&	$47.6$	&	$49$	&	$24.5$	&	$24.5$	&	$24.3$	&	$24.5$	&	$24.5$	&	$24.3$	&	$25$\\
$25$	&	$50.1$	&	$50.2$	&	$50.3$	&	$50.3$	&	$51$	&	$25.7$	&	$25.7$	&	$25.5$	&	$25.7$	&	$25.7$	&	$25.5$	&	$26$\\
$30$	&	$59.7$	&	$59.8$	&	$59.7$	&	$59.6$	&	$61$	&	$30.5$	&	$30.5$	&	$30.2$	&	$30.6$	&	$30.5$	&	$30.2$	&	$31$\\
$35$	&	$70.5$	&	$71.2$	&	$71.1$	&	$70.7$	&	$71$	&	$35.3$	&	$35.2$	&	$35.0$	&	$35.4$	&	$35.4$	&	$35.2$	&	$36$\\
$40$	&	$79.3$	&	$79.3$	&	$79.3$	&	$79.2$	&	$81$	&	$40.6$	&	$40.6$	&	$40.3$	&	$40.6$	&	$40.5$	&	$40.2$	&	$41$\\
$45$	&	$90.7$	&	$91.2$	&	$91.1$	&	$90.7$	&	$91$	&	$45.2$	&	$44.9$	&	$44.7$	&	$45.2$	&	$45.1$	&	$44.8$	&	$46$\\
$50$	&	$100.1$	&	$100.1$	&	$100.2$	&	$100.1$	&	$101$	&	$50.1$	&	$50.3$	&	$50.1$	&	$50.1$	&	$50.2$	&	$50.1$	&	$51$\\
$55$	&	$109.2$	&	$109.6$	&	$109.6$	&	$109.6$	&	$111$	&	$54.9$	&	$54.8$	&	$54.5$	&	$54.4$	&	$54.4$	&	$54.3$	&	$56$\\
$60$	&	$120.7$	&	$120.6$	&	$120.7$	&	$120.7$	&	$121$	&	$60.7$	&	$60.6$	&	$60.2$	&	$61.1$	&	$61.0$	&	$60.5$	&	$61$\\
\hline
\end{tabular}
\end{table}

Details of the empirical convergence orders $p$ calculated for problem (\ref{eq:harm_osc}) are presented in Table~\ref{tab:conv_orders_harm_osc} for $N = 1,\, \ldots,\, 60$. Comparison of the empirical convergence orders $p$ with expected values $p_{\rm G}$ and $p_{\rm L}$ shows very good agreement, which demonstrates the applicability of the theory of the ADER-DG method developed in this paper. Similar data on the empirical convergence orders calculated for problem (\ref{eq:pend}) are presented in Table~\ref{tab:conv_orders_pend}. These empirical convergence orders $p$ also demonstrate a very good agreement with the expected values $p_{\rm G}$ and $p_{\rm L}$.

\begin{figure}[h!]
\centering
\includegraphics[width=0.24\textwidth]{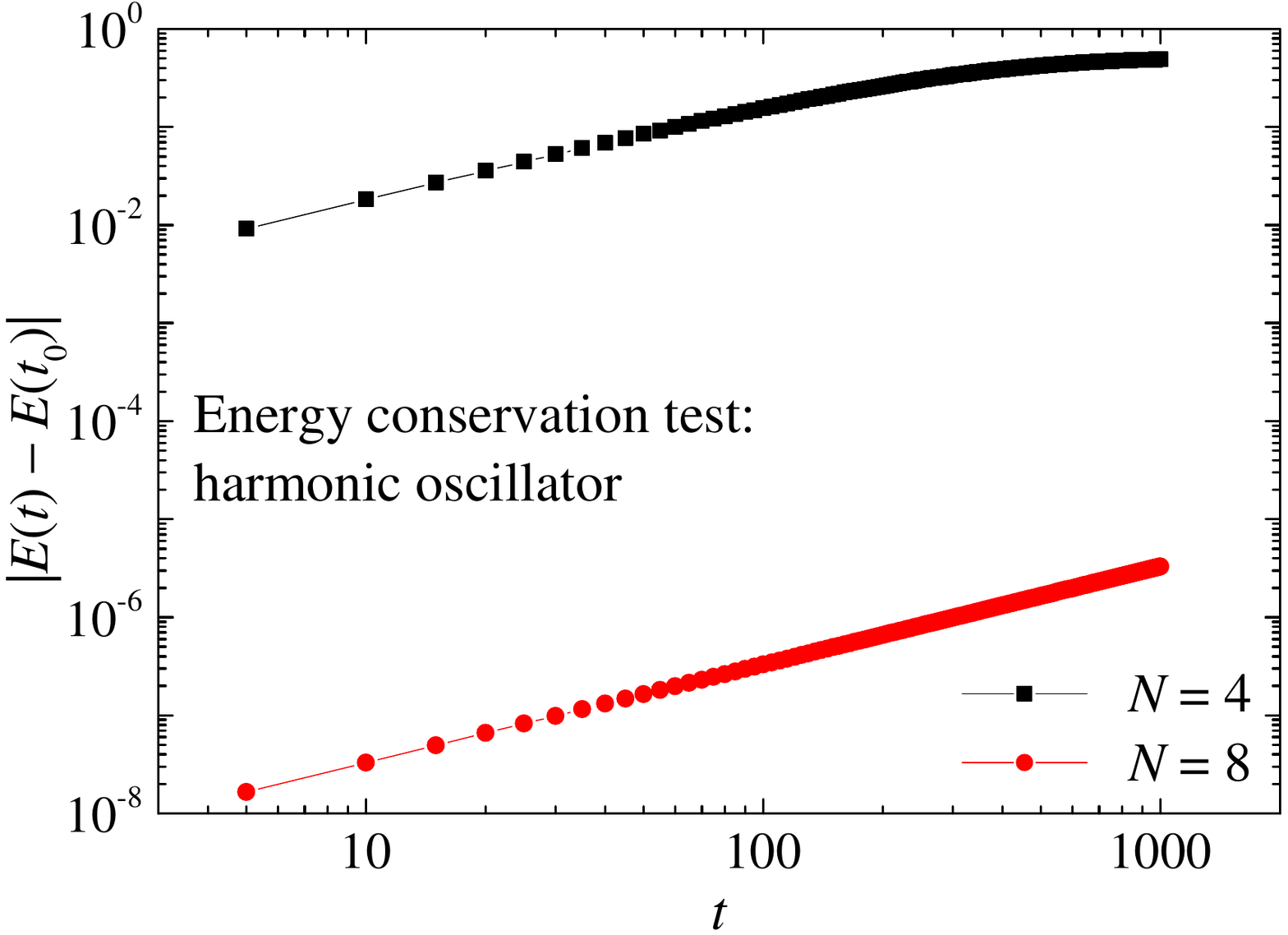}
\includegraphics[width=0.24\textwidth]{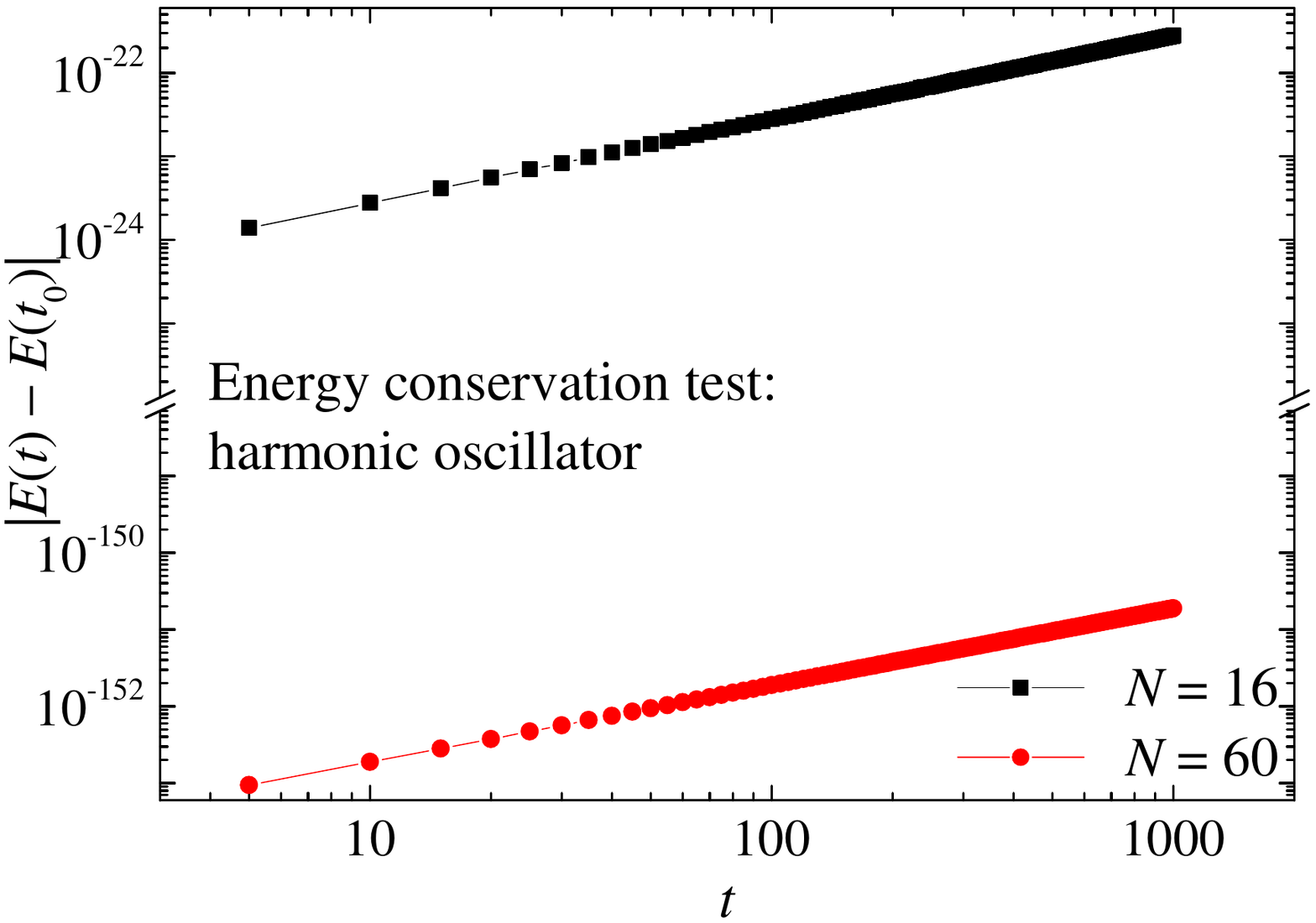}
\includegraphics[width=0.24\textwidth]{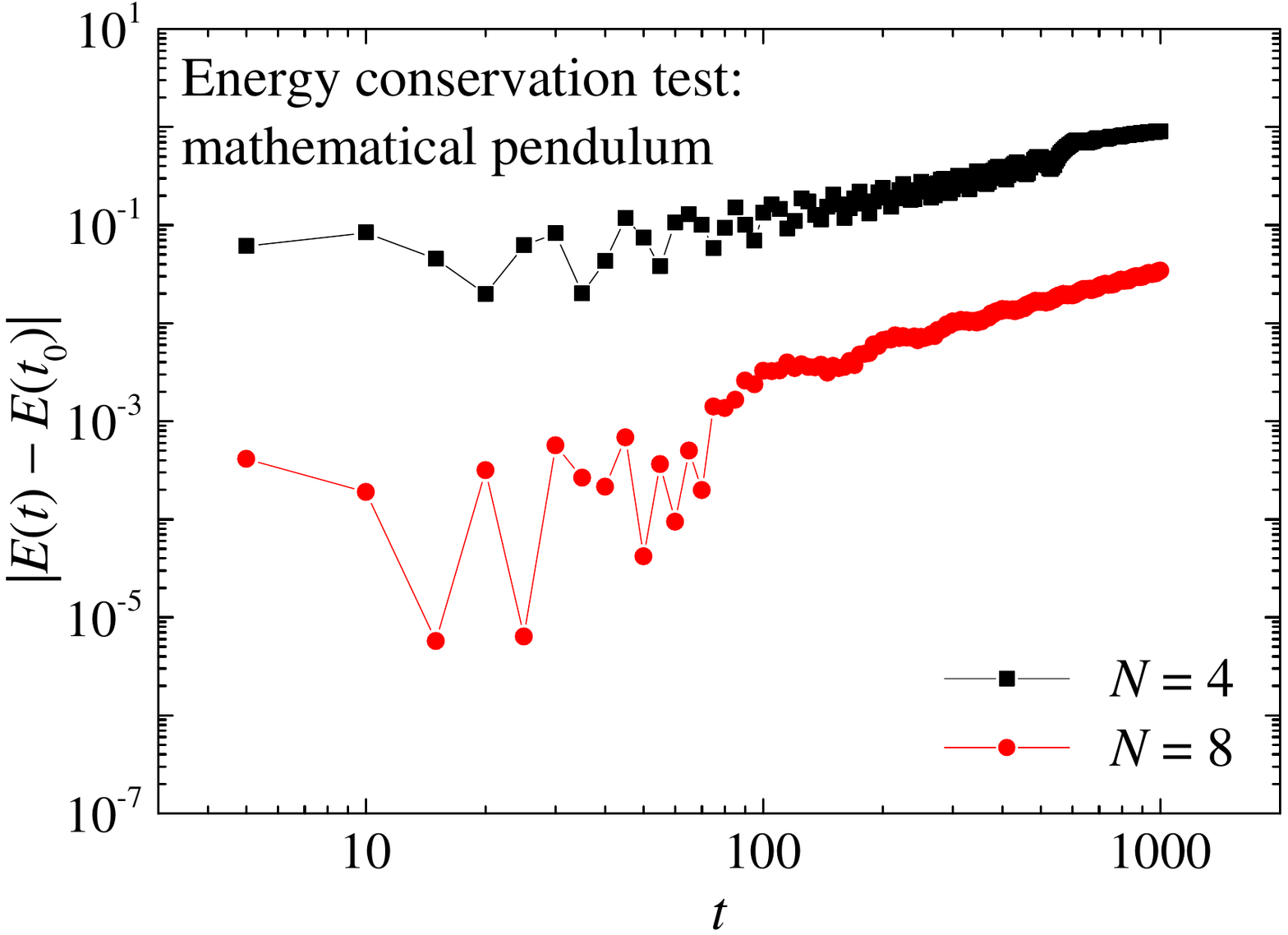}
\includegraphics[width=0.24\textwidth]{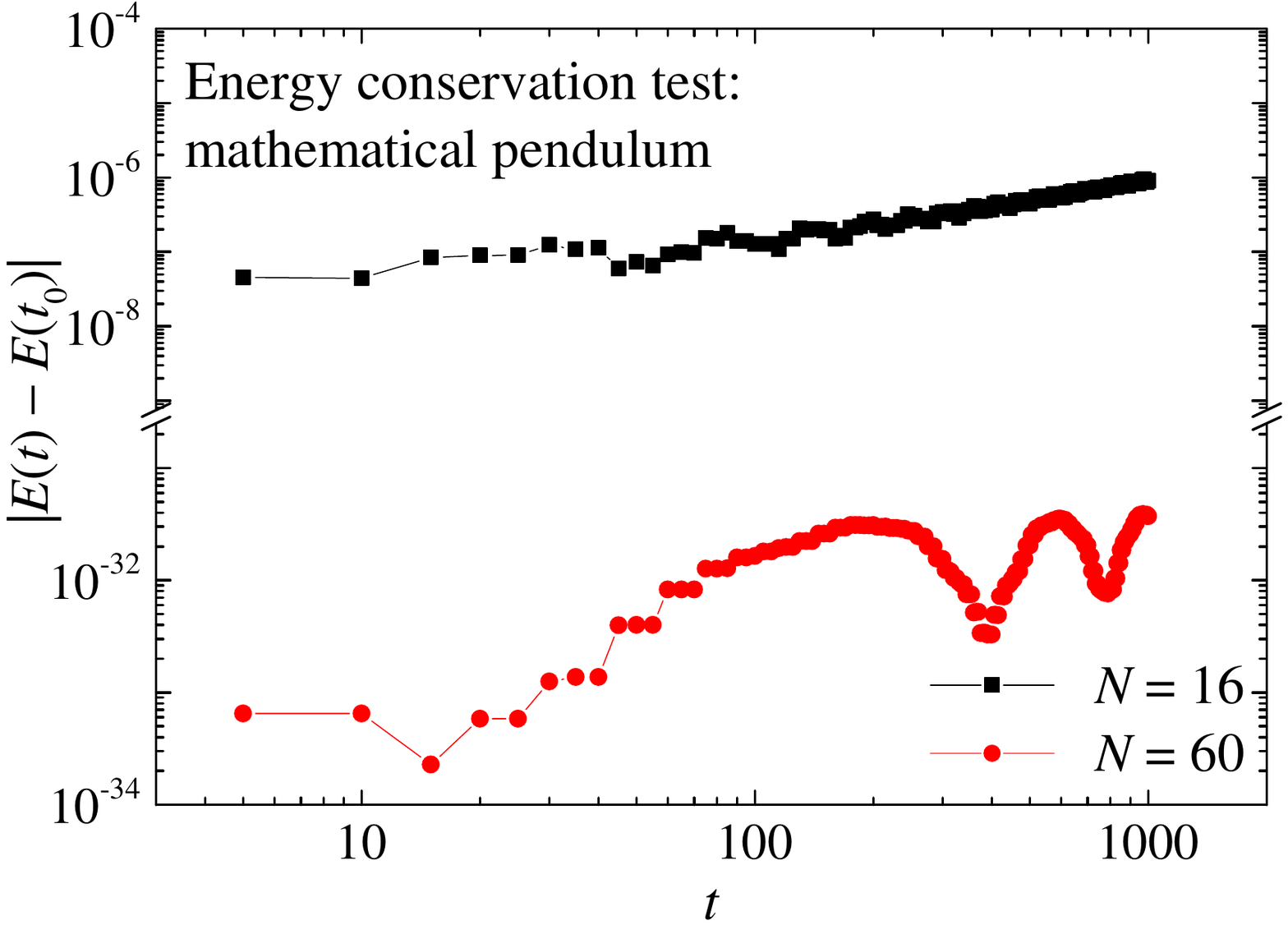}\\
\caption{\label{fig:e_errs}%
Energy conservation error $|E(t)-E(t_{0})|$ in the numerical solution obtained using the ADER-DG method with $N = 4$, $8$, $16$, $60$ for problems (\ref{eq:harm_osc}) and (\ref{eq:pend}) for a wide range $t\in[0, 1000]$.
}
\end{figure}

\corrtext{The ADER-DG numerical method is not symplectic method (the matrix $\mathrm{M}\not\equiv0$ (\ref{eq:q_m_matrices_def}) due to Lemma~\ref{lemma:q_m_matrices_posdef}), so the energy conservation law $E(t) = \mathrm{const}$ or its discrete analog is not expected to hold. However, due to the very high accuracy of the numerical solution, it can be expected that the energy conservation error can be made very small, at least comparable to or less than the precision of representing real numbers as double-precision floating-point numbers. To demonstrate this property, a numerical solution of problems (\ref{eq:harm_osc}) and (\ref{eq:pend}) is obtained for a wide range $t\in[0, 1000]$. Fig.~\ref{fig:e_errs} shows the dependencies of the energy conservation error $|E(t)-E(t_{0})|$ on $t$, calculated with a large discretization step ${\Delta t} = 5$, for $N = 4$, $8$, $16$ and $60$. The obtained results clearly demonstrate that the error in the energy conservation law can be very small, including less than the precision of double-precision floating-point numbers.}

\corrtext{The stability properties are verified by solving four ODEs. Linear $A$- and $AN$-stability are verified by solving the following two ODEs:
\begin{align}
&\dot{u} = -100 u,\hspace{4mm} u(0) = 1,\quad t\in[0,\, 5],\label{eq:a_test}\\
&\dot{u} = -\frac{100 u}{1+t},\hspace{3mm} u(0) = 1,\quad t\in[0,\, 20],\label{eq:an_test}
\end{align}
Nonlinear $B$- and $BN$-stability are verified by solving the two contractive ODEs, for which two different initial conditions are chosen, corresponding to the two solutions $u(t)$ and $v(t)$:
\begin{align}
&\dot{u} = -u^{3},\hspace{7.5mm} u(0) = 1,\qquad \dot{v} = -v^{3},\hspace{7.5mm} v(0) = \frac12,\qquad t\in[0,\, 100],\label{eq:b_test}\\
&\dot{u} = -tu^{3},\hspace{6.2mm} u(0) = 1,\qquad \dot{v} = -tv^{3},\hspace{6.2mm} v(0) = \frac12,\qquad t\in[0,\, 100].\label{eq:bn_test}
\end{align}
Exact analytical solutions of these ODEs are obtained by trivial integration.}

\begin{figure}[h!]
\centering
\includegraphics[width=0.24\textwidth]{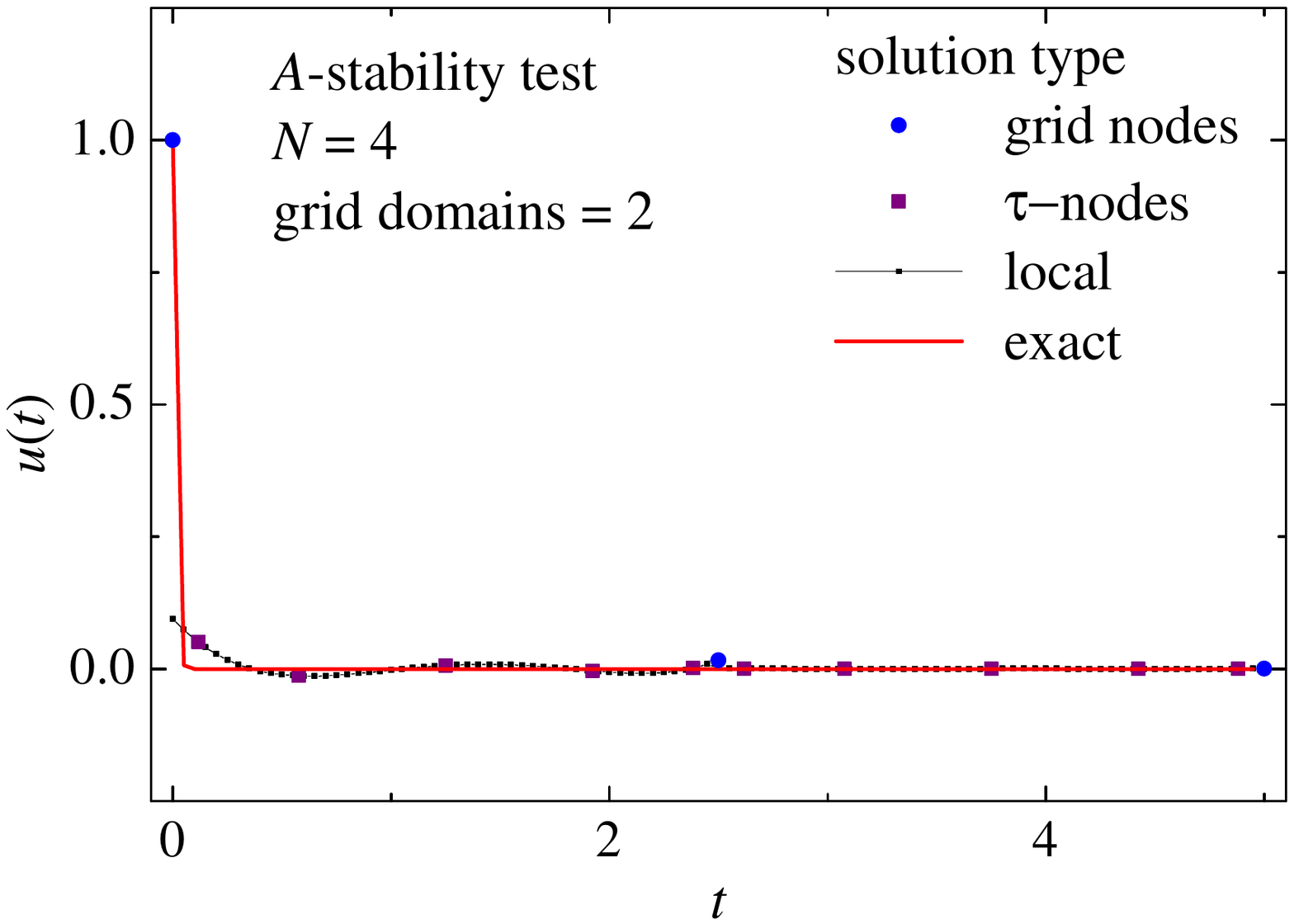}
\includegraphics[width=0.24\textwidth]{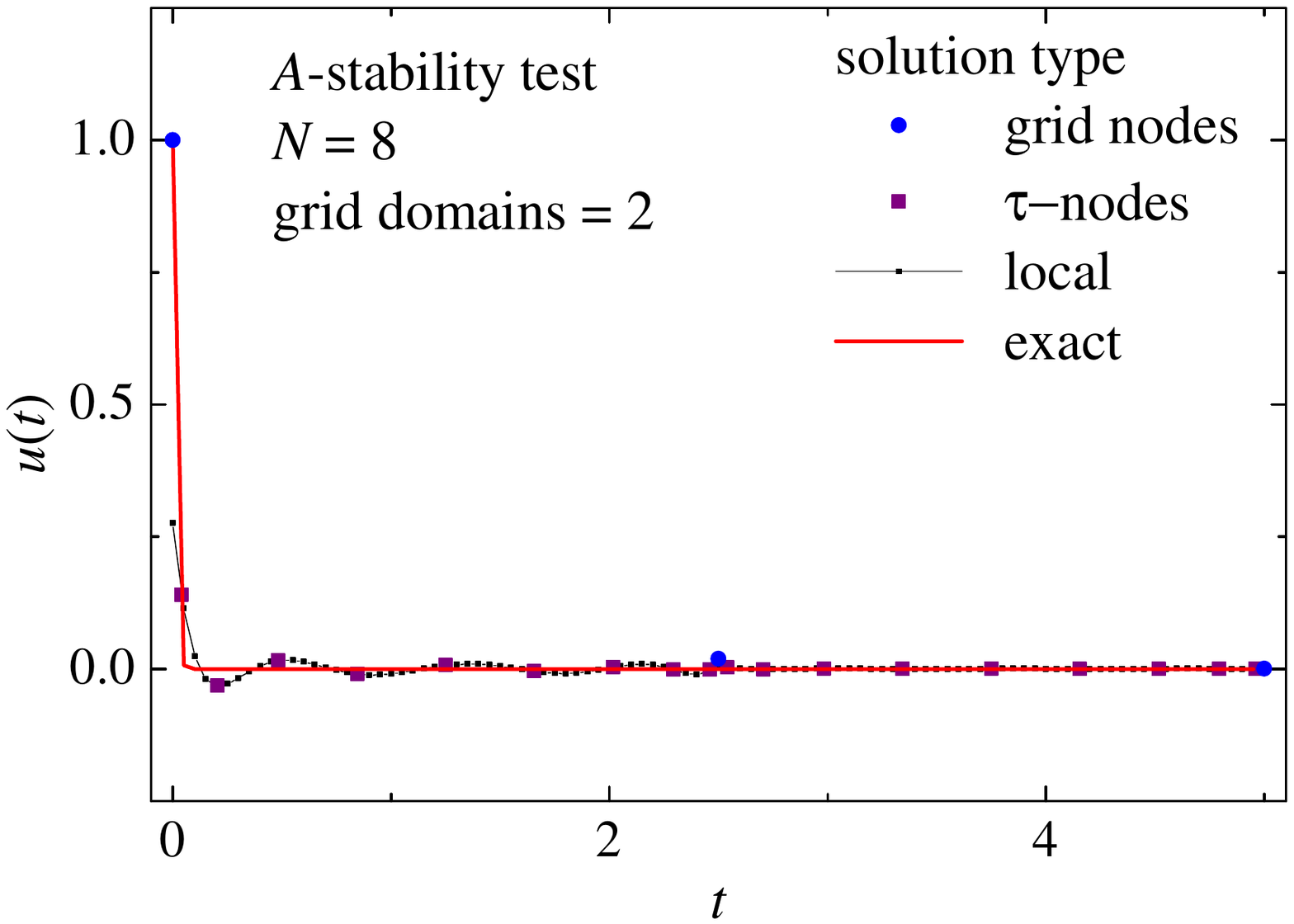}
\includegraphics[width=0.24\textwidth]{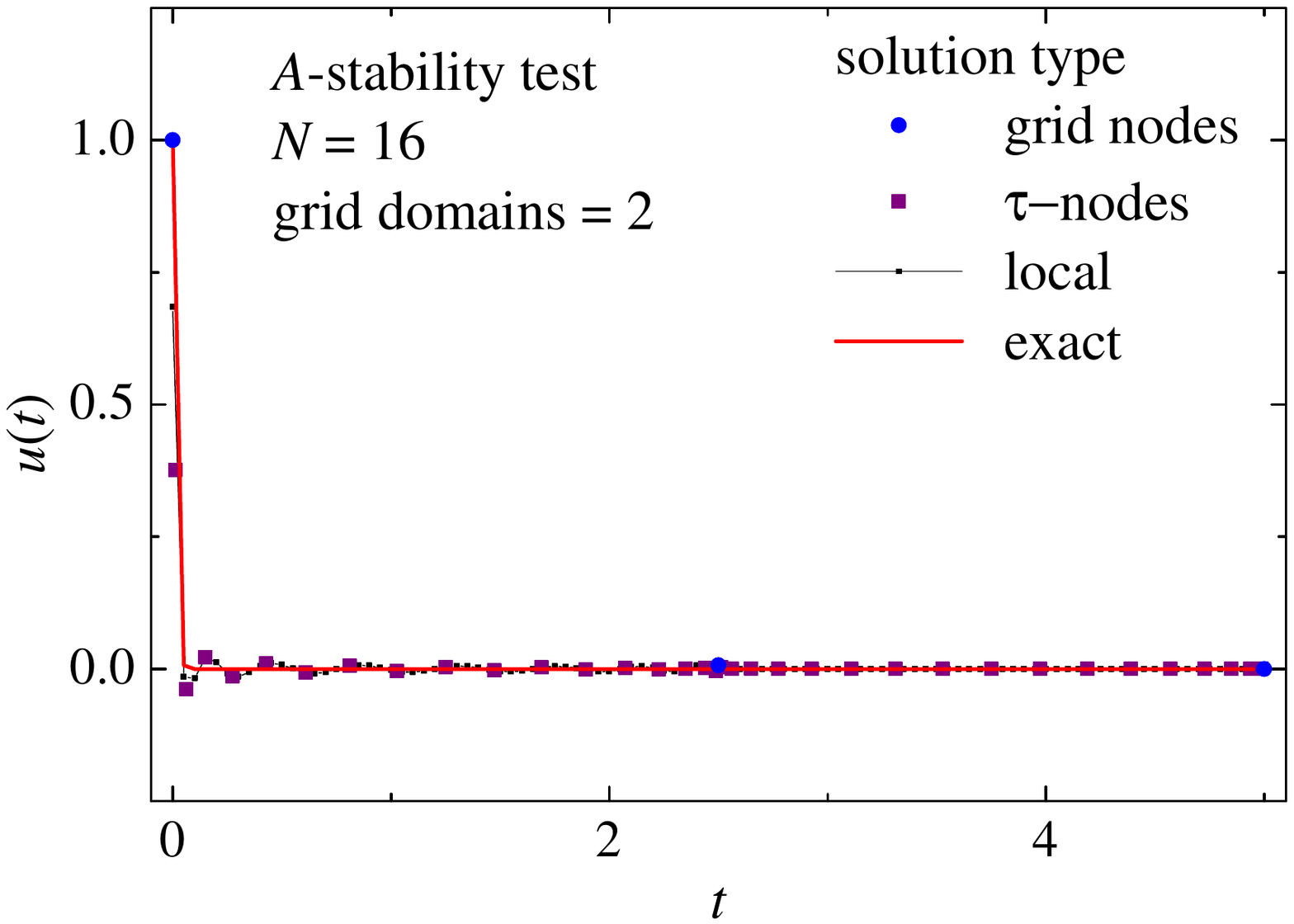}
\includegraphics[width=0.24\textwidth]{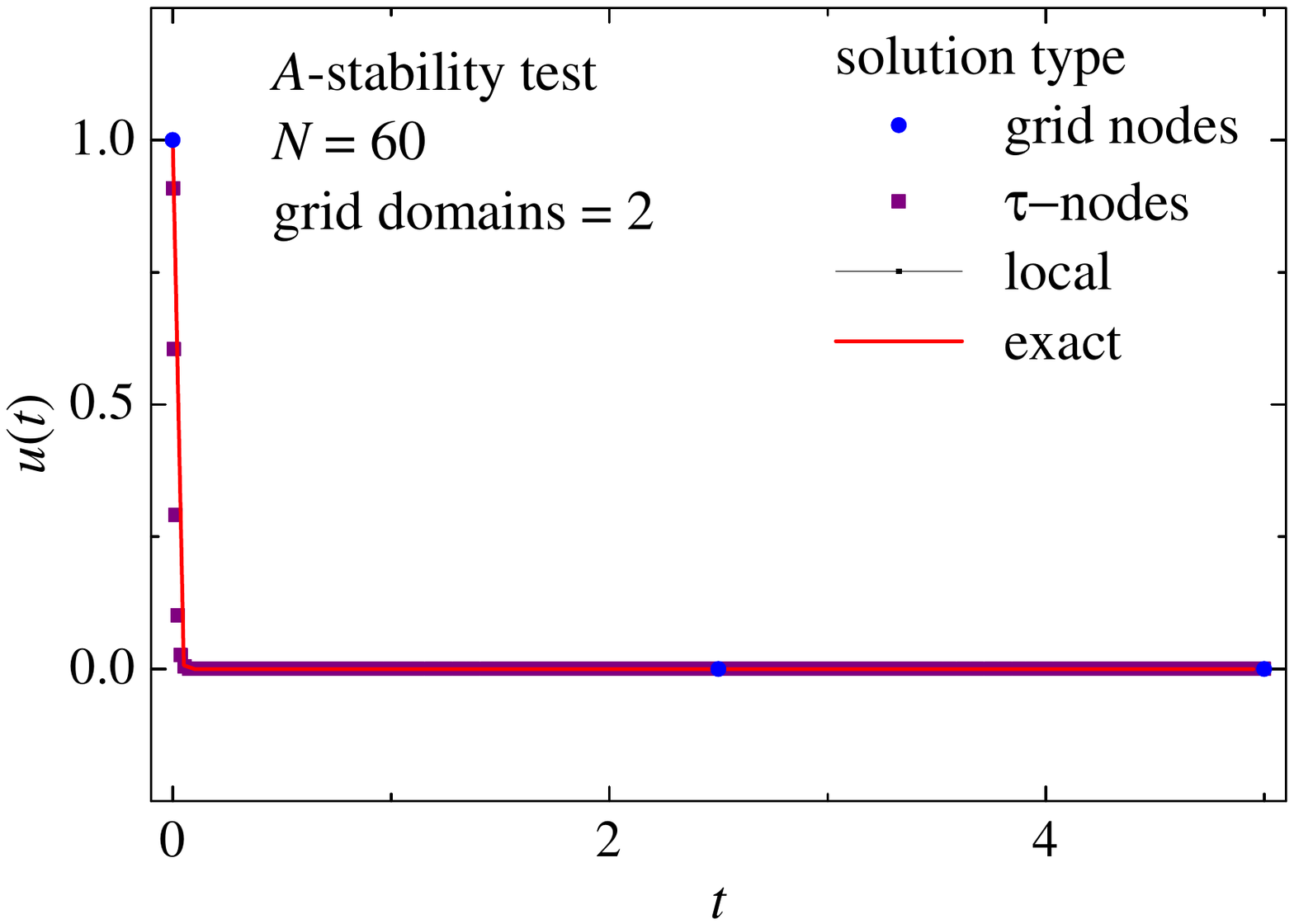}\\
\includegraphics[width=0.24\textwidth]{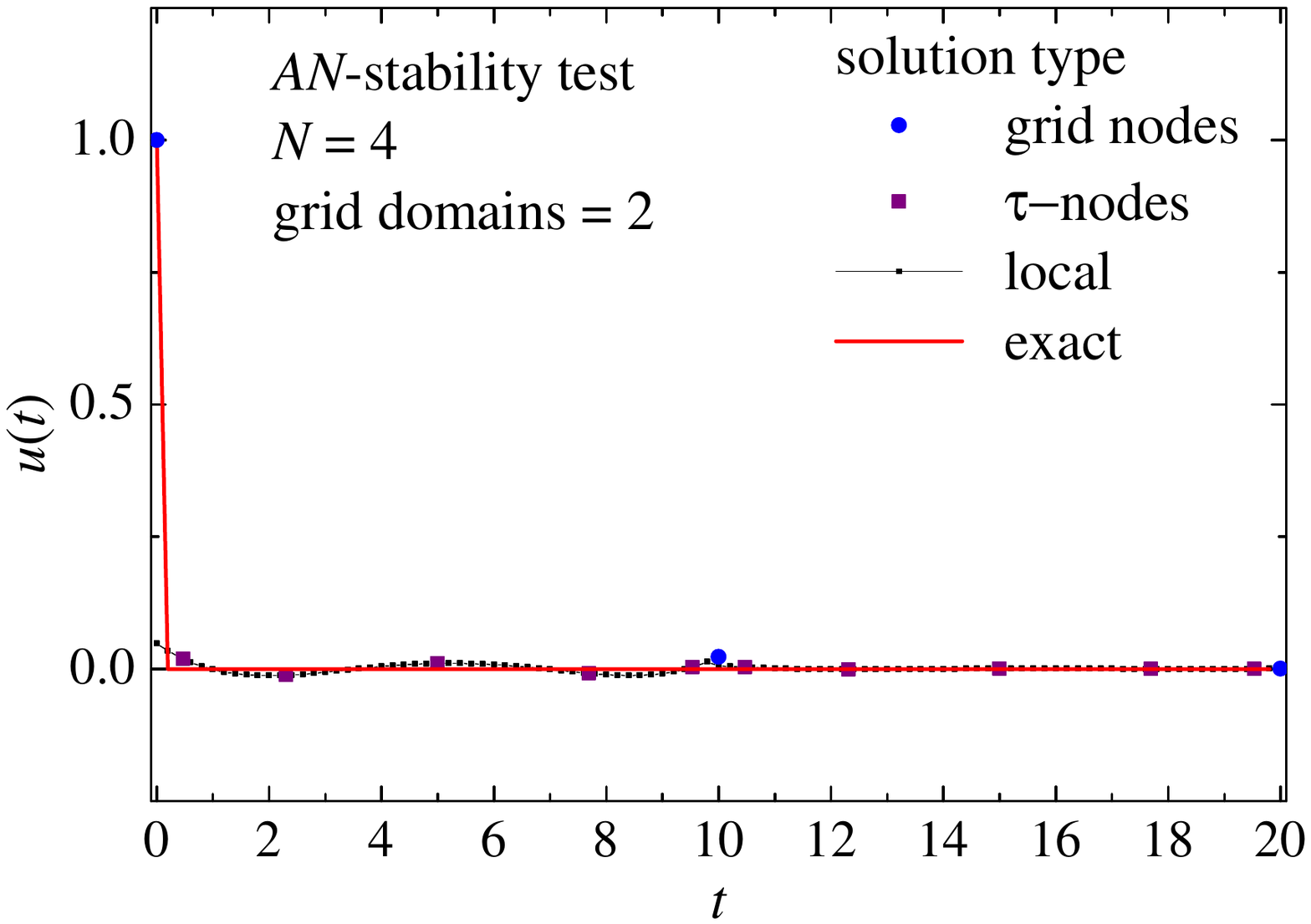}
\includegraphics[width=0.24\textwidth]{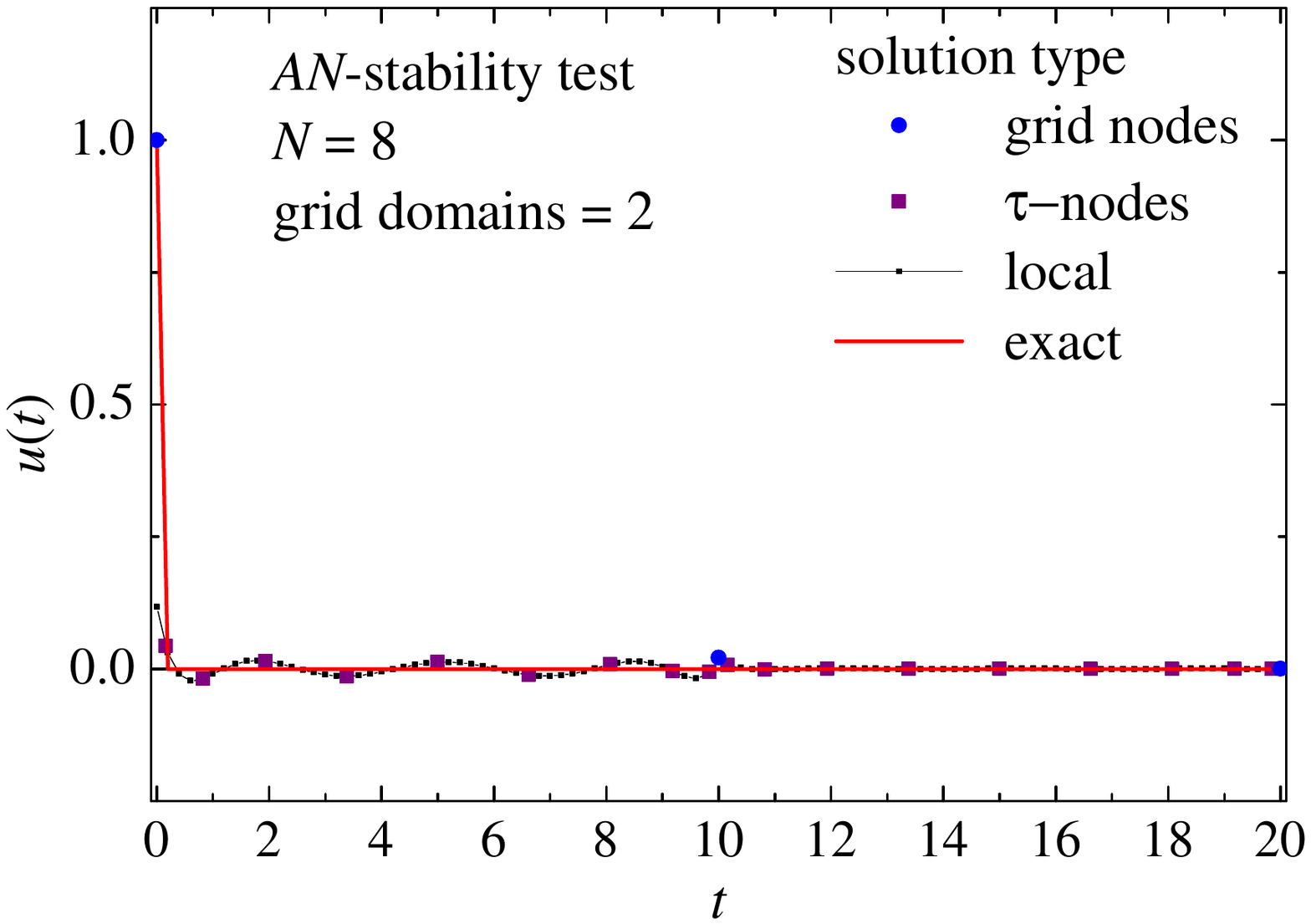}
\includegraphics[width=0.24\textwidth]{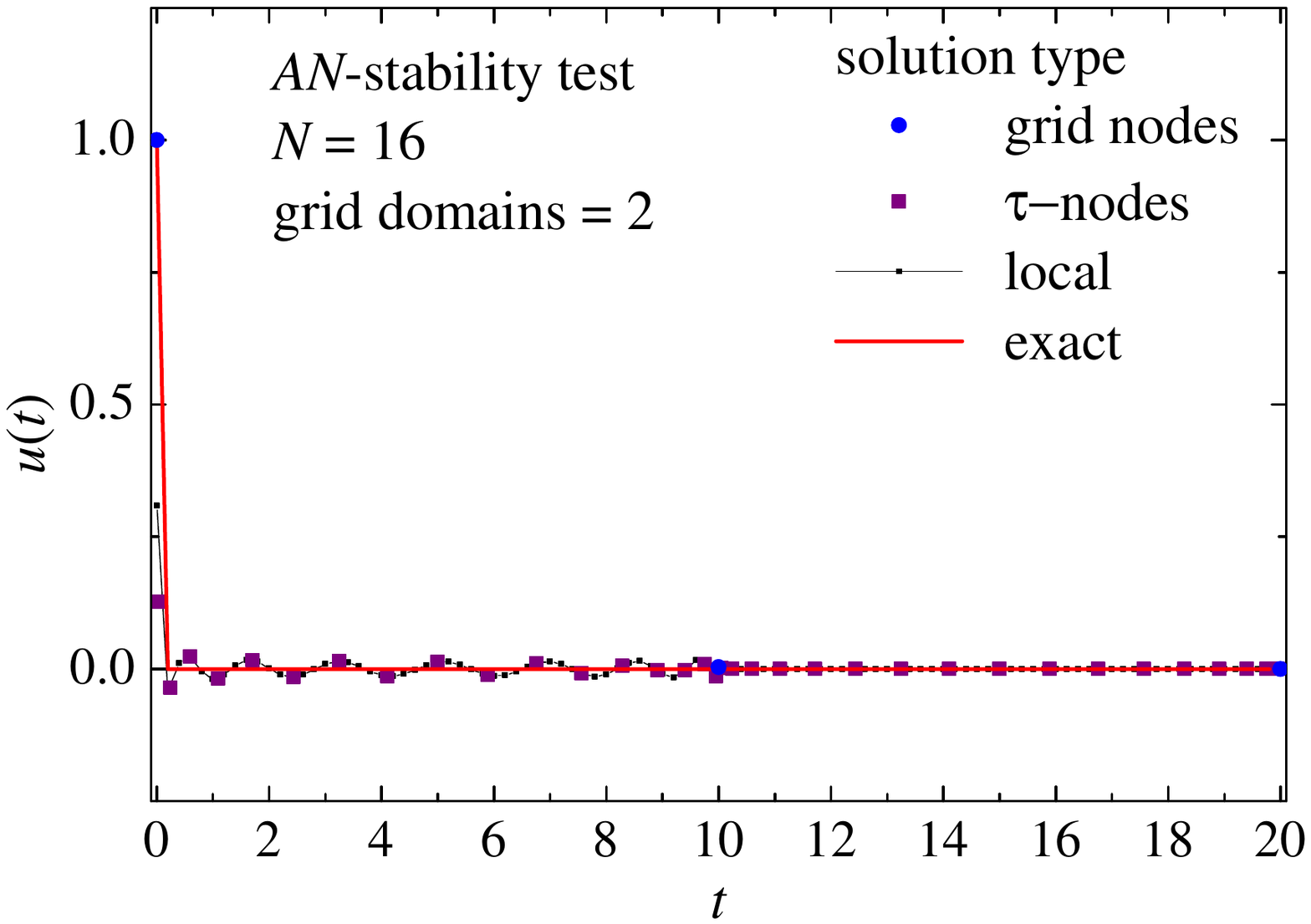}
\includegraphics[width=0.24\textwidth]{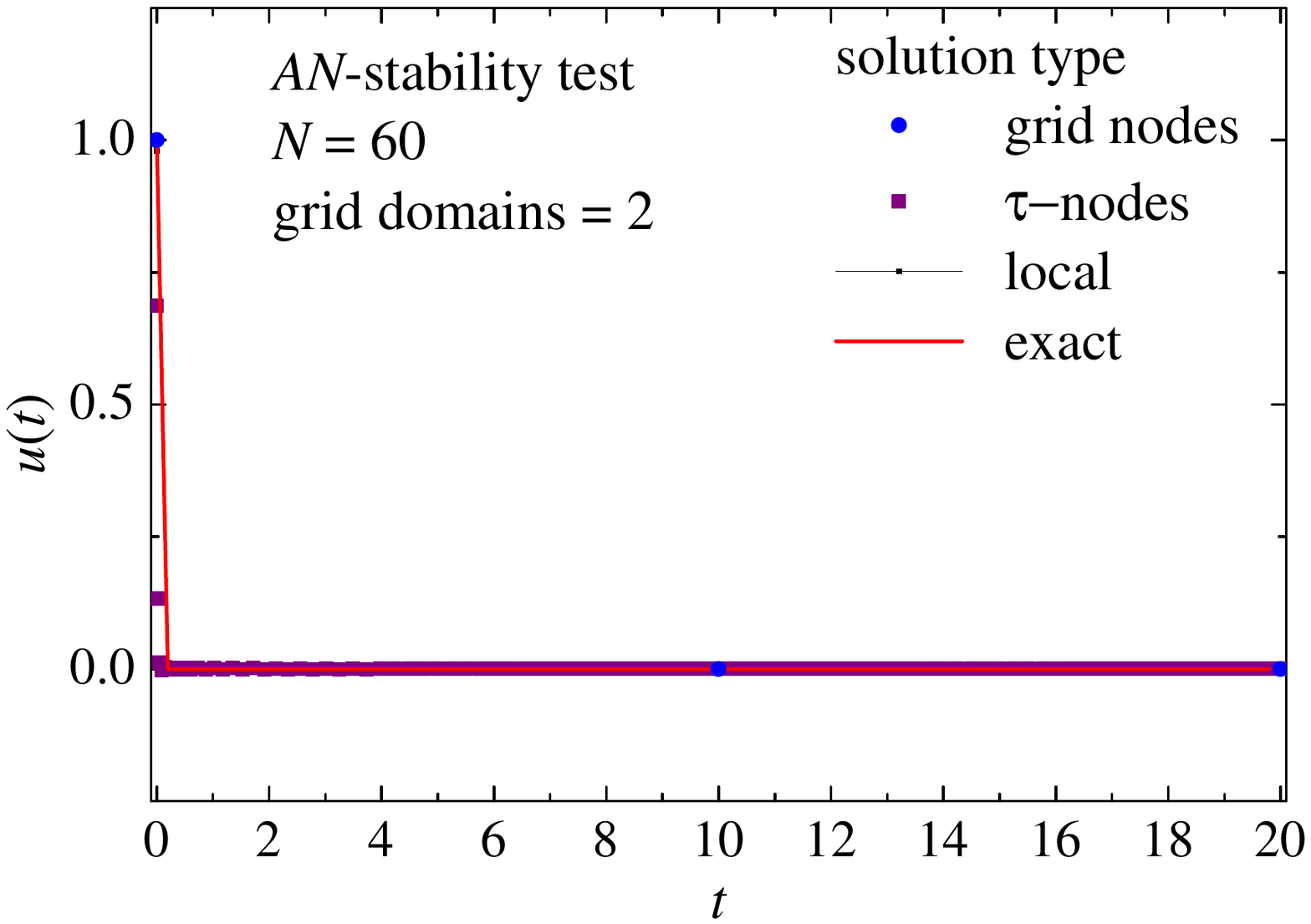}\\
\caption{\label{fig:astab_and_anstab}%
Numerical solution $u(t)$, obtained by the ADER-DG method with $N = 8$, $16$, $60$ on a grid with $2$ discretization domains for $A$-stability test (\ref{eq:a_test}) and $AN$-stability test (\ref{eq:an_test}).
}
\end{figure}
\begin{figure}[h!]
\centering
\includegraphics[width=0.4\textwidth]{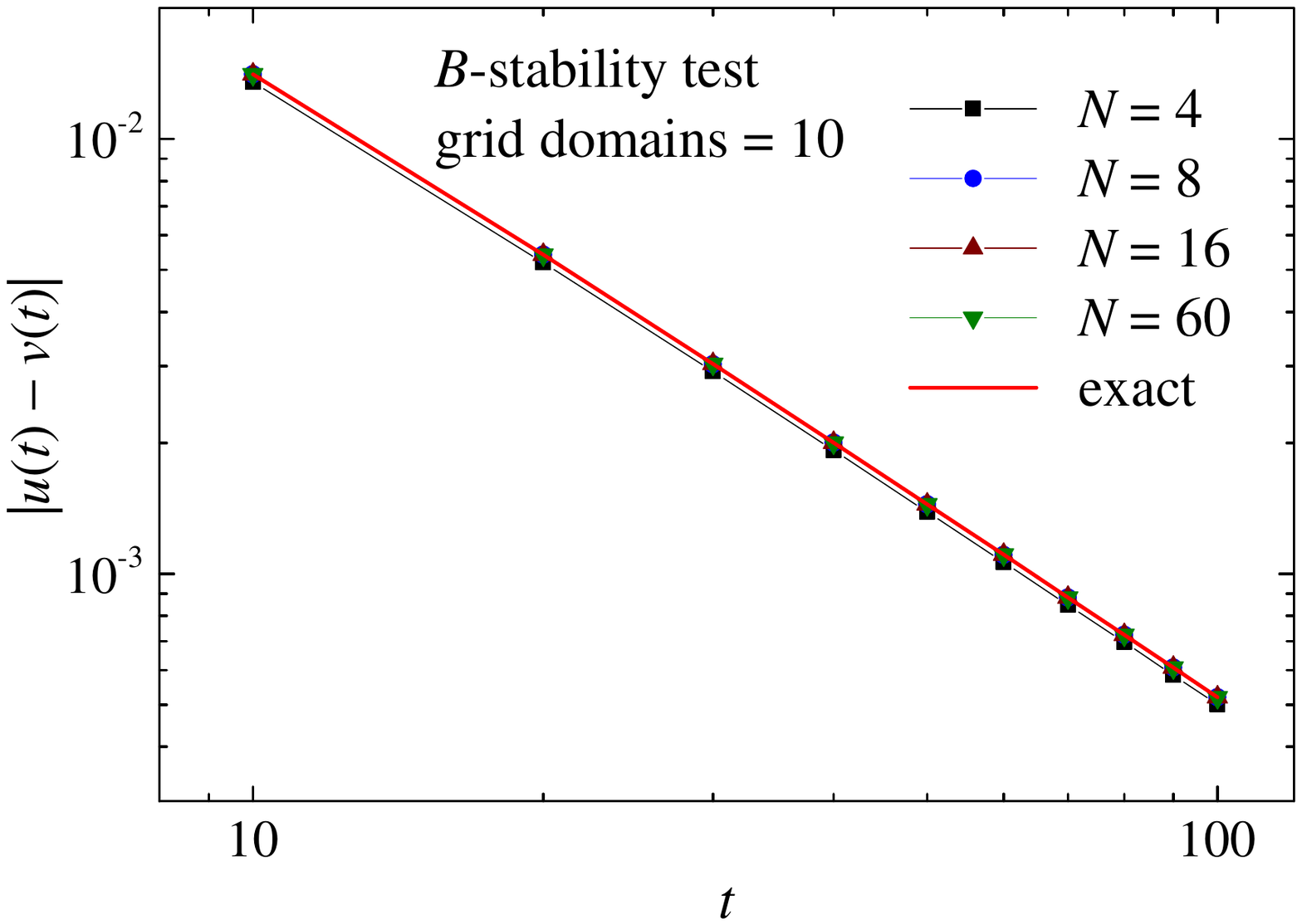}\hspace{10mm}
\includegraphics[width=0.4\textwidth]{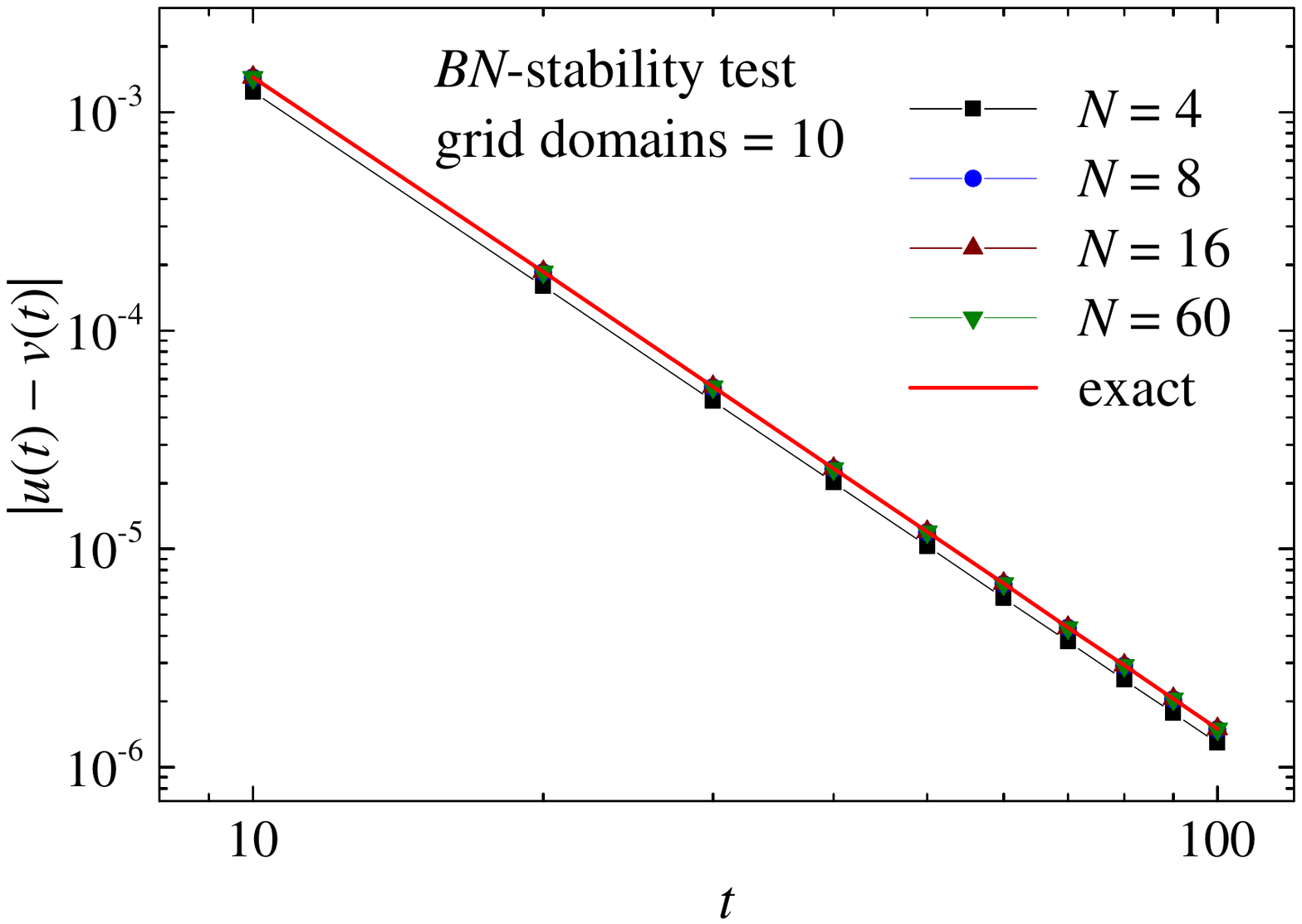}\\
\caption{\label{fig:bstab_and_bnstab}%
Norm of the difference $|u(t) - v(t)|$ between the solutions $u(t)$ and $v(t)$, obtained by the ADER-DG method with $N = 8$, $16$, $60$ for $B$-stability test (\ref{eq:b_test}) and $BN$-stability test (\ref{eq:bn_test}).
}
\end{figure}

\corrtext{The numerical solutions of problems (\ref{eq:a_test}) and (\ref{eq:an_test}) obtained by the ADER-DG method with $N = 4$, $8$, $16$, $60$ are showed in Fig.~\ref{fig:astab_and_anstab}. The numerical solutions at the nodes $\mathbf{u}_{n}$ demonstrate the monotonicity of the numerical solution, which demonstrates $A$-stability for problem (\ref{eq:a_test}) and $AN$-stability for problem (\ref{eq:an_test}). The effectiveness of applying the ADER-DG method to solving essentially stiff problems, which require at least $L$-stability, was demonstrated in~\cite{ader_dg_ivp_ode, ader_dg_ivp_dae}. Fig.~\ref{fig:bstab_and_bnstab} shows the norm of the difference $|u(t) - v(t)|$ of the numerical solutions $u(t)$ and $v(t)$ at the nodes of problems (\ref{eq:b_test}) and (\ref{eq:bn_test}) obtained by the ADER-DG method with $N = 4$, $8$, $16$, $60$. The obtained results clearly demonstrate the monotonic decay of the norm of the difference of the numerical solutions and well agreement with the exact analytical dependence even for a large discretization step ${\Delta t} = 10$. This demonstrates $B$-stability for problem (\ref{eq:b_test}) and $BN$-stability for problem (\ref{eq:bn_test}).}

\section{Conclusion}

\corrtext{In conclusion, it should be noted that in this paper a rigorous study of the linear and nonlinear stability properties of the ADER-DG method are clarified and proved. It is proved that the ADER-DG method is $AN$-stable, whereas previously only $A$-stability was known~\cite{ader_proofs_2025}. The $L$-stability property is emphasized, in contrast to the GL method. It is proved that the ADER-DG method has the property of nonlinear $BN$-stability, and as a consequence $B$-stable, and it is algebraically stable.}

\corrtext{Several other relations useful for an application and implementation of the ADER-DG method are proved. The relationship between the ADER-DG method and the RK method it generates with other RK methods, such as Radau methods, is revealed.}

\corrtext{Examples of the application of the ADER-DG method for solving the ODE system are presented, which demonstrated good agreement with the expected theoretical results. The empirical results~\cite{ader_dg_ivp_ode} are proven. Demonstration application examples are presented, demonstrating compliance with the expected convergence orders of the ADER-DG method. Demonstration examples of linear $A$- and $AN$-stability and nonlinear $B$- and $BN$-stability are presented.}

\corrtext{It is planned to construct a theory of the ADER-DG method as applied to PDE systems based on the developed theory.}

\section*{Acknowledgments}
The author would like to thank the anonymous reviewers for their remarks that helped to improve the quality and readability of this article.
The author would like to thank A.P.~Popova for help in correcting the English text.


\begin{thebibliography}{10}

\bibitem{Babuska_book_2001}
{\sc I.~Babu{\v{s}}ka and T.~Strouboulis}, {\em The Finite Element Method and
  Its Reliability, Numerical Mathematics and Scientific Computation}, Clarendon
  Press, Oxford, 2001.

\bibitem{dg_ivp_ode_1}
{\sc M.~Baccouch}, {\em Analysis of a posteriori error estimates of the
  discontinuous {Galerkin} method for nonlinear ordinary differential
  equations}, Appl. Numer. Math., 106 (2016), pp.~129--153.

\bibitem{dg_ivp_ode_2}
{\sc M.~Baccouch}, {\em A posteriori error estimates and adaptivity for the
  discontinuous {Galerkin} solutions of nonlinear second-order initial-value
  problems}, Appl. Numer. Math., 121 (2017), pp.~18--37.

\bibitem{dg_ivp_ode_3}
{\sc M.~Baccouch}, {\em Superconvergence of the discontinuous {Galerkin} method
  for nonlinear second-order initial-value problems for ordinary differential
  equations}, Appl. Numer. Math., 115 (2017), pp.~160--179.

\bibitem{dg_ivp_ode_4}
{\sc M.~Baccouch}, {\em Analysis of optimal superconvergence of the local
  discontinuous {Galerkin} method for nonlinear fourth-order boundary value
  problems}, Numerical Algorithms, 86 (2021), pp.~1615--1650.

\bibitem{dg_ivp_ode_5}
{\sc M.~Baccouch}, {\em The discontinuous {Galerkin} method for general
  nonlinear third-order ordinary differential equations}, Appl. Numer. Math.,
  162 (2021), pp.~331--350.

\bibitem{dg_ivp_ode_6}
{\sc M.~Baccouch}, {\em Superconvergence of an ultra-weak discontinuous
  {Galerkin} method for nonlinear second-order initial-value problems},
  International Journal of Computational Methods, 20(2) (2023), p.~2250042.

\bibitem{ivp_ode_taylor_series_2017}
{\sc A.~Baeza, S.~Boscarino, P.~Mulet, G.~Russo, and D.~Zorio}, {\em
  Approximate {Taylor} methods for {ODEs}}, Computers and Fluids, 159 (2017),
  pp.~156--166.

\bibitem{ader_dg_seiemic_underwater}
{\sc C.~Bassi, S.~Busto, and M.~Dumbser}, {\em High order {ADER-DG} schemes for
  the simulation of linear seismic waves induced by nonlinear dispersive
  free-surface water waves}, Appl. Numer. Math., 158 (2020), p.~236.

\bibitem{Butcher_book_2016}
{\sc J.~Butcher}, {\em Numerical Methods for Ordinary Differential Equations},
  Wiley, UK, 2016.

\bibitem{Cockburn_base_3}
{\sc B.~Cockburn, S.~Hou, , and C.-W. Shu}, {\em {{TVB} {Runge-Kutta} local
  projection discontinuous {Galerkin} finite element method for conservation
  laws. IV. The multidimensional case}}, Math. Comp., 54 (1990), pp.~545--581.

\bibitem{Cockburn_base_2}
{\sc B.~Cockburn, S.-Y. Lin, and C.-W. Shu}, {\em {{TVB} {Runge-Kutta} local
  projection discontinuous {Galerkin} finite element method for conservation
  laws III: One-dimensional systems}}, J. Comput. Phys., 84 (1989),
  pp.~90--113.

\bibitem{Cockburn_base_1}
{\sc B.~Cockburn and C.-W. Shu}, {\em {{TVB} {Runge-Kutta} local projection
  discontinuous {Galerkin} finite element method for conservation laws. II.
  General framework}}, Math. Comp., 52 (1989), pp.~411--435.

\bibitem{Cockburn_base_5}
{\sc B.~Cockburn and C.-W. Shu}, {\em {The {Runge-Kutta} local projection
  $P^1$-discontinuous-{Galerkin} finite element method for scalar conservation
  laws}}, ESAIM: M2AN, 25 (1991), pp.~337--361.

\bibitem{Cockburn_base_4}
{\sc B.~Cockburn and C.-W. Shu}, {\em {{TVB} {Runge-Kutta} local projection
  discontinuous {Galerkin} Method for Conservation Laws V: Multidimensional
  Systems}}, J. Comput. Phys., 141 (1998), pp.~199--224.

\bibitem{Dekker_Verwer_1984}
{\sc K.~Dekker and J.~Verwer}, {\em Stability of Runge-Kutta Methods for Stiff
  Nonlinear Differential Equations}, North-Holland, Amsterdam, 1984.

\bibitem{Delfour_1986}
{\sc M.~Delfour and F.~Dubeau}, {\em Discontinuous polynomial approximations in
  the theory of one-step, hybrid and multistep methods for nonlinear ordinary
  differential equations}, Math. Comp., 47 (1986), p.~169.

\bibitem{Delfour_1981}
{\sc M.~Delfour, W.~Hager, and F.~Trochu}, {\em Discontinuous {Galerkin}
  methods for ordinary differential equations}, Math. Comp., 36 (1981),
  pp.~455--473.

\bibitem{PNPM_DG_2010}
{\sc M.~Dumbser}, {\em Arbitrary high order ${P}_{N}{P}_{M}$ schemes on
  unstructured meshes for the compressible {Navier-Stokes} equations},
  Computers \& Fluids, 39 (2010), pp.~60--76.

\bibitem{ader_stiff_1}
{\sc M.~Dumbser, C.~Enaux, and E.~Toro}, {\em Finite volume schemes of very
  high order of accuracy for stiff hyperbolic balance laws}, J. Comput. Phys.,
  227 (2008), p.~3971.

\bibitem{PNPM_DG_2009}
{\sc M.~Dumbser and O.~Zanotti}, {\em Very high order ${P}_{N}{P}_{M}$ schemes
  on unstructured meshes for the resistive relativistic mhd equations}, J.
  Comput. Phys., 228 (2009), p.~6991.

\bibitem{ader_dg_gr_z4_2024}
{\sc M.~Dumbser, O.~Zanotti, E.~Gaburro, and I.~Peshkov}, {\em A well-balanced
  discontinuous {Galerkin} method for the first-order {Z4} formulation of the
  {Einstein-Euler} system}, J. Comp. Phys., 504 (2024), p.~112875.

\bibitem{ader_dg_wb_shwater_2022}
{\sc E.~Fernandez, M.~Diaz, M.~Dumbser, and T.~de~Luna}, {\em An arbitrary high
  order well-balanced {ADER-DG} numerical scheme for the multilayer
  shallow-water model with variable density}, J. Sci. Comput., 90 (2022),
  p.~52.

\bibitem{ader_dg_PNPM}
{\sc E.~Gaburro and M.~Dumbser}, {\em A posteriori subcell finite volume
  limiter for general {$P_N P_M$} schemes: Applications from gasdynamics to
  relativistic magnetohydrodynamics}, J. Sci. Comput., 86 (2021), p.~37.

\bibitem{dg_entropy}
{\sc E.~Gaburro, P.~Offner, M.~Ricchiuto, and D.~Torlo}, {\em High order
  entropy preserving {ADER-DG} schemes}, Applied Mathematics and Computation,
  440 (2023), p.~127644.

\bibitem{Hairer_book_1}
{\sc E.~Hairer, S.~N{\o}rsett, and G.~Wanner}, {\em Solving Ordinary
  Differential Equations I: Nonstiff Problems}, Springer-Verlag, Berlin,
  Heidelberg, 1987.

\bibitem{Hairer_book_2}
{\sc E.~Hairer and G.~Wanner}, {\em Solving Ordinary Differential Equations II:
  Stiff and Differential-Algebraic Problems}, Springer-Verlag, Berlin,
  Heidelberg, 1996.

\bibitem{ader_improving_2024}
{\sc M.~Han~Veiga, L.~Micalizzi, and D.~Torlo}, {\em On improving the
  efficiency of {ADER} methods}, Appl. Math. Comput., 466 (2024), p.~128426.

\bibitem{dec_vs_ader_2021}
{\sc M.~Han~Veiga, P.~{\"O}ffner, and D.~Torlo}, {\em {DeC} and {ADER}:
  Similarities, differences and a unified framework}, Journal of Scientific
  Computing, 87 (2021), p.~2.

\bibitem{ader_stiff_2}
{\sc A.~Hidalgo and M.~Dumbser}, {\em {ADER} schemes for nonlinear systems of
  stiff advection-diffusion-reaction equations}, J. Sci. Comput., 48 (2011),
  p.~173.

\bibitem{ivp_ode_taylor_series_soft_2005}
{\sc A.~Jorba and M.~Zou}, {\em A software package for the numerical
  integration of {ODEs} by means of high-order {Taylor} methods}, Experiment.
  Math., 14 (2005), pp.~99--117.

\bibitem{dg_entropy_add}
{\sc S.-C. Klein}, {\em Stabilizing discontinuous {Galerkin} methods using
  {Dafermos'} entropy rate criterion: I --- one-dimensional conservation laws},
  J. Sci. Comput., 95 (2023), p.~55.

\bibitem{dec_vs_ader_2023}
{\sc L.~Micalizzi, D.~Torlo, and W.~Boscheri}, {\em Efficient iterative
  arbitrary high-order methods: an adaptive bridge between low and high order},
  Commun. Appl. Math. Comput., 7 (2023), p.~40.

\bibitem{ader_proofs_2025}
{\sc P.~{\"O}ffner, L.~Petri, and D.~Torlo}, {\em Analysis for implicit and
  implicit-explicit {ADER} and {DeC} methods for ordinary differential
  equations, advection-diffusion and advection-dispersion equations}, Appl.
  Numer. Math., 212 (2025), pp.~110--134.

\bibitem{ader_stiff_3}
{\sc I.~Popov}, {\em Space-time adaptive {ADER-DG} finite element method with
  {LST-DG} predictor and a posteriori sub-cell {WENO} finite-volume limiting
  for simulation of non-stationary compressible multicomponent reactive flows},
  J. Sci. Comput., 95 (2023), p.~44.

\bibitem{ader_dg_ivp_ode}
{\sc I.~Popov}, {\em Arbitrary high order {ADER-DG} method with local {DG}
  predictor for solutions of initial value problems for systems of first-order
  ordinary differential equations}, J. Sci. Comput., 100 (2024), p.~22.

\bibitem{ader_stiff_4}
{\sc I.~Popov}, {\em Space-time adaptive {ADER-DG} finite element method with
  {LST-DG} predictor and a posteriori sub-cell {ADER-WENO} finite-volume
  limiting for multidimensional detonation waves simulation}, Computers \&
  Fluids, 284 (2024), p.~106425.

\bibitem{ader_eff_blas}
{\sc I.~Popov}, {\em The effective use of {BLAS} interface for implementation
  of finite-element {ADER-DG} and finite-volume {ADER-WENO} methods}, Commun.
  Comput. Phys., 38(5) (2025), pp.~1237--1330.

\bibitem{ader_dg_ivp_dae}
{\sc I.~Popov}, {\em High order {ADER-DG} method with local {DG} predictor for
  solutions of differential-algebraic systems of equations}, J. Sci. Comput.,
  102 (2025), p.~48.

\bibitem{lasl_rep_dg_1973}
{\sc W.~Reed and T.~Hill}, {\em Triangular mesh methods for the neutron
  transport equation}, Tech. Report LA-UR-73-479, Los Alamos Scientific
  Laboratory, 1973.

\bibitem{ader_init_1}
{\sc V.~Titarev and E.~Toro}, {\em {ADER: arbitrary high order Godunov
  approach}}, J. Sci. Comput., 17 (2002), p.~609.

\bibitem{ader_init_2}
{\sc V.~Titarev and E.~Toro}, {\em {ADER schemes for three-dimensional
  nonlinear hyperbolic systems}}, J. Comput. Phys., 204 (2005), p.~715.

\bibitem{ader_rev_2024}
{\sc E.~Toro, V.~Titarev, M.~Dumbser, A.~Iske, C.~Goetz, C.~Castro,
  G.~Montecinos, and D.~R.}, {\em The {ADER} approach for approximating
  hyperbolic equations to very high accuracy}, in Hyperbolic Problems: Theory,
  Numerics, Applications. Volume I. HYP 2022. {SEMA} {SIMAI} Springer Series,
  vol 34., C.~Pares, M.~Castro, T.~de~Luna, and M.~Muñoz-Ruiz, eds., Springer,
  Cham, 2024.

\bibitem{Wahlbin_lectures_1995}
{\sc L.~Wahlbin}, {\em Superconvergence in {Galerkin} Finite Element Methods},
  Springer-Verlag, Berlin Heidelberg, 1995.

\bibitem{ader_dg_seiemic}
{\sc S.~Wolf, M.~Galis, C.~Uphoff, A.-A. Gabriel, M.~P., D.~Gregor, and
  M.~Bader}, {\em An efficient {ADER-DG} local time stepping scheme for {3D}
  {HPC} simulation of seismic waves in poroelastic media}, J. Comp. Phys., 455
  (2022), p.~110886.

\bibitem{ader_dg_ideal_flows}
{\sc O.~Zanotti, F.~Fambri, M.~Dumbser, and A.~Hidalgo}, {\em Space-time
  adaptive {ADER} discontinuous {Galerkin} finite element schemes with a
  posteriori sub-cell finite volume limiting}, Computers \& Fluids, 118 (2015),
  p.~204.

\end{thebibliography}

\end{document}